\documentclass[11pt]{article}
\usepackage[]{amsmath,amssymb,amsfonts,latexsym,amsthm,enumerate,fullpage}
\usepackage[latin9]{inputenc}
\usepackage{amsmath}
\usepackage{amssymb}
\usepackage{breqn}
\usepackage{url}
\usepackage{graphicx}
\usepackage[pdfstartview=FitH]{hyperref}
\usepackage{amsthm}
\usepackage{hyperref}
\usepackage{url}
\usepackage{enumitem}
\usepackage{authblk}
\usepackage{bbm}
\usepackage{verbatim}
\usepackage{cleveref}
\usepackage{tikz}
\usepackage{todonotes}
\usepackage{orcidlink}

\newcommand{\im}{\operatorname{Im}}
\newcommand{\Ker}{\operatorname{Ker}}

\newcommand{\Aut}{\operatorname{Aut}}

\newcommand{\supp}{\operatorname{supp}}
\newcommand{\Cone}{\operatorname{Cone}}

\newcommand{\Rad}{\operatorname{Rad}}

\newcommand{\Sym}{\operatorname{Sym}}
\newcommand{\sign}{\operatorname{sign}}

\newcommand{\flag}{\operatorname{Flag}}
\newcommand{\orifl}{\operatorname{Orifl}}
\newcommand{\st}{\operatorname{st}}
\newcommand{\lk}{\operatorname{lk}}
\newcommand{\mat}{\operatorname{Mat}}

\newcommand{\N}{\mathbb{N}}
\newcommand{\R}{\mathbb{R}}
\newcommand{\Z}{\mathbb{Z}}
\newcommand{\EE}{\mathcal{E}}
\newcommand{\II}{\mathbbm{1}}
\newcommand{\FF}{\mathcal{F}}
\newcommand{\HH}{\mathcal{H}}
\renewcommand{\epsilon}{\varepsilon}
\newcommand*{\uak}{\mathcal{U}_A(K)}

\newcommand*{\CC}{\mathcal{CC}}

\newcommand*{\chev}{\operatorname{Chev}}

\newcommand{\Ao}{\mathring A}
\newcommand*{\cobound}{\operatorname{cb}}
\newcommand*{\cosys}{\operatorname{cs}}
\newcommand*{\SL}{\operatorname{SL}}

\renewcommand{\tilde}{\widetilde}

\newcommand{\cF}{\mathcal F}

\newcommand{\cW}{\mathcal W}

\newcommand{\El}{\ell}

\newcommand{\diagnode}[1]{\fill #1 circle (.1);}

\DeclareMathOperator*{\freeprod}{\raisebox{-2pt}{\scalebox{2}{$\ast$}}}

\makeatletter

\begin{document}
\newtheorem{theorem}{Theorem}[section]
\newtheorem{lemma}[theorem]{Lemma}
\newtheorem{definition}[theorem]{Definition}
\newtheorem{claim}[theorem]{Claim}

\newtheorem{conjecture}[theorem]{Conjecture}
\newtheorem{remark}[theorem]{Remark}
\newtheorem{proposition}[theorem]{Proposition}
\newtheorem{corollary}[theorem]{Corollary}
\newtheorem{observation}[theorem]{Observation}
\newtheorem{fact}[theorem]{Fact}
\theoremstyle{definition}
\newtheorem{example}[theorem]{Example}

\newcommand{\subscript}[2]{$#1 _ #2$}

\author{Izhar Oppenheim  \footnote{The author is partially supported by ISF grant no.  242/24.} \orcidlink{0000-0002-1849-1851}}

\affil{Department of Mathematics, Ben-Gurion University of the Negev, Be'er Sheva 84105, Israel, izharo@bgu.ac.il}

\author{Inga Valentiner-Branth \footnote{The author is supported by the FWO and the F.R.S.--FNRS under the Excellence of Science (EOS) program (project ID~40007542).} \footnote{The authors contributed equally to this work.} \orcidlink{0000-0001-7744-8733}
}
\affil{Department of Mathematics, Computer Science and Statistics, Ghent University, Krijgslaan 281 - S9, 9000 Ghent, Belgium, Inga.ValentinerBranth@ugent.be}

\title{New cosystolic high-dimensional expanders from KMS groups}

\date{September 9, 2026}

\maketitle

\begin{abstract}
Cosystolic expansion is a high-dimensional generalization of the Cheeger constant for simplicial complexes. Originally, this notion was motivated by the fact that it implies the topological overlapping property,  but more recently it was shown to be connected to problems in theoretical computer science such as  list agreement expansion and agreement expansion in the low soundness regime. 

There are only a few constructions of high-dimensional cosystolic expanders and,  in dimensions larger than $2$,  the only known constructions prior to our work were (co-dimension 1)-skeletons of quotients of affine buildings.   In this paper,  we give the first coset complex construction of cosystolic expanders for an arbitrary dimension.  Our construction is more symmetric and arguably more elementary than the previous constructions relying on quotients of affine buildings. 

The coset complexes we consider arise from finite quotients of Kac--Moody--Steinberg (KMS) groups, called KMS complexes.  KMS complexes were introduced in recent work by Grave de Peralta and Valentiner-Branth where it was shown that they are local spectral expanders.  Our result is that KMS complexes, satisfying some minor condition, give rise to infinite families of bounded degree cosystolic expanders of arbitrary dimension and for any finitely generated Abelian coefficient group.
	
This result is achieved by observing that proper links of KMS complexes are joins of opposite complexes in spherical buildings.  In order to show that these opposite complexes are coboundary expanders,  we develop a new method for constructing cone functions by iteratively adding sets of vertices. Using the prior local-to-global results, we obtain cosystolic expansion for the (co-dimension 1)-skeletons of the KMS complexes. 
\end{abstract}



\section{Introduction}

High-dimensional expanders (HDXs) are high-dimensional analogues of expander graphs. The basic idea behind HDXs is to extend notions of expansion from graphs to higher-dimensional objects, namely simplicial complexes. Two main notions of expansion are commonly used in the definition of HDXs: spectral expansion (more precisely, local spectral expansion; see the exact definition in \Cref{Local spectral expansion subsec}) and coboundary/cosystolic expansion, which are higher-dimensional analogues of the Cheeger constant (see the exact definition in \Cref{cosystolic/coboundary expansion subsec}). Unlike the case of expander graphs, however, these two notions are not equivalent.

In this paper, we focus on coboundary/cosystolic expansion rather than spectral expansion. The original motivation for studying cosystolic expansion came from the result of \cite{DKW}, showing that such complexes satisfy the topological overlapping property (see the definition in \cite[Introduction]{DKW}). Since then, coboundary/cosystolic expansion has found numerous applications, including coding theory and property testing \cite{KL, GK, BafM}, quantum codes \cite{EKZ, KT}, and Probabilistically Checkable Proofs (PCPs) \cite{DDL, BMVY}.



In our work below,  we construct new examples of bounded degree cosystolic HDXs.  Our construction relies on KMS complexes (which are coset complexes introduced in \cite{hdxfromkms}) and on a novel idea of constructing a cone function in an iterative process.

\subsection{KMS complexes}

In \cite{KOhdx18, KO2023high},  Kaufman and Oppenheim gave a construction of spectral HDXs based on the idea of coset complexes (see below).  The construction of Kaufman and Oppenheim was based on the group $\SL_{n+1} (\mathbb{F}_p [t])$,  and was later generalized by O'Donnell and Pratt \cite{ODP} to other Chevalley groups.  

In \cite{hdxfromkms},  another coset complex construction of spectral HDXs was given using Kac--Moody--Steinberg groups.  Namely,  in \cite{hdxfromkms} it was shown that given a generalized Cartan matrix satisfying some conditions and a finite field $\mathbb{F}_q$,  there is an infinite family of finite coset complexes that are spectral HDXs. The construction of KMS complexes involves some technicalities and we refer the reader to \Cref{The KMS complexes sec} for the full details and an explicit example. Our results are restricted to a certain class of KMS complexes characterized as follows. We call a generalized Cartan matrix $A=(A_{ij})_{i,j \in I}$ of rank $m$ $n$-classical (for $n \leq m$) if every submatrix $(A_{ij})_{i,j \in J}, J\subset I, \lvert J \rvert \leq n$ is of classical (i.e. type $A_k,B_k,C_k,D_k$) type or if it is reducible, all irreducible parts are of classical type. A KMS complex is called $n$-classical if its underlying generalized Cartan matrix is $n$-classical. 

While KMS complexes can be seen as a variant of the Chevalley group constructions of \cite{KO2023high,  ODP},  they seem to have a significant advantage that is relevant to our result below.  Namely,  the links of vertices of KMS complexes are associated to opposite complexes which are well-studied simplicial complexes.  Opposite complexes are large subcomplexes of spherical buildings and were studied independently to answer questions of finiteness properties of arithmetic groups (see for instance Abramenko \cite{Abr}).  The connection of opposite complexes to KMS complexes is that every link of a KMS complex is either an opposite complex or a join of opposite complexes.  Prior work on opposite complexes in \cite{Abr} suggests that,  in some sense, they mimic the properties of spherical buildings.  In our work below,  this point is explored and we prove that much like spherical buildings,  opposite complexes are coboundary expanders.

\subsection{Main results}

A link of a pure $n$-dimensional simplicial complex is called \emph{proper} if its dimension is at least $1$ and at most $n-1$ (i.e. it is the link of a face of dimension between $0$ and $n-2$).  As noted above,  our main result is that, under some restrictions,  proper links of the KMS complexes constructed in \cite{hdxfromkms} (see \Cref{The KMS complexes sec} below) are coboundary expanders:

\begin{theorem}
\label{cobound exp of links intro thm}
Let $n \geq 3$ be an integer.  There exists $\varepsilon >0$ such that for every prime power $q > 2^{2n-1}$ and every $n$-dimensional, $n$-classical KMS complex $X$ constructed over $k=\mathbb{F}_q$,  all the proper links of $X$ are $\varepsilon$-coboundary expanders \footnote{When we say that a simplicial complex $Y$ is an $\varepsilon$-coboundary expander,  we mean that for every $0 \leq j \leq \dim (Y) -1$,  its $j$-th coboundary expansion is bounded from below by $\varepsilon$ -  see exact definition in \Cref{cosystolic/coboundary expansion subsec}} (with respect to any finitely generated Abelian coefficient group).   
\end{theorem}

The theorem relies on the fact that proper links of $n$-classical KMS complexes are joins of opposite complexes in spherical buildings of classical type. Given an explicit $n$-classical KMS group with the extra restriction that no subdiagram of type $D_k$ appears in the Dynkin diagram, it is possible to track $\varepsilon$ and get an explicit (but very rough) bound. Tracking an explicit bound for $\varepsilon$ in the $D_n$ case could probably be done by very carefully analysing the proof, but we did not attempt to do so.

In \cite{hdxfromkms}, it was shown that $n$-dimensional KMS complexes over $\mathbb{F}_q$ are $\frac{2}{\sqrt{q}}$-local spectral expanders for every $q$ that is large enough with respect to $n$ (the KMS complexes are constructed under some assumptions on the corresponding generalized Cartan matrices, see exact formulation in Theorem \ref{KMS spectral thm} below).  Combining this fact with the local-to-global results of \cite{EK,  KM,  DD-cosys} yields the following theorem. 
\begin{theorem}
\label{cosys exp of KMS intro thm}
For an integer $n \geq 3$,  let $\lbrace X_i \rbrace_i$ be a family of $n$-dimensional, $n$-classical KMS complexes over $\mathbb{F}_q$ (see details in \Cref{The KMS complexes sec}) and let $\lbrace Y_i \rbrace_i$ be the family of $(n-1)$-dimensional skeletons of $\lbrace X_i \rbrace$.   There are $\varepsilon >0,  \mu >0, q_0 \in \mathbb{N}$ such that  if $q \geq q_0$ every $Y_i$ is an $(\varepsilon,  \mu)$-cosystolic expander \footnote{When we say that a simplicial complex $Y$ is a $(\varepsilon,  \mu)$-cosystolic expander,  we mean that for every $0 \leq j \leq \dim (Y) -1$,  its $j$-th cosystolic expansion is bounded from below by $\varepsilon$ and that the size $j$-cosystoles is bounded from below by $\mu$ - see exact definition in \Cref{cosystolic/coboundary expansion subsec}} (with respect to any finitely generated Abelian group).   
\end{theorem}

By \cite{DKW},  it follows that the family of complexes $\lbrace Y_i \rbrace_i$ defined in the Theorem above has the topological overlapping property (see exact formulation in \cite{DKW}).

\subsection{Our contributions}

There are only a few examples of cosystolic HDXs (with a uniformly bounded degree).  Indeed,  prior to our work,  the known constructions were:
\begin{itemize}
\item Constructions stemming from quotients of affine buildings (e.g.,  Ramanujan complexes \cite{LSV}).  These constructions give cosystolic HDXs of every dimension with respect to any finitely generated Abelian coefficient group via the work of Evra and Kaufman \cite{EK} and the generalizations of Kaufman and Mass \cite{KM} and of Dikstein and Dinur \cite{DD-cosys}.  However,  all known such constructions are less symmetric and less elementary than the construction arising from coset complexes (see below). 
\item Coset complexes arising from Chevalley groups. There are two constructions of this flavour and both only yield $2$-dimensional cosystolic HDXs.  First,  Kaufman and Oppenheim \cite{KO-cobound} showed that their coset complex construction of spectral HDXs gives rise to $2$-dimensional cosystolic HDXs with respect to any group coefficients.  Second,  in a more recent work,  O'Donnell and Singer \cite{ODS} showed that for coset complexes arising from $B_3$ Chevalley groups,  the $2$-skeletons of links of vertices are coboundary expanders with respect to $\mathbb{F}_2$ coefficients and that these complexes give rise to $2$-dimensional cosystolic HDXs with respect to $\mathbb{F}_2$ coefficients.
\end{itemize}

Our construction is the first construction of cosystolic HDXs arising from coset complexes that holds beyond dimension $2$.  Furthermore,  our results are more general than previous results on coset complexes arising from Chevalley groups,  since our results hold for finitely generated Abelian coefficient groups and for KMS complexes of the types $\widetilde{A}_n,  \widetilde{B}_n,  \widetilde{C}_n,  \widetilde{D}_n$. 

The method of our proof (as described below) is also novel in the setting of HDX -- we introduce a new method of constructing a cone function by iteratively adding sets of vertices to a subcomplex until we added all the vertices of the entire complex.

\subsection{Proof idea} 

As noted above,  \Cref{cosys exp of KMS intro thm} follows from \Cref{cobound exp of links intro thm} via the results of \cite{EK,  KM,  DD-cosys}.  Thus, we will focus on the idea of the proof of \Cref{cobound exp of links intro thm}, i.e., on the idea of proving that links of KMS complexes are coboundary expanders.  Since these links are symmetric,  it is enough to show that they have a cone function with bounded radius.  

In the constructions of cosystolic expanders stemming from quotients of affine buildings,  a cone function was constructed in \cite{LMMExpBuilding} by utilizing the fact that the links were spherical buildings that have a rich apartment structure.  In our construction,  the links of KMS complexes do not have such an apartment structure and thus we need a different method for constructing a cone function with bounded radius.  The method that we use is based on the idea of constructing a cone function by starting with a subcomplex with a cone function, iteratively adding sets of vertices and extending the cone function to that set.  

This idea stems from the works of \cite{Quill,  Vogt,  AApaper}, where this method was used to prove certain simplicial complexes are homotopy equivalent to a wedge of spheres.   In \cite{Abr},  Abramenko showed how to apply this idea to prove that certain subcomplexes of spherical classical buildings are homotopy equivalent to a wedge of spheres.  The subcomplexes treated in \cite{Abr} are exactly the links in our $n$-classical KMS complexes discussed above and we show that the same method can be ``quantified'' in order to produce a cone function with bounded radius for the links of our  $n$-classical KMS complexes.

The exact formulation of our result is quite technical (see Theorem \ref{adding vertices to cone thm} below), and we will explain the idea only on the level of vanishing of homology (as noted above,  a cone function can be thought of as a quantification of vanishing of homology with respect to a base point).   Let $X$ be an $n$-dimensional simplicial complex with a full subcomplex $X'$.   Given a vertex $w \in X$, we denote $X' \cup \lbrace w \rbrace$ to be the full subcomplex of $X$ spanned by the vertices of $X'$ and $w$.  The link of $w$ relative to $X'$ is the intersection $X_w \cap X'$,  where $X_w$ is the link of $w$ in $X$.  We will show that vanishing of reduced homology for $X'$ and the link of $w$ relative to $X'$ implies vanishing of reduced homology for $X' \cup \lbrace w \rbrace$. Indeed, let $A =X', B = \{w\} * (X_w \cap X')$ (more precisely the geometric realization of the two subcomplexes, but we will not consider this subtlety at this point). Then 
$$X' \cup \lbrace w \rbrace = X' \cup \lbrace \sigma \in X : w \in \sigma,  \sigma \setminus \lbrace w \rbrace \in X' \rbrace = A \cup B,$$
and $A \cap B = X_w \cap X'$ (note that the intersection of $A$ and $B$ is in the topological sense, which can be seen as $A\cap B = \{\tau \cap \sigma \mid \tau \in A, \sigma \in B\}$). 
Furthermore, $B$ can be contracted onto $\{w\}$. Hence we have $\widetilde{H}_i (A) = \widetilde{H}_i (B) =  0$ for $0 \leq i \leq n-1$ and $\widetilde{H}_i (A \cap B) = 0$ for $0 \leq i \leq n-2$. Thus, the Mayer-Vietoris Theorem implies that $\widetilde{H}_i (X)=0$ for $0 \leq i \leq n-1$. 

An important observation is that in the procedure above,  instead of adding a single vertex $w$ to $X'$,  one can add a set of vertices $W$ to $X'$ as long as no two vertices in $W$ are connected by an edge in $X$ and for every vertex $w \in W$, the link of $w$ relative to $X'$ has vanishing of reduced homology.   Theorem \ref{adding vertices to cone thm} is a quantification of this idea in the language of cone functions, where the cone radius is bounded as a function of the cone radius of $X'$ and the cone radii of the links of $w \in W$ relative to $X'$. 

As noted above,   applying this idea to links of KMS complexes follows the work of Abramenko \cite{Abr} on opposite complexes of spherical buildings of classical type.  Let $X$ be an $n$-dimensional, $n$-classical KMS complex constructed over $\mathbb{F}_q$. Under this assumption, the links of vertices in $X$ are subcomplexes of classical spherical buildings known as opposite complexes, or joins of such. Roughly speaking, Abramenko showed that if $q$ is sufficiently large,  then the links of $X$ can be constructed by the procedure of adding sets of vertices described above, that the number of steps in this procedure is independent of $q$ and that the (relative) links of the added vertices are joins of opposite complexes of spherical buildings of lower dimension.  This gives rise to an inductive bound on the cone radii and thus (by the fact that the links are symmetric) on the coboundary expansion of opposite complexes in classical spherical buildings (and thus of the proper links in $n$-classical KMS complexes).

\subsection{Organization}
The paper is organized as follows. In \Cref{Preliminaries sec},  we recall the notions of local spectral, coboundary and cosystolic high-dimensional expansion, explain the relevant (co)homological background, define cone functions and show some basic properties. This includes the fact that a cone function with coefficients over $\Z$ gives rise to cone functions with coefficients in arbitrary finitely generated Abelian groups, and how one can construct a cone function for the join of two simplicial complexes, given cone functions for each of them separately. 
In \Cref{Constructing a cone function by adding vertices idea}, we describe our technique of extending cone functions which is a key tool in our work.  In \Cref{The KMS complexes sec}, we recall the construction of the KMS complexes, and in \Cref{sec: describing buildings} we first describe the connection of KMS complexes to buildings and then explain a geometric way to construct classical buildings. 
Finally, in \Cref{sec: cbe for subcomplexes of spherical buildings} we combine the previous observations to prove \Cref{cobound exp of links intro thm}.

\section*{Acknowledgements}
The second author is grateful to Pierre-Emmanuel Caprace, Tom De Medts and Timoth\'ee Marquis for their helpful comments and advice. 
We are also grateful for the careful and detailed feedback received during the review process, which significantly improved both the presentation and the content of the paper.

\section{Preliminaries}
\label{Preliminaries sec}

\subsection{Weighted simplicial complexes}

Let $X$ be a finite $n$-dimensional simplicial complex. A simplicial complex $X$ is called {\em pure $n$-dimensional} if every face in $X$ is contained in some face of dimension $n$.  The set of all $k$-faces of $X$ is denoted $X(k)$, and we will be using the convention in which $X(-1) = \{\emptyset\}$.  For every $0 \leq k \leq n$,  the \textit{$k$-skeleton of $X$} is the simplicial complex $\bigcup_{i=-1}^k X(i)$.  We will say that $X$ is \textit{connected} if its $1$-skeleton is a connected graph.

For a complex $X$,  a subcomplex $X' \subseteq X$ is called a \textit{full subcomplex},  if for every $v_0,\dots,v_k \in X'$,  if $\lbrace v_0,\dots,v_k \rbrace \in X$, then $\lbrace v_0,\dots,v_k \rbrace \in X'$.  Equivalently,  $X'$ is a full subcomplex if there is a set $S$ that is a subset of the vertices of $X$ such that $X ' = \lbrace \sigma \in X : \sigma \subseteq S \rbrace$.

Given a pure $n$-dimensional simplicial complex $X$,  the \textit{weight function} $ w : \bigcup_{k=-1}^n X(k) \rightarrow \mathbb{R}_+$ is defined to be
$$\forall \tau \in X(k), w(\tau) = \frac{\vert \lbrace \sigma \in X(n) : \tau \subseteq \sigma \rbrace \vert}{{n+1 \choose k+1} \vert X(n) \vert}.$$

We note that $w$ is normalized such that for each $-1 \leq k \leq n$, the function $w$ can be thought of as a probability function on $X(k)$.

Given a finite pure $n$-dimensional complex $X$ and $\tau \in X$, \textit{the link of $\tau$} is the subcomplex $X_{\tau}$ defined as 
$$X_{\tau} = \lbrace \sigma \in X : \tau \cap \sigma = \emptyset,  \tau \cup \sigma \in X \rbrace.$$
To clarify the notation in some situations, we will also write $\lk_X(\tau)$ for the link of $\tau$ in $X$.
We note that if $\tau \in X(k)$,  then $X_\tau$ is a pure $(n-k-1)$-dimensional complex. We call $X_\tau$ a proper link, if $1 \leq \dim(X_\tau)\leq n-1$.

We denote the weight function $w_\tau$ to be the weight function on $X_\tau$ defined as above.
By slight abuse of notation, we will write $v \in X(0)$ to state that $v$ is a vertex of $X$.

\subsection{Local spectral expansion}
\label{Local spectral expansion subsec}

Let $X$ be a finite pure $n$-dimensional simplicial complex with $n \geq 1$ and $w$ the weight function on $X$ defined above.  We define the stochastic matrix of the random walk on (the $1$-skeleton of) $X$ to be the matrix indexed by $X(0) \times X(0)$ and defined as 
$$M (\lbrace v \rbrace,  \lbrace u \rbrace) = \begin{cases}
 \frac{w( \lbrace u,v \rbrace) }{\sum_{\lbrace u' ,v \rbrace \in X(1)} w( \lbrace u' ,v \rbrace)} & \lbrace v,u \rbrace \in X(1) \\
 0 & \lbrace v,u \rbrace \notin X(1)
\end{cases}.$$

For $\lambda <1$,  we say that $X$ is a \textit{(one-sided) $\lambda$-spectral expander} if $X$ is connected and the second largest eigenvalue of $M$ is $\leq \lambda$.  We say that $X$ is a  \textit{(one-sided) $\lambda$-local spectral expander} if for every $-1 \leq k \leq n-2$ and every $\tau \in X(k)$,  the simplicial complex $X_{\tau}$ is a (one-sided) $\lambda$-spectral expander.  We note that every (one-sided) $\lambda$-local spectral expander is a (one-sided) $\lambda$-spectral expander, since the link of the empty set is $X$ itself.

\subsection{Cohomological notations and cosystolic/coboundary expansion}
\label{cosystolic/coboundary expansion subsec}

Let $X$ be a pure $n$-dimensional simplicial complex and $\mathbb{A}$ an Abelian group.  

We denote $\overrightarrow{X} (k)$ to be the set of oriented $k$-simplices and for $\sigma \in \overrightarrow{X} (k)$,  we denote $- \sigma$ to be the simplex with reversed orientation from $\sigma$ on the same vertex set.  For $\lbrace v_0,\dots, v_k \rbrace \in X(k)$,  we denote $[v_0,\dots,v_k ] \in \overrightarrow{X} (k)$ to be the oriented simplex with the orientation induced from the ordering of the vertices.  Note that for every permutation $\gamma \in \Sym \lbrace 0,\dots, k\rbrace$ it holds that 
$$[v_{\gamma (0)},\dots,v_{\gamma (k)}] = (-1)^{\sign (\gamma)} [v_0,\dots,v_k].$$

We denote $C^k (X ; \mathbb{A})$ to be the set of $k$-cochains defined as follows: $C^k (X ;  \mathbb{A})$ is the set of functions $\phi : \overrightarrow{X} (k) \rightarrow \mathbb{A}$ that are anti-symmetric, i.e., for every $[v_0,\dots,v_k]$, $\phi (- [v_0,\dots,v_k]) = - \phi ([v_0,\dots,v_k])$.  We recall that $C^k (X ;  \mathbb{A})$ is an Abelian group with respect to pointwise addition.

Let $w$ be the weight function on $X$ defined above.  For $\phi \in C^k (X ; \mathbb{A})$, we define 
\begin{align*}
\supp(\phi) &= \lbrace \lbrace v_0,\dots,v_k \rbrace \in X(k) : \phi ([v_0,\dots,v_k]) \neq 0_{\mathbb{A}} \rbrace , \\
\Vert \phi \Vert &= \sum_{\sigma \in \supp(\phi)} w(\sigma).
\end{align*}

Given a subgroup $K \subseteq C^k (X ; \mathbb{A})$ and a cochain $\phi \in C^k (X ; \mathbb{A})$,  we define 
$$\Vert \phi - K \Vert = \min_{\psi \in K} \Vert \phi - \psi \Vert.$$

The \textit{coboundary map} $d_k : C^k (X ; \mathbb{A}) \rightarrow C^{k+1} (X ; \mathbb{A})$ is the map defined by mapping $\phi \in C^k(X;\mathbb{A})$ to $d_k \phi$ given by
$$(d_k \phi) ([v_0,\dots,v_{k+1}]) = \sum_{i=0}^{k+1} (-1)^i \phi ([v_0,\dots,\hat{v_i},\dots,v_{k+1}]),$$

where $[v_0,\dots,\hat{v_i},\dots,v_{k+1}]$ is the simplex obtained from $[v_0, \dots ,v_{k+1}]$ by removing the vertex $v_i$.
We denote $B^{k} (X ; \mathbb{A}) = \im (d_{k-1})$ and $Z^k  (X ; \mathbb{A}) = \Ker (d_k)$ and recall that these are Abelian subgroups of $C^k (X ; \mathbb{A})$ and that $B^{k} (X ; \mathbb{A}) \subseteq Z^k  (X ; \mathbb{A})$.  The $k$-cohomology $H^{k} (X; \mathbb{A})$ is defined as $H^{k} (X; \mathbb{A}) = Z^k  (X ; \mathbb{A}) / B^{k} (X ; \mathbb{A}) $. 

The $k$-coboundary/cosystolic expansion constant of $X$ (with coefficients in $\mathbb{A}$) are defined to be
$$h^k_{\cobound} (X ; \mathbb{A}) = \min_{\phi \in C^k (X ; \mathbb{A})\setminus B^k(X;\mathbb{A})} \frac{\Vert d_k \phi \Vert}{\Vert \phi - B^{k} (X ; \mathbb{A}) \Vert} ,$$
$$h^k_{\cosys} (X ; \mathbb{A}) = \min_{\phi \in C^k (X ; \mathbb{A})\setminus Z^k(X;\mathbb{A})} \frac{\Vert d_k \phi \Vert}{\Vert \phi - Z^{k} (X ; \mathbb{A}) \Vert} .$$

We note that by definition $h^k_{\cosys} (X ; \mathbb{A})  \geq h^k_{\cobound} (X ; \mathbb{A})$ and that by definition $h^k_{\cosys} (X ; \mathbb{A}) >0$.  We also note that $h^k_{\cosys} (X ; \mathbb{A})  = h^k_{\cobound} (X ; \mathbb{A})$ if and only if  $H^{k} (X; \mathbb{A}) = 0$ and that if $H^{k} (X; \mathbb{A}) \neq 0$,  then $h^k_{\cobound} (X ; \mathbb{A}) =0$.

\begin{remark}
We note that in some papers,  $h^k_{\cosys} (X ; \mathbb{A})$ is also called the \textit{$k$-th cocycle expansion} (e.g.,  see \cite{KKL}).  
\end{remark}

For a finite pure $n$-dimensional simplicial complex $X$ and a constant $\varepsilon>0$,  we will say that $X$ is a \textit{$(\varepsilon,  \mathbb{A})$-coboundary expander} if for every $0 \leq k \leq n-1$,  $h^k_{\cobound} (X ; \mathbb{A}) \geq \varepsilon$.  Also,  for constants $\varepsilon >0,  \mu >0$,  we will say that \textit{$X$ is a $(\varepsilon,  \mu,  \mathbb{A})$-cosystolic expander} if for every $0 \leq k \leq n-1$,  $h^k_{\cosys} (X ; \mathbb{A}) \geq \varepsilon$ and for every $\phi \in Z^k (X, \mathbb{A}) \setminus B^k (X, \mathbb{A})$,  $\Vert \phi \Vert \geq \mu$ (the parameter $\min \lbrace \Vert \phi \Vert : \phi \in Z^k (X, \mathbb{A}) \setminus B^k (X, \mathbb{A}) \rbrace$ is sometimes called the \textit{$k$-cosystole with respect to $\mathbb{A}$},  but we will not use this terminology).  

In \cite{KM, DD-cosys} the following local to global results were proven:
\begin{theorem}
Let $n \geq 3$ be an integer,  $\varepsilon ' >0$ a constant and $\mathbb{A}$ be a finitely generated Abelian group.  There are constants  $0< \lambda <1, \varepsilon >0,  \mu >0$ such that the following holds: For every finite pure $n$-dimensional simplicial complex $X$ if
\begin{itemize}
\item the simplicial complex $X$ is a (one-sided) $\lambda$-local spectral expander,
\end{itemize}
and
\begin{itemize}
\item for every $0 \leq k \leq n-2$ and every $\tau \in X (k)$,  it holds that $X_\tau$ is a $(\varepsilon ',  \mathbb{A})$-coboundary expander,
\end{itemize}
then the $(n-1)$-skeleton of $X$ is a $(\varepsilon,  \mu,  \mathbb{A})$-cosystolic expander. 
\end{theorem}

\subsection{Homological notations and cone functions}

Let $X$ be a pure $n$-dimensional simplicial complex and $\mathbb{A}$ an Abelian group.

Keeping the notations regarding oriented simplices above,  we denote $C_k (X ; \mathbb{A})$ to be the set of $k$-chains defined as follows: $C_k (X ;  \mathbb{A})$  is the set of all formal sums of the form 
$$\sum_{\sigma \in \overrightarrow{X} (k)} t_{\sigma}  \sigma,$$
such that $t_{\sigma} \in \mathbb{A}$ and $t_{- \sigma} = - t_{\sigma}$.   We note that $C_k (X ; \mathbb{A})$ is an Abelian group or equivalently a $\mathbb{Z}$-module.  For $A \in C_k (X ; \mathbb{A})$,  we define 
$$\supp (A) = \lbrace \lbrace v_0,\dots,v_k \rbrace : t_{[v_0,\dots,v_k]} \neq 0_{\mathbb{A}} \rbrace.$$ 

In the sequel, it will be useful to work with a ``basis'' for $C_k (X ;  \mathbb{A})$.  For any $\sigma_0 \in \overrightarrow{X} (k)$ and any $a \in \mathbb{A}$,  we define $a \mathbbm{1}_{\sigma_0} \in C_k (X ; \mathbb{Z})$ to be 
$$a \mathbbm{1}_{\sigma_0} =\sum_{\sigma \in \overrightarrow{X} (k)} t_{\sigma}^{\sigma_0}  \sigma$$
with 
$$t_\sigma^{\sigma_0} = \begin{cases}
a & \sigma = \sigma_0 \\
-a & \sigma = - \sigma_0 \\
0_{\mathbb{A}} & \text{otherwise}
\end{cases}.$$
With this notation,  every element in $C_k (X ; \mathbb{A})$ is a formal sum of elements of the form $a \mathbbm{1}_\sigma,  \sigma \in \overrightarrow{X} (k),  a \in \mathbb{A}$.

The \textit{boundary map} $\partial_k : C_k(X ; \mathbb{A}) \rightarrow C_{k-1}(X ; \mathbb{A})$ is the map given by
$$\partial_k \left( \sum_{[v_0,...,v_k] \in \overrightarrow{X} (k)} t_{[v_0,\dots,v_k]}  \mathbbm{1}_{[v_0,\dots,v_k]} \right)=  \sum_{[v_0,\dots,v_k] \in \overrightarrow{X} (k)}  \sum_{i=0}^{k} (-1)^i t_{[v_0,\dots,v_k]}  \mathbbm{1}_{[v_0,\dots,\hat{v_i},\dots,v_k]}.$$


\begin{definition}[Cone function]
Let $X$ be a finite simplicial complex.  Let $k \in \mathbb{N} \cup \lbrace 0, -1 \rbrace$  be a constant, $\mathbb{A}$ an Abelian group and $v$ be a vertex of $X$.  A $(k,  \mathbb{A})$-cone function with apex $v$ is a function $\Cone_{X} : \bigoplus_{j=-1}^{k} C_{j} (X ; \mathbb{A}) \rightarrow \bigoplus_{j= 0}^{k+1} C_{j} (X ; \mathbb{A})$ that fulfills the following conditions:
\begin{enumerate}
\item For every $a \in \mathbb{A}$,  $\Cone_{X} (a \mathbbm{1}_{\emptyset}) = a \mathbbm{1}_{[v]}$.
\item For every $-1 \leq j \leq k$,  $\Cone_X (C_{j} (X ; \mathbb{A})) \subseteq C_{j+1} (X ; \mathbb{A})$ and
$$\left.  \Cone_{X} \right\vert_{C_{j} (X ; \mathbb{A})} : C_{j} (X ; \mathbb{A}) \rightarrow C_{j+1} (X ; \mathbb{A})$$ 
is a homomorphism of $\Z$-modules.
\item For every $0 \leq j \leq k$ and every $A \in C_j (X ; \mathbb{A})$,
$$\partial_{j+1} \Cone_{X} (A) +  \Cone_{X} (\partial_j A) = A.$$
\end{enumerate}
Below,  we will refer to this equation as \textit{the cone equation}.
\end{definition}

\begin{definition}[Radius of a cone]
Let $X$ be an $n$-dimensional simplicial complex and $-1 \leq k \leq n-1$.  Given a $(k,  \mathbb{A})$-cone function $\Cone_{X}$,  we define for every $-1 \leq j \leq k$,  the $j$-th radius of $\Cone_{X}$ to be 
$$\Rad_j ( \Cone_{X} ) = \max_{\sigma \in \overrightarrow{X} (j),  a \in \mathbb{A}} \vert \supp (\Cone_X (a \mathbbm{1}_\sigma) )\vert.$$
Define the $(k,  \mathbb{A})$-th cone radius of $X$ to be
$$\Rad_k (X,  \mathbb{A}) = \min \lbrace \Rad_k ( \Cone_{X}  ) : \Cone_{X}  \text{ is a } (k, \mathbb{A}) \text{-cone function} \rbrace.$$
If no $(k,  \mathbb{A})$-cone function of $X$ exists, we define $\Rad_k (X,  \mathbb{A}) = \infty$.
\end{definition}

\begin{remark}
Note that in the case where $n=0$ and $X$ is just a non-empty set of vertices,  it holds by definition that $\Rad_{-1} (X,  \mathbb{A}) =1$.  
\end{remark}

\begin{example} \label{exmpl: cone for graph}
Let $\Gamma$ be a simplicial complex of dimension 1, i.e. a graph. Assume that $\Gamma$ is connected. Then we can construct a $(0,\Z)$-cone function for $\Gamma$ in the following way. 
Pick an arbitrary vertex $v \in \Gamma(0)$ and set $\Cone_\gamma(\II_\emptyset) = \II_{[v]}$.
Furthermore, we define $\Cone_\Gamma(\II_{[v]}) = 0$. 
In this case we have
$$\partial \Cone_\Gamma(\II_{[v]}) = 0 = \II_{[v]} - \Cone_\Gamma(\II_\emptyset) = \II_{[v]} - \Cone_\Gamma(\partial \II_{[v]})$$
and the cone equation is satisfied in this case. 

Next, let $w \in \Gamma(0), w \neq v$. Since $\Gamma$ is connected, there exists a path $w=w_0,\dots, w_n = v$ from $w$ to $v$. We define 
$$\Cone_\Gamma(\II_{[w]}) = \sum_{i=1}^n \II_{[w_{i-1},w_i]}.$$
Then 
\begin{align*}
\partial \Cone_\Gamma(\II_{[w]}) = \sum_{i=1}^n \II_{[w_{i-1}]} - \II_{[w_i]} = \II_{[w_0]} - \II_{[w_n]} = \II_{[w]} - \II_{[v]} = \II_{[w]} - \Cone_\Gamma(\II_\emptyset) = \II_{[w]} - \Cone_\Gamma(\partial \II_{[w]})
\end{align*}
which shows that the cone equation is satisfied in all cases, since we defined $\Cone_\Gamma$ on the generators of $C_0(\Gamma,\Z)$ and extend it to a homomorphism.

Additionally, note that if we choose the path from $w$ to $v$ to be of minimal length we have
$$\Rad_{-1}(\Cone_\Gamma) = 1, \qquad \Rad_0(\Cone_\Gamma) \leq \max_{w \in \Gamma(0)} \operatorname{dist}(v,w)$$
where $\operatorname{dist}(v,w)$ is the minimum length of a path between $v,w$.  We note that $ \Rad_0(\Cone_\Gamma)$ is upper-bounded by the diameter of $\Gamma$ in the usual graph theoretic sense.  Moreover,  when considering $\Rad_0 (\Gamma,  \mathbb{Z})$,  we minimize over all the choices of apex vertices $v$ and thus $\Rad_0 (\Gamma,  \mathbb{Z})$ is upper-bounded by the radius of $\Gamma$ in the usual graph theoretic sense.
\end{example}

We will show that a bound on the cone radii of $\mathbb{Z}$-cones bounds the cone radii for every finitely generated Abelian group.  This result is probably well known to experts, but we could not find a reference for it, so we include a sketch of the proof here (leaving some of the details for the reader).  

\begin{proposition}
\label{general abelian group cone prop}
Let $X$ be an $n$-dimensional simplicial complex and $-1 \leq k \leq n-1$.  For every finitely generated Abelian group $\mathbb{A}$,  it holds that $\Rad_k (X,  \mathbb{A}) \leq \Rad_k (X,  \mathbb{Z})$.
\end{proposition}

\begin{proof}
If $(k,  \mathbb{Z})$-cone functions of $X$ do not exist, then $\Rad_k (X,  \mathbb{Z}) = \infty$ and the inequality holds trivially.  Assume that there exists a $(k,  \mathbb{Z})$-cone function.  

By the fundamental theorem of finitely generated Abelian groups, there are prime powers $q_1,\dots,q_l$ and an integer $m \in \mathbb{N} \cup \lbrace 0 \rbrace$ such that $\mathbb{A} = \left( \bigoplus_{j=1}^l \mathbb{Z} / (q_j \mathbb{Z}) \right) \oplus \mathbb{Z}^m$.  We will denote the elements of $\mathbb{A}$ as $m+l$ tuples $(a_1,\dots,a_l, a_{l+1},\dots,a_{l+m})$,  where for every $1 \leq j \leq l$,  $a_j \in \mathbb{Z} / (q_j \mathbb{Z})$ and $a_{l+1},\dots,a_{l+m} \in \mathbb{Z}$.  We note that for every $-1 \leq i \leq n$,  
$$C_i (X, \mathbb{A}) = \left(\bigoplus_{j=1}^l C_i (X, \mathbb{Z} / (q_j \mathbb{Z})) \right) \oplus C_i (X, \mathbb{Z})^m.$$ 
Explicitly,  for every $\sigma \in \overrightarrow{X} (i)$  and every $(a_1,\dots,a_{l+m}) \in \mathbb{A}$, we identify $(a_1,\dots,a_{l+m}) \mathbbm{1}_{\sigma} \in C_i (X, \mathbb{A}) $ with $\bigoplus_{j=1}^{l+m} a_j \mathbbm{1}_\sigma \in \left(\bigoplus_{j=1}^l C_i (X, \mathbb{Z} / (q_j \mathbb{Z})) \right) \oplus C_i (X, \mathbb{Z})^m$.

For $1 \leq j \leq l$,  let $\varphi_j : \mathbb{Z} \rightarrow \mathbb{Z}/ (q_j \mathbb{Z})$ be the quotient map.  For every $-1 \leq i \leq n$,  $\varphi_j$ extends to a surjective homomorphism $\varphi_j : C_i (X, \mathbb{Z}) \rightarrow C_i (X, \mathbb{Z} / (q_j \mathbb{Z}))$ (this is a homomorphism of $\mathbb{Z}$-modules).  

Let $\Cone_X$ be a $(k,  \mathbb{Z})$-cone function.  For every $1 \leq j \leq l$,  we define the function
$$\Cone_X^{j } :  \bigoplus_{i=-1}^{k} C_{i} (X ; \mathbb{Z} / (q_j \mathbb{Z})) \rightarrow \bigoplus_{i= 0}^{k+1} C_{i} (X ; \mathbb{Z} / (q_j \mathbb{Z}))$$ 
as follows:
$$\Cone_X^{j} (\varphi_j (A)) = \varphi_j (\Cone_X (A)),  \forall -1 \leq i \leq k,  \forall A \in C_{i} (X ; \mathbb{Z}).$$

We note that $\Cone_X^{j}$ is well-defined: Let $A, A' \in  C_{i} (X ; \mathbb{Z})$ with $\varphi_j (A) = \varphi_j (A')$.  Then $A-A' \in \Ker (\varphi_j)$,  which implies that there is $A '' \in C_{i} (X ; \mathbb{Z})$ such that $A - A' = q_j A''$.  Thus
\begin{align*}
\Cone_X^{j} (\varphi_j (A)) -  \Cone_X^{j} (\varphi_j (A ')) = 
\varphi_j (\Cone_X (A)) - \varphi_j (\Cone_X (A ')) = \\
\varphi_j (\Cone_X (A -A ')) =
 \varphi_j (\Cone_{X}^{\mathbb{Z}}  (q_j A''))= \varphi_j (q_j \Cone_{X}^{\mathbb{Z}}  ( A'')) = 0_{\mathbb{Z} / (q_j \mathbb{Z})}
\end{align*}
and it follows that $\Cone_X^{j} (\varphi_j (A)) = \Cone_X^{j} (\varphi_j (A '))$.   Also,  by the fact that $\varphi_j$ is surjective on every $ C_i (X, \mathbb{Z} / (q_j \mathbb{Z}))$ it follows that $\Cone_X^{j}$ is defined on every $A \in C_{i} (X ; \mathbb{Z} / (q_j \mathbb{Z}))$.  

We leave it to the reader to verify that $\Cone_X^{j}$ is a $(k,   \mathbb{Z} / (q_j \mathbb{Z}))$-cone function and note that for every $i$, every $\sigma \in \overrightarrow{X} (i)$ and every $a \in \mathbb{Z} / (q_j \mathbb{Z})$ it holds that 
$$\supp (\Cone_X^j (a \mathbbm{1}_\sigma)) \subseteq \supp (\Cone_X ( \mathbbm{1}_\sigma)).$$

For every $l+1 \leq j \leq m+l$,  we denote $\Cone_X^{j} = \Cone_X$ (this is a $(k,\mathbb{Z})$-cone function for every $l+1 \leq j \leq l+m$).  For every $-1 \leq i \leq k$ and every 
$$\bigoplus_{j=1}^{l+m} A_j \in \left(\bigoplus_{j=1}^l C_i (X, \mathbb{Z} / (q_j \mathbb{Z})) \right) \oplus C_i (X, \mathbb{Z})^m,$$
we define
$$\Cone_X^{\mathbb{A}} (\bigoplus_{j=1}^{l+m} A_j) = \bigoplus_{j=1}^{l+m} \Cone_X^j (A_j).$$
Extending this map $\mathbb{Z}$-linearly yields a map
$$\Cone_X^{\mathbb{A}} : \bigoplus_{i=-1}^{k} \left( \left(\bigoplus_{j=1}^l C_i (X, \mathbb{Z} / (q_j \mathbb{Z})) \right) \oplus C_i (X, \mathbb{Z})^m \right) \rightarrow \bigoplus_{i= 0}^{k+1} \left(  \left(\bigoplus_{j=1}^l C_i (X, \mathbb{Z} / (q_j \mathbb{Z})) \right) \oplus C_i (X, \mathbb{Z})^m \right),$$
i.e., a map
$$\Cone_X^{\mathbb{A}} : \bigoplus_{i=-1}^{k} C_i (X, \mathbb{A}) \rightarrow \bigoplus_{i= 0}^{k+1} C_i (X, \mathbb{A}).$$

We conclude the proof by observing that $\Cone_X^{\mathbb{A}}$ is a $(k, \mathbb{A})$-cone function that for every $(a_1,\dots,a_{l+m}) \in \mathbb{A}$,  every $-1 \leq i \leq k$ and every $\sigma \in \overrightarrow{X} (i)$,  it holds that 
$$\supp (\Cone_X^{\mathbb{A}} ((a_1,\dots,a_{l+m}) \mathbbm{1}_\sigma) \subseteq \supp (\Cone_X (\mathbbm{1}_\sigma)).$$

Since $\Cone_X$ was an arbitrary $(k,\mathbb{Z})$-cone function,  it follows that $\Rad_k (X,  \mathbb{A}) \leq \Rad_k (X,  \mathbb{Z})$.

\end{proof}

In light of this Proposition,  in the sequel,  we will be interested in bounding $\mathbb{Z}$-radii for cone functions.  Below, we will omit the $\mathbb{Z}$ notation when discussing cone functions: we will write $k$-cone functions in lieu of $(k, \mathbb{Z})$-cone functions,  $\Rad_k (X)$ in lieu of $\Rad_k (X,  \mathbb{Z})$ etc.

\begin{observation}
\label{Cone function defined on simplices obs}
We observe that for every $j$,  $C_j (X ; \mathbb{Z})$ is a free $\mathbb{Z}$-module and thus the condition in the definition of the cone function that 
$$\left.  \Cone_{X} \right\vert_{C_{j} (X ; \mathbb{Z})} : C_{j} (X ; \mathbb{Z}) \rightarrow C_{j+1} (X ; \mathbb{Z})$$ 
is a homomorphism is equivalent to the condition that this map is a $\mathbb{Z}$-linear map of $\mathbb{Z}$-modules.  

By this linearity,  the cone equation is equivalent to the condition:
$$\partial_{j+1} \Cone_{X} (\mathbbm{1}_{\sigma}) +  \Cone_{X} (\partial_j \mathbbm{1}_{\sigma}) = \mathbbm{1}_{\sigma} , \forall \sigma \in \overrightarrow{X} (j).$$
In other words,  in order to verify that a linear function is a cone function,  it is enough to check the cone equation on $\mathbbm{1}_{\sigma}$ for every oriented simplex $\sigma$.  

Moreover,  when one defines a cone function,  it is enough to define it as a linear function that fulfills the cone equation on  $\mathbbm{1}_{\sigma}$ on the subset of oriented simplices that contain at least one orientation for every simplex.
\end{observation}

The work of Gromov \cite{Gromov} and subsequent works \cite{KO-cobound, KozM} showed that for symmetric simplicial complexes,  a bound on the cone radius gives rise to a bound on the coboundary expansion.  For a simplicial complex $X$,  we denote $\Aut (X)$ to be the simplicial automorphisms from $X$ to itself.  In \cite{KO-cobound} the following was proven:

\begin{theorem}
Let $n \geq 1$, $\mathcal{R} >0$ be constants and $\mathbb{A}$ an Abelian group.   For every finite pure $n$-dimensional simplicial complex,  if 
\begin{itemize}
\item the group $\Aut (X)$ acts transitively on $X (n)$,
\end{itemize}
and
\begin{itemize}
\item for every $0 \leq k \leq n-1$ it holds that $\Rad_k (X,  \mathbb{A})  \leq \mathcal{R}$,
\end{itemize}
then for every $0 \leq k \leq n-1$,  $h^k_{\cobound} (X,  \mathbb{A}) \geq \frac{1}{\mathcal{R} {n+1 \choose k+1}}$.
\end{theorem}

Combining this Theorem with Proposition \ref{general abelian group cone prop} yields the following:
\begin{theorem}
Let $n \geq 1$, $\mathcal{R} >0$ be constants.   For every finite pure $n$-dimensional simplicial complex,  if 
\begin{itemize}
\item The group $\Aut (X)$ acts transitively on $X (n)$
\end{itemize}
and
\begin{itemize}
\item For every $0 \leq k \leq n-1$ it holds that $\Rad_k (X,  \mathbb{Z})  \leq \mathcal{R}$
\end{itemize}
Then for every $0 \leq k \leq n-1$ and every finitely generated Abelian group $\mathbb{A}$ it holds that $$h^k_{\cobound} (X,  \mathbb{A}) \geq \frac{1}{\mathcal{R} {n+1 \choose k+1}}.$$
\end{theorem}

\begin{remark}
The statement in \cite{KO-cobound} is for $\mathbb{A} = \mathbb{F}_2$, but the proof given there generalizes almost verbatim for a general finitely generated Abelian group.
\end{remark}

\subsection{Cone radius of joins}
\label{Cone radius of joins subsec}

In this section,  we denote $\sqcup$ to be a disjoint union.  Given two non-empty finite simplicial complexes $Y_1$ and $Y_2$ with vertex sets $V(Y_1)$ and $V(Y_2)$,  the \textit{join} of $Y_1$ and $Y_2$,  denoted $Y_1 * Y_2$,  is the simplicial complex with the vertex set $V (Y_1) \sqcup V(Y_2)$ and simplices $\sigma \sqcup  \tau$ for every $\sigma \in Y_1$ and $\tau \in Y_2$.  We note that  $Y_1 * Y_2$ includes all the simplices of the form $\tau = \emptyset \sqcup \tau$ and $\sigma = \sigma \sqcup \emptyset$ for every $\sigma \in Y_1$ and $\tau \in Y_2$. 

Below,  we will need the following definition (which will also be used in the next section): Given a simplicial complex $X$ and a vertex $u$ in $X$,  we define the following operation;
For every $0 \leq j \leq k-1$ and every $[v_0,\dots,v_{j}] \in \overrightarrow{X_u} (j)$,  we define 
$$[u,  \mathbbm{1}_{[v_0,\dots,v_{j}]}] = \mathbbm{1}_{[u,  v_0,\dots,v_{j}]}.$$
Extending this definition linearly,  for every $A \in C_{j} (X_u ; \mathbb{Z})$,  we define $[u,A] \in  C_{j+1} (X ; \mathbb{Z})$.  We observe that for every  $[v_0,\dots,v_{j}] \in \overrightarrow{X_u} (j)$, 
$$\partial_{j+1} \mathbbm{1}_{[u,  v_0,\dots,v_{j}]} = \mathbbm{1}_{[v_0,\dots,v_{j}]} - [u,  \partial_{j} \mathbbm{1}_{[v_0,\dots,v_j]}]$$
and thus by linearity,  for every $A \in C_{j} (X_u ; \mathbb{Z})$,
\begin{equation}
\label{[w,A] eq}
\partial_{j+1} [u, A] = A - [u,  \partial_{j} A]. 
\end{equation}

\begin{proposition}
\label{cone of a join - single vertex prop}
Let $Y$ be a non-empty finite simplicial complex and $v$ a vertex that is not in $Y$. Then for every $k \in \mathbb{N} \cup \lbrace 0\rbrace$,  there is a $k$-cone function $\Cone_{\lbrace v \rbrace * Y}$ such that for every $-1 \leq j \leq k$, $\Rad_j (\Cone_{\lbrace v \rbrace * Y}) \leq 1$.  
\end{proposition}

\begin{proof}
Fix $k \in \mathbb{N} \cup \lbrace 0 \rbrace$ and define a function $\Cone_{\lbrace v \rbrace * Y}$ as follows: 
\begin{itemize}
\item $\Cone_{\lbrace v \rbrace * Y} (\emptyset) = \mathbbm{1}_{[v]}$.
\item For every $0 \leq j \leq k$ and every $[ v_0,\dots,v_j ] \in \overrightarrow{Y} (j)$,  
$$\Cone_{\lbrace v \rbrace * Y} (\mathbbm{1}_{[v_0,\dots, v_j]}) = [v,  \mathbbm{1}_{[v_0,\dots, v_j]}] = \mathbbm{1}_{[v, v_0,\dots, v_j]}.$$
\item For every $0 \leq j \leq k$ and every $\lbrace v_1,\dots,v_j \rbrace \in Y$,  
$$\Cone_{\lbrace v \rbrace * Y} (\mathbbm{1}_{[v,  v_1,\dots, v_j]}) = 0.$$
\end{itemize}
Extend the function $\Cone_{\lbrace v \rbrace * Y} $ linearly.

We check that the cone equation holds for the above function.  In the case of  $[ v_0,\dots,v_j ] \in \overrightarrow{Y} (j)$ this is due to \eqref{[w,A] eq}.  For $\lbrace v_1,\dots,v_j \rbrace \in Y (j-1)$,  we note that 
\begin{align*}
& \partial_{j+1} \Cone_{\lbrace v \rbrace * Y} (\mathbbm{1}_{[v,  v_1,\dots, v_j]}) + \Cone_{\lbrace v \rbrace * Y} (\partial_{j}  \mathbbm{1}_{[v,  v_1,\dots, v_j]}) = \\
& 0 + \Cone_{\lbrace v \rbrace * Y} (  \mathbbm{1}_{[v_1,\dots, v_j]}) -  \Cone_{\lbrace v \rbrace * Y} ( [v,  \partial_{j-1} [v_1,\dots,v_j]) = \\
& 0 + \mathbbm{1}_{[v,v_1,\dots, v_j]} + 0 = \mathbbm{1}_{[v,v_1,\dots, v_j]}
\end{align*}
as needed.  Thus, by Observation \ref{Cone function defined on simplices obs},  the above function is a $k$-cone function and it is easy to see that $\Rad_j (\Cone_{\lbrace v \rbrace * Y} ) \leq 1$ for every $-1 \leq j \leq k$.
\end{proof}

\begin{remark}
The reader should note that geometrically $\lbrace v \rbrace * Y$ is a cone (in the usual geometric definition) with apex $v$ and base $Y$.
\end{remark}

\begin{remark}
We note that in the above construction,  the cone function is well-defined for every $k$,  i.e., it is well-defined even in the cases where $k \geq \dim (\lbrace v \rbrace * Y)$.
\end{remark}

Given a simplicial complex $X$ with $ \lbrace v_0,\dots,v_{j_1} \rbrace,  \lbrace u_0,\dots,u_{j_2} \rbrace \in X$ such that 
 $ \lbrace v_0,\dots,v_{j_1} \rbrace \cap \lbrace u_0,\dots,u_{j_2} \rbrace  = \emptyset$ and $\lbrace v_0,\dots,v_{j_1} \rbrace \cup \lbrace u_0,\dots,u_{j_2} \rbrace   \in X$,  we define for $\sigma = [v_0,\dots,v_{j_1}]$ and $\tau =  [u_0,\dots,u_{j_2}]$, 
 $$[\sigma,  \tau ] = [v_0,\dots,v_{j_1}, u_0,\dots,u_{j_2}] \in \overrightarrow{X} (j_1+j_2+1)$$
 and
 $$[\mathbbm{1}_\sigma,  \mathbbm{1}_{\tau}] = \mathbbm{1}_{[\sigma,  \tau ] }.$$
For $X = Y_1 * Y_2$,  we note that for every $\sigma \in \overrightarrow{Y_1} (j_1)$ and $\tau \in \overrightarrow{Y_2} (j_2)$,  $[\mathbbm{1}_\sigma,  \mathbbm{1}_{\tau}]$ is well-defined.  Also note that for $\sigma$ and $\tau$ as above, 
$$\partial_{j_1+j_2+1} [\mathbbm{1}_\sigma,  \mathbbm{1}_{\tau}] = [\partial_{j_1} \mathbbm{1}_\sigma,  \mathbbm{1}_{\tau}] + (-1)^{j_1+1}  [\mathbbm{1}_\sigma,  \partial_{j_2} \mathbbm{1}_{\tau}].$$
Extending this equation linearly on  $C_{j_1} (Y_1 ; \mathbb{Z})$ and  $C_{j_2} (Y_2 ; \mathbb{Z})$ yields that for every $A_1 \in C_{j_1} (Y_1 ; \mathbb{Z})$ and every  $A_2 \in C_{j_2} (Y_2 ; \mathbb{Z})$,
 \begin{equation}
\label{[A1,A2] eq}
\partial_{j_1+j_2+1} [A_1,  A_2] = [\partial_{j_1} A_1,  A_2] + (-1)^{j_1+1}  [A_1,  \partial_{j_2} A_2].
\end{equation}
We note that this equation is actually a generalization of \eqref{[w,A] eq} above.

Furthermore, observe that $\supp([A_1,A_2]) = \supp(A_1) \times \supp(A_2)$.

The main idea of the following Proposition is taken from P. Wild's PhD Thesis \cite{Wild} and we claim no originality here.
\begin{proposition}
\label{cone of a join - general prop}
Let $Y_1,  Y_2$ be finite simplicial complexes of dimensions $n_1,  n_2$ respectively.   Assume that for $i=1,2$ there is an $n_i$-cone function $\Cone_{Y_i}$ and constants $R_j^{(i)} \in \mathbb{N},   -1 \leq j \leq n_i-1$ such that $\Rad_j (\Cone_{Y_i}) \leq R_j^i$.  Then there is an $(n_1+n_2)$-cone function $\Cone_{Y_1*Y_2}$ such that for every $-1 \leq j \leq n_1 + n_2$,  
$$\Rad_{j} (\Cone_{Y_1*Y_2}) \leq \begin{cases}
\max_{-1 \leq j_1 \leq j} R_{j_1}^{(1)} & j < n_1 \\
\max \lbrace \max_{-1 \leq j_1 \leq n_1-1} R_{j_1}^{(1)},  ((n_1 +1) R_{n_1-1}^{(1)} +1) R_{j-n_1-1}^{(2)} \rbrace & n_1 \leq j \leq n_1 +n_2  
\end{cases}.$$
\end{proposition}

\begin{proof}
Let $v$ be the apex of $\Cone_{Y_1}$.  Define 
$$\Cone_{Y_1 * Y_2} (\mathbbm{1}_{\emptyset}) = \mathbbm{1}_{[v]} .$$  
Let $-1 \leq j_1 \leq n_1$ and $-1 \leq j_2 \leq n_2$ such that $0 \leq j_1 + j_2+1 \leq n_1 + n_2$.  For $\sigma_i \in \overrightarrow{X} (j_i),  i=1,2$,  we define 
\begin{equation}
\label{Cone_(Y1*Y_2) eq1}
\Cone_{Y_1 * Y_2} (\mathbbm{1}_{[\sigma_1,  \sigma_2]}) = 
\begin{cases}
[\Cone_{Y_1} (\mathbbm{1}_{\sigma_1}),  \mathbbm{1}_{\sigma_2}] & \dim (\sigma_1) < n_1 \\
(-1)^{n_1+1} [\mathbbm{1}_{\sigma_1} -  \Cone_{Y_1} (\partial_{n_1} \mathbbm{1}_{\sigma_1}),  \Cone_{Y_2} (\mathbbm{1}_{\sigma_2})] & \dim (\sigma_1) = n_1
\end{cases}
\end{equation}
and extend it linearly.

We note that by linearity, 
\begin{equation}
\label{Cone_(Y1*Y_2) eq2}
\Cone_{Y_1 * Y_2} (\partial_{j_1+j_2+1} \mathbbm{1}_{[\sigma_1,  \sigma_2]}) =^{\eqref{[A1,A2] eq}} 
 \Cone_{Y_1 * Y_2} ( [\partial_{j_1} \mathbbm{1}_{\sigma_1},  \mathbbm{1}_{\sigma_2}]) + (-1)^{j_1+1} \Cone_{Y_1 * Y_2} (  [\mathbbm{1}_{\sigma_1},  \partial_{j_2} \mathbbm{1}_{\sigma_2}]).
\end{equation}

We will check that the cone equation is fulfilled by dealing with two separate cases. 

Assume first that $ \dim (\sigma_1) < n_1$.  Then
\begin{align*}
& \partial_{j_1+j_2+2} \Cone_{Y_1 * Y_2} (\mathbbm{1}_{[\sigma_1,  \sigma_2]})  =^{\eqref{Cone_(Y1*Y_2) eq1}} 
\partial_{j_1+j_2+2}  [\Cone_{Y_1} (\mathbbm{1}_{\sigma_1}),  \mathbbm{1}_{\sigma_2}] =^{\eqref{[A1,A2] eq}} \\
& [\partial_{j_1 +1} \Cone_{Y_1} (\mathbbm{1}_{\sigma_1}),  \mathbbm{1}_{\sigma_2}] + (-1)^{j_1+2} [ \Cone_{Y_1} (\mathbbm{1}_{\sigma_1}), \partial_{j_2}  \mathbbm{1}_{\sigma_2}] = \\
& [\mathbbm{1}_{\sigma_1} - \Cone_{Y_1} (\partial_{j_1} \mathbbm{1}_{\sigma_1}),  \mathbbm{1}_{\sigma_2}] + (-1)^{j_1+2} [ \Cone_{Y_1} (\mathbbm{1}_{\sigma_1}), \partial_{j_2}  \mathbbm{1}_{\sigma_2}] = \\
& [\mathbbm{1}_{\sigma_1},  \mathbbm{1}_{\sigma_2}] - \left( [\Cone_{Y_1} (\partial_{j_1} \mathbbm{1}_{\sigma_1}),  \mathbbm{1}_{\sigma_2}] + (-1)^{j_1+1} [ \Cone_{Y_1} (\mathbbm{1}_{\sigma_1}), \partial_{j_2}  \mathbbm{1}_{\sigma_2}] \right) =^{\eqref{Cone_(Y1*Y_2) eq2}} \\
& \mathbbm{1}_{[\sigma_1, \sigma_2]} - \Cone_{Y_1 * Y_2} (\partial_{j_1+j_2+1} \mathbbm{1}_{[\sigma_1, \sigma_2]}),
\end{align*}
thus the cone equation holds as needed.

Next,  assume that $ \dim (\sigma_1) = n_1$.  Then by \eqref{Cone_(Y1*Y_2) eq2},  
\begin{align*}
& \Cone_{Y_1 * Y_2} (\partial_{n_1+j_2+1} \mathbbm{1}_{[\sigma_1,  \sigma_2]}) =^{\eqref{Cone_(Y1*Y_2) eq2}} \\
 & \Cone_{Y_1 * Y_2} ( [\partial_{n_1} \mathbbm{1}_{\sigma_1},  \mathbbm{1}_{\sigma_2}]) + (-1)^{n_1+1} \Cone_{Y_1 * Y_2} (  [\mathbbm{1}_{\sigma_1},  \partial_{j_2} \mathbbm{1}_{\sigma_2}]) =^{\eqref{Cone_(Y1*Y_2) eq1}} \\
& [\Cone_{Y_1} (\partial_{n_1} \mathbbm{1}_{\sigma_1}),  \mathbbm{1}_{\sigma_2}] + (-1)^{n_1+1} (-1)^{n_1+1} [\mathbbm{1}_{\sigma_1} -  \Cone_{Y_1} (\partial_{n_1} \mathbbm{1}_{\sigma_1}),  \Cone_{Y_2} (\partial_{j_2}  \mathbbm{1}_{\sigma_2})]  = \\
& [\Cone_{Y_1} (\partial_{n_1} \mathbbm{1}_{\sigma_1}),  \mathbbm{1}_{\sigma_2}] + [\mathbbm{1}_{\sigma_1} -  \Cone_{Y_1} (\partial_{n_1} \mathbbm{1}_{\sigma_1}),  \Cone_{Y_2} (\partial_{j_2}  \mathbbm{1}_{\sigma_2})], 
\end{align*}
i.e., 
\begin{equation}
\label{Cone_(Y1*Y_2) eq3}
\Cone_{Y_1 * Y_2} (\partial_{n_1+j_2+1} \mathbbm{1}_{[\sigma_1,  \sigma_2]}) = [\Cone_{Y_1} (\partial_{n_1} \mathbbm{1}_{\sigma_1}),  \mathbbm{1}_{\sigma_2}] + [\mathbbm{1}_{\sigma_1} -  \Cone_{Y_1} (\partial_{n_1} \mathbbm{1}_{\sigma_1}),  \Cone_{Y_2} (\partial_{j_2}  \mathbbm{1}_{\sigma_2})].
\end{equation}

We also note that 
\begin{equation}
\label{fake partial^2 eq}
\partial_{n_1} \left( \mathbbm{1}_{\sigma_1} -  \Cone_{Y_1} (\partial_{n_1} \mathbbm{1}_{\sigma_1}) \right)) = \partial_{n_1} \mathbbm{1}_{\sigma_1} -  \left(\partial_{n_1} \mathbbm{1}_{\sigma_1} -  \Cone_{Y_1} (\partial_{n_1-1} \partial_{n_1} \mathbbm{1}_{\sigma_1}) \right) =0.
\end{equation}

Then
\begin{align*}
& \partial_{n_1+j_2+2} \Cone_{Y_1 * Y_2} (\mathbbm{1}_{[\sigma_1,  \sigma_2]})  =^{\eqref{Cone_(Y1*Y_2) eq1}} \\
& \partial_{n_1+j_2+2} \left( (-1)^{n_1+1} [\mathbbm{1}_{\sigma_1} -  \Cone_{Y_1} (\partial_{n_1} \mathbbm{1}_{\sigma_1}),  \Cone_{Y_2} (\mathbbm{1}_{\sigma_2})] \right) =^{\eqref{[A1,A2] eq}} \\
& (-1)^{n_1+1} [\partial_{n_1} \left( \mathbbm{1}_{\sigma_1} -  \Cone_{Y_1} (\partial_{n_1} \mathbbm{1}_{\sigma_1}) \right),  \Cone_{Y_2} (\mathbbm{1}_{\sigma_2})] + \\
& (-1)^{n_1+1} (-1)^{n_1+1} [ \mathbbm{1}_{\sigma_1} -  \Cone_{Y_1} (\partial_{n_1} \mathbbm{1}_{\sigma_1}) ,  \partial_{j_2+1} \Cone_{Y_2} (\mathbbm{1}_{\sigma_2})] =^{\eqref{fake partial^2 eq}} \\
& [ \mathbbm{1}_{\sigma_1} -  \Cone_{Y_1} (\partial_{n_1} \mathbbm{1}_{\sigma_1}),  \partial_{j_2+1} \Cone_{Y_2} (\mathbbm{1}_{\sigma_2})] =^{\text{the cone equation for } \Cone_{Y_2}} \\
& [ \mathbbm{1}_{\sigma_1} -  \Cone_{Y_1} (\partial_{n_1} \mathbbm{1}_{\sigma_1}),  \mathbbm{1}_{\sigma_2}- \Cone_{Y_2} ( \partial_{j_2} \mathbbm{1}_{\sigma_2})] = \\
&  [ \mathbbm{1}_{\sigma_1} -  \Cone_{Y_1} (\partial_{n_1} \mathbbm{1}_{\sigma_1}),  \mathbbm{1}_{\sigma_2}] -  [ \mathbbm{1}_{\sigma_1} -  \Cone_{Y_1} (\partial_{n_1} \mathbbm{1}_{\sigma_1}), \Cone_{Y_2} ( \partial_{j_2} \mathbbm{1}_{\sigma_2})]  = \\
& [\mathbbm{1}_{\sigma_1},  \mathbbm{1}_{\sigma_2}]  - \left( [\Cone_{Y_1} (\partial_{n_1} \mathbbm{1}_{ \sigma_1}),  \mathbbm{1}_{\sigma_2}] + [ \mathbbm{1}_{\sigma_1} -  \Cone_{Y_1} (\partial_{n_1} \mathbbm{1}_{\sigma_1}), \Cone_{Y_2} ( \partial_{j_2} \mathbbm{1}_{\sigma_2})]  \right) =^{\eqref{Cone_(Y1*Y_2) eq3}} \\
& \mathbbm{1}_{[\sigma_1, \sigma_2]} - \Cone_{Y_1 * Y_2} (\partial_{j_1+j_2+1} \mathbbm{1}_{[\sigma_1, \sigma_2]}),
\end{align*}
thus the cone equation holds as needed.

Let $0 \leq j \leq n_1+n_2$.  Then every $\sigma \in \overrightarrow{X} (j)$ is of the form $\sigma = [\sigma_1, \sigma_2]$ where $\sigma_i \in \overrightarrow{X} (j_i),  i=1,2$ and $j_1 + j_2 +1 = j$.  Note that $j_2 \geq -1$ and thus $j_1 \leq j$.   In particular,  if $j < n_1$, then $j_1 < n_1$. 

For such $j_1, j_2$,  it holds that if $j_1 < n_1$,  then 
$$\vert \supp (\Cone_{Y_1 * Y_2} (\mathbbm{1}_{[\sigma_1,  \sigma_2]})) \vert \leq R_{j_1}^{(1)}.$$

If $j_1 = n_1$,  then $j_2 = j - n_1 -1$ and 
\begin{align*}
\vert \supp (\Cone_{Y_1 * Y_2} (\mathbbm{1}_{[\sigma_1,  \sigma_2]})) \vert \leq 
\vert \supp ( \mathbbm{1}_{\sigma_1} -  \Cone_{Y_1} (\partial_{n_1} \mathbbm{1}_{\sigma_1})) \vert R_{j_2}^{(2)} \leq 
((n_1 +1) R_{n_1-1}^{(1)} +1) R_{ j - n_1 -1}^{(2)}.
\end{align*}

Combining these bounds yields the bound stated above for $\Rad_j (\Cone_{Y_1 * Y_2} )$.  
\end{proof}

\section{Constructing a cone function by adding vertices}
\label{Constructing a cone function by adding vertices idea}

The idea is to construct a cone function for a given complex $X$, by starting with a given subcomplex for which we have a cone function and adding a set of vertices which is ``well-behaved''.  The exact formalism of this idea is given in Theorem \ref{adding vertices to cone thm} that can be thought of as a special version of the Mayer-Vietoris theorem for cone functions.

\begin{theorem}
\label{adding vertices to cone thm}
Let $X$ be an $n$-dimensional simplicial complex and $X '$ be a full subcomplex.  Assume there is an $(n-1)$-cone function $\Cone_{X'}$ and  constants $R_k ' \in \mathbb{N},  -1 \leq j \leq n-1$ such that for every $-1 \leq j \leq n-1$,  $\Rad_j (\Cone_{X'})  \leq R_j '$. 
Also let $W \subseteq X(0)$ be a set of vertices such that the following holds:
\begin{enumerate}
\item $W \cap X' = \emptyset$. \label{add vert thm: as 1}
\item For every $w_1,  w_2 \in W$,  $\lbrace w_1, w_2 \rbrace \notin X (1)$.  \label{add vert thm: as 2}
\item For every $w \in W$,  $X_w \cap X'$ is a non-empty simplicial complex, i.e.,  for every $w \in W$ there exists  $v \in X' (0)$ such that $\lbrace w,v \rbrace \in X(1)$. \label{add vert thm: as 3}
\item There are $(n-2)$-cone functions $\Cone_{X_w \cap X'},  w \in W$ and constants $R_j'' \in \mathbb{N},  -1 \leq k \leq n-2$ such that for every $w \in W$ and every $-1 \leq j \leq n-2$,  $\Rad_j  (\Cone_{X_w \cap X'})  \leq R_j ''.$ \label{add vert thm: as 4}
\end{enumerate} 
Let $X' \cup W$ be the full subcomplex of $X$ spanned by $X' (0) \cup W$.  Then there is an $(n-1)$-cone function $\Cone_{X' \cup W}$ such that for every $-1 \leq j \leq n-1$,  $\Rad_j (\Cone_{X' \cup W}) \leq R_{j-1} '' (R_j '+1)$.
\end{theorem}

\begin{proof}
We will define $\Cone_{X' \cup W}$ by ``adding'' the cone functions of $X_w \cap X',  w \in W$ to $\Cone_{X'}$.  

We note that by linearity and Observation \ref{Cone function defined on simplices obs},  it is enough to define the cone function on every $\mathbbm{1}_{[v_0,\dots,v_j]}$ for every $0 \leq j \leq n-1$ and every $[v_0,\dots,v_j] \in \overrightarrow{X' \cup W} (j)$.  We also note that by our assumptions,  every $[v_0,\dots,v_j] \in \overrightarrow{X' \cup W} (j)$ is either in  $\overrightarrow{X'} (j)$ or that it has a single vertex in $W$.  Thus,  it is enough to define the cone function on $\mathbbm{1}_{[v_0,\dots,v_j]}$ in the following two cases: either $[v_0,\dots,v_j] \in \overrightarrow{X'} (j)$ or $v_0 = w \in W$ and $[ v_1,\dots,v_j ] \in \overrightarrow{X_w \cap X'} (j-1)$.

After the preceding discussion,  we define an $(n-1)$-cone function as a linear continuation of the following: First,  we define $\Cone_{X' \cup W}  (\emptyset) = \mathbbm{1}_{[v]}$ where $v$ is the apex of the cone function $\Cone_{X'}$.  Second, for every $0 \leq j \leq n-1$ and every $ [v_0,\dots,v_j] \in \overrightarrow{X'} (j)$,  we define
$$\Cone_{X' \cup W}  (\mathbbm{1}_{[v_0,\dots,v_j]}) = \Cone_{X'}  (\mathbbm{1}_{[v_0,\dots,v_j]}).$$
Last,  for every $0 \leq j \leq n-1$,  every $w \in W$ and every $[v_1,\dots,v_j] \in \overrightarrow{X_{w} \cap X'} (j-1)$,  we define
\begin{align*}
\Cone_{X' \cup W}  (\mathbbm{1}_{[w, v_1,\dots,v_j]}) = \Cone_{X'}  (\Cone_{X_w \cap X'} (\mathbbm{1}_{[v_1,\dots,v_j]})) - [w,  \Cone_{X_w \cap X'}  (\mathbbm{1}_{[v_1,\dots,v_j]})].
\end{align*}

We need to verify that the cone equation holds for every $\mathbbm{1}_{[w, v_1,\dots,v_j]}$.

We will start by computing $ \Cone_{X' \cup W}  (\partial_j \mathbbm{1}_{[w, v_1,\dots,v_j]}) $:
\begin{align*}
& \Cone_{X' \cup W}  (\partial_j \mathbbm{1}_{[w, v_1,\dots,v_j]})  =^{\eqref{[w,A] eq}} \\
&  \Cone_{X' \cup W}  (\mathbbm{1}_{[v_1,\dots,v_j]}) -   \Cone_{X' \cup W}  ([w, \partial_{j-1} \mathbbm{1}_{[v_1,..,.v_j]}])  = \\
&  \Cone_{X'}  (\mathbbm{1}_{[v_1,\dots,v_j]}) -  \Cone_{X'}  (\Cone_{X_w \cap X'} (\partial_{j-1} \mathbbm{1}_{[v_1,\dots,v_j]})) +\\
& [w,  \Cone_{X_w \cap X'}  (\partial_{j-1} \mathbbm{1}_{ [v_1,\dots,v_j]})] =^{\text{the cone equation for } \Cone_{X_w \cap X'}} \\
&  \Cone_{X'}  (\mathbbm{1}_{[v_1,\dots,v_j]}) -  \Cone_{X'} ( \mathbbm{1}_{[v_1,\dots,v_j]}) + \Cone_{X'}  (\partial_j \Cone_{X_w \cap X'} ( \mathbbm{1}_{[v_1,\dots,v_j]})) +\\
& [w,  \Cone_{X_w \cap X'}  (\partial_{j-1} \mathbbm{1}_{ [v_1,\dots,v_j]})] = \\
& \Cone_{X'}  (\partial_j \Cone_{X_w \cap X'} ( \mathbbm{1}_{[v_1,\dots,v_j]})) + [w,  \Cone_{X_w \cap X'}  (\partial_{j-1} \mathbbm{1}_{ [v_1,\dots,v_j]})].
\end{align*}


Using this computation,  we will verify the cone equation for $\Cone_{X' \cup W} $ (below we use the cone equation for $\Cone_{X'}$ and $\Cone_{X_w \cap X'}$  without noting it):
\begin{align*}
& \partial_{j+1} \Cone_{X' \cup W}  (\mathbbm{1}_{[w, v_1,\dots,v_j]}) = \\
& \partial_{j+1}  \Cone_{X'}  (\Cone_{X_w \cap X'}  (\mathbbm{1}_{[v_1,\dots,v_j]})) - \partial_{j+1} [w,  \Cone_{X_w \cap X'}  (\mathbbm{1}_{[v_1,\dots,v_j]})]  = \\
&- \Cone_{X'}  (\partial_j \Cone_{X_w \cap X'}  (\mathbbm{1}_{[v_1,\dots,v_j]})) +\Cone_{X_w \cap X'}  (\mathbbm{1}_{[v_1,\dots,v_j]}) - \\
& \left( \Cone_{X_w \cap X'}  (\mathbbm{1}_{[v_1,\dots,v_j]}) - [w, \partial_j \Cone_{X_w \cap X'}  (\mathbbm{1}_{[v_1,\dots,v_j]})] \right) = \\
& - \Cone_{X'}  (\partial_j \Cone_{X_w \cap X'}  (\mathbbm{1}_{[v_1,\dots,v_j]})) -  [w,  \Cone_{X_w \cap X'}  (\partial_{j-1} \mathbbm{1}_{[v_1,\dots,v_j]})] +  [w,  \mathbbm{1}_{[v_1,\dots,v_j]}] = \\
& - \Cone_{X' \cup W}  (\partial_j \mathbbm{1}_{[w, v_1,\dots,v_j]})   + \mathbbm{1}_{[w,v_1,\dots,v_j]}
\end{align*}
as needed.

In order to conclude the proof,  we will show that $\Rad_j^{n-1} (\Cone_{X' \cup W}) \leq R_{j-1} '' (R_j ' +1)$.  For $[v_0,\dots,v_j] \in \overrightarrow{X'} (j)$, 
$$\vert \supp (\Cone_{X' \cup W} (\mathbbm{1}_{[v_0,\dots,v_j]} )) \vert = \vert \supp (\Cone_{X' } (\mathbbm{1}_{[v_0,\dots,v_j]} )) \vert \leq R_j ' \leq  R_{j-1} '' (R_j ' +1)  .$$

For $[w,v_1,\dots,v_j] \in \overrightarrow{X' \cup W} (j)$ where $w \in W$, 
\begin{align*}
& \vert \supp (\Cone_{X' \cup W} (\mathbbm{1}_{[w, v_0,\dots,v_j]} )) \vert  \leq \\
& \vert \supp (\Cone_{X'}  (\Cone_{X_w \cap X'} (\mathbbm{1}_{[v_1,\dots,v_j]})) ) \vert + \vert \supp ([w,  \Cone_{X_w \cap X'}  (\mathbbm{1}_{[v_1,\dots,v_j]})]) \vert \leq \\
& R_j ' R_{j-1} '' + R_{j-1} '' = R_{j-1} '' (R_j '  +1 )
\end{align*}
as needed.
\end{proof}

\section{The Kac--Moody--Steinberg complexes}
\label{The KMS complexes sec}

In this section, we describe the families of simplicial complexes, for which we will prove that they give rise to cosystolic expanders. The main source, what we call KMS complexes, are infinite families of coset complexes over finite quotients of so-called Kac--Moody--Steinberg groups. They were introduced in \cite{hdxfromkms}, where it was shown that they are local spectral high-dimensional expanders. Our methods also apply to the complexes introduced in \cite{DLYZ}, since the links of these complexes are isomorphic to links of certain KMS complexes.

\subsection{Root systems and Chevalley groups}
We start by recalling some facts about generalized Cartan matrices, the associated root systems and Chevalley groups. These are the building blocks of KMS groups. 

Most of the structure of a KMS group is encoded in its generalized Cartan matrix, which also gives rise to a Dynkin diagram and a root system.

\begin{definition} \label{def: GCM and DyDi}
A generalized Cartan matrix (GCM) is a matrix $A = (A_{ij})_{i,j \in I} \in \mat_n(\Z)$ such that $A_{ii}=2$ for all $i \in I$, $A_{ij} \leq 0$ for all $i \neq j \in I$ and $A_{ij}= 0 \iff A_{ji}=0$. We call $\operatorname{rank}(A) := \vert I \vert$ the rank of $A$.

Every GCM $A=(A_{i,j})_{i,j \in I}$ corresponds to a Dynkin diagram in the following way. As vertex set, we take the index set $I$ and two vertices $i,j$ are connected by $|A_{i,j}|$ edges if $A_{i,j}A_{j,i} \leq 4$ and $|A_{i,j}| \geq |A_{j,i}|$, and these edges are equipped with an arrow pointing towards $i$ if $|A_{i,j}|>1$ (a bi-directed edge is possible). If $A_{i,j}A_{j,i}>4$, the vertices $i$ and $j$ are connected by a bold-faced labelled edge labelled with the ordered pair of integers $|A_{i,j}|,|A_{j,i}|$.

A GCM is \emph{irreducible} if there exists no non-trivial partition $I = I_1 \cup I_2$ such that $A_{i_1,i_2} = 0$ for all $i_1 \in I_1, i_2 \in I_2$.
\end{definition}

\begin{figure}

\begin{align*}
	\text{Classical types } & \begin{cases}
	A_n &
	\begin{tikzpicture}[line width=1pt, scale=.8, baseline=-.6ex]
		\draw (0,0) -- (2,0);
		\draw[dotted] (2,0) -- (4,0);
		\draw (4,0) -- (6,0);
		\diagnode{(0,0)}
		\diagnode{(1,0)}
		\diagnode{(2,0)}
		\diagnode{(4,0)}
		\diagnode{(5,0)}
		\diagnode{(6,0)}
	\end{tikzpicture} \\[1ex]
	B_n &
	\begin{tikzpicture}[line width=1pt, scale=.8, baseline=-.6ex]
		\draw (0,0) -- (2,0);
		\draw[dotted] (2,0) -- (4,0);
		\draw (4,0) -- (5,0);
		\draw[double distance=2.3pt] (5,0) -- (6,0);
		\draw[line width=.7pt] (5.3,.3) -- (5.6,0) -- (5.3,-.3);
		\diagnode{(0,0)}
		\diagnode{(1,0)}
		\diagnode{(2,0)}
		\diagnode{(4,0)}
		\diagnode{(5,0)}
		\diagnode{(6,0)}
	\end{tikzpicture} \\[1ex]
	C_n &
	\begin{tikzpicture}[line width=1pt, scale=.8, baseline=-.6ex]
		\draw (0,0) -- (2,0);
		\draw[dotted] (2,0) -- (4,0);
		\draw (4,0) -- (5,0);
		\draw[double distance=2.3pt] (5,0) -- (6,0);
		\draw[line width=.7pt] (5.7,.3) -- (5.4,0) -- (5.7,-.3);
		\diagnode{(0,0)}
		\diagnode{(1,0)}
		\diagnode{(2,0)}
		\diagnode{(4,0)}
		\diagnode{(5,0)}
		\diagnode{(6,0)}
	\end{tikzpicture} \\[1ex]
	D_n &
	\begin{tikzpicture}[line width=1pt, scale=.8, baseline=-.6ex]
		\draw (0,0) -- (2,0);
		\draw[dotted] (2,0) -- (4,0);
		\draw (4,0) -- (5,0);
		\draw (5,0) -- (6,.4);
		\draw (5,0) -- (6,-.4);
		\diagnode{(0,0)}
		\diagnode{(1,0)}
		\diagnode{(2,0)}
		\diagnode{(4,0)}
		\diagnode{(5,0)}
		\diagnode{(6,.4)}
		\diagnode{(6,-.4)}
	\end{tikzpicture}
	\end{cases} \\[2ex]
	\text{Exceptional types } & \begin{cases}
	E_6, E_7, E_8 &
	\begin{tikzpicture}[line width=1pt, scale=.8, baseline=-.6ex]
		\draw (0,0) -- (4,0);
		\draw[dotted] (4,0) -- (6,0);
		\draw (2,0) -- (2,-1);
		\diagnode{(0,0)}
		\diagnode{(1,0)}
		\diagnode{(2,0)} \diagnode{(2,-1)}
		\diagnode{(3,0)}
		\diagnode{(4,0)}
		\diagnode{(5,0)}
		\diagnode{(6,0)}
	\end{tikzpicture} \\[5ex]
	F_4 & 	\begin{tikzpicture}[line width=1pt, scale=.8, baseline=-.6ex]
		\draw (0,0) -- (1,0);
		\draw[double distance=2.3pt] (1,0) -- (2,0);
		\draw (2,0) -- (3,0);
		\draw[line width=.7pt] (1.3,.3) -- (1.6,0) -- (1.3,-.3);
		\diagnode{(0,0)}
		\diagnode{(1,0)}
		\diagnode{(2,0)}
		\diagnode{(3,0)} \end{tikzpicture} \\[1ex]
	G_2 &	\begin{tikzpicture}[line width=1pt, scale=.8, baseline=-.6ex]
		\draw[double distance=3pt] (0,0) -- (1,0);
		\draw (0,0) -- (1,0);
		\draw[line width=.7pt] (0.3,.3) -- (0.6,0) -- (0.3,-.3);
		\fill (0,0) circle (.11);
		\fill (1,0) circle (.11); \end{tikzpicture} \\[1ex] \end{cases}
\end{align*}	
\caption{The Dynkin diagrams of irreducible spherical root systems\label{table:dynkin}}	
\end{figure}
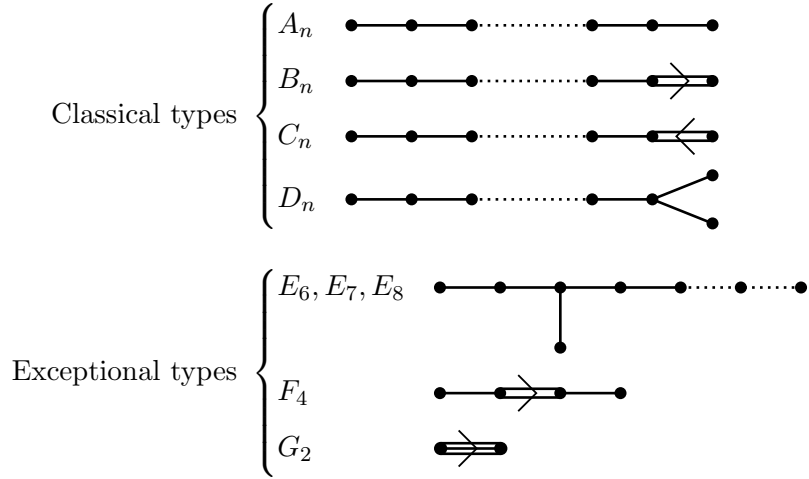

\begin{definition}
Let $A$ be a $(d+1)\times (d+1)$ generalized Cartan matrix with index set $I=\{0,\dots,d\}$, and let $J \subseteq I$.
\begin{enumerate}[label=(\roman*)]
\item $A$ gives rise to a Coxeter matrix $M$ of the same size (different GCMs can give rise to the same Coxeter matrix) which in turn defines a (finite or infinite) Coxeter group (see \Cref{subsec: types of buildings}).
    \item  The subset  $J$ is called \emph{spherical} if $A_J = (A_{i,j})_{i,j \in J}$ is of spherical type, meaning that the associated Coxeter group is finite (see e.g. \cite[Chapter 6.4.1]{bourbaki2008lie}).
Given $n  \geq 2$, $A$ is \emph{$n$-spherical} if every subset $J \subseteq I$ of size $n$ is spherical. 

\item We denote by $Q_A$ the set of spherical subsets of $I$ associated to the generalized Cartan matrix $A$.

\item A generalized Cartan matrix $A$ is \emph{purely $n$-spherical} if every spherical subset $J \subseteq I$ is contained in a spherical subset of size $n$.  (In particular, no set of size $n+1$ is spherical for a purely $n$-spherical generalized Cartan matrix.) Note that any such GCM is irreducible (if it could be decomposed, each irreducible component needs to be spherical and thus the whole GCM would be spherical).

\item Analogously, we introduce the following terminology and call a GCM \emph{$n$-classical} if each irreducible factor of $A_J = (A_{i,j})_{i,j \in J}$ is of classical type $A_k,B_k,C_k,D_k$ (see e.g. \cite[Chapter 6.4.1]{bourbaki2008lie}) for each $J \subseteq I$ with $\lvert J \rvert \leq n$. 

\item To a generalized Cartan matrix $A$ we can associate the following sets:
\begin{itemize}
    \item a set of simple roots $\Pi = \{\alpha_i \mid i \in I\}$,
    \item a set of real roots $\Phi \subseteq \bigoplus_{i \in I} \Z \alpha_i$,
    \item two sets, one of positive and one of negative real roots $\Phi^+= \bigoplus_{i \in I}\N \alpha_i \cap \Phi, \Phi^- = - \Phi^+$ .
\end{itemize}   
Note that $\Phi = \Phi^+ \sqcup \Phi^-$. More details can be found e.g. in \cite[Chapter 3.5]{marquis2018introduction}.
\end{enumerate}

\end{definition}
\begin{remark}
\begin{enumerate}[label=(\roman*)]
\item Irreducible spherical diagrams have been classified. All possible diagrams are given in \Cref{table:dynkin}. 
\item All rank $n+1$, $n$-spherical diagrams have been classified and are precisely the affine diagrams and the so-called compact hyperbolic diagrams. The latter have rank at most 5, where the diagram of rank 5 is exactly the one that we will have to exclude to get $n$-classicality below. See e.g. \cite[Chapter 6.9]{HumphreyReflGrp} for the classification result. 

The non-spherical, rank $n+1$, $n$-classical diagrams, with $n \geq 3$ are precisely the non-spherical rank $n+1$ diagrams which are $n$-spherical but not one of the following
\begin{itemize}
\item affine diagrams $\tilde{E}_6,\tilde{E}_7,\tilde{E}_8,\tilde{F}_4,\tilde{G}_2$,
\item twisted affine diagrams $E_6^{(2)}, D_4^{(3)}$ (in the notation of \cite[Table 5.1]{marquis2018introduction})
\item the compact hyperbolic diagram \begin{tikzpicture}[line width=1pt, scale=.8, baseline=-.6ex]
		\draw (0,0) -- (1,0);
		\draw[double distance=2.3pt] (1,0) -- (2,0);
		\draw (2,0) -- (3,0);
		\draw[line width=.7pt] (1.3,.3) -- (1.6,0) -- (1.3,-.3);
		\draw (0,0) -- (1.5,1);
		\draw (3,0) -- (1.5,1);
		\diagnode{(0,0)}
		\diagnode{(1,0)}
		\diagnode{(2,0)}
		\diagnode{(3,0)}
		\diagnode{(1.5,1)} 
		\end{tikzpicture}.
\end{itemize}
  
\end{enumerate}

\end{remark}

\begin{definition} \label{def: chevalley group}
Corresponding to any irreducible, spherical GCM $\Ao$ with root system $\mathring \Phi$ of rank at least 2, and any finite field $K$, there is an associated universal (or simply connected) Chevalley group, denoted $\chev_{\Ao}(K)$. Abstractly, it is generated by symbols $x_\alpha(s)$ for $\alpha \in \mathring \Phi$ and $s \in K$, subject to the relations
$$
\begin{aligned}
x_\alpha(s) x_\alpha(u) & =x_\alpha(s+u) \\
{\left[x_\alpha(s), x_\beta(u)\right] } & =\prod_{i, j>0} x_{i \alpha+j \beta}\left(C_{i j}^{\alpha, \beta} s^i u^j\right) \quad(\text {for } \alpha+\beta \neq 0) \\
h_\alpha(s) h_\alpha(u) & =h_\alpha(s u) \quad(\text { for } s,u \neq 0), \\
\text { where } \quad h_\alpha(s) & =n_\alpha(s) n_\alpha(-1) \\
\text { and } \quad n_\alpha(s) & =x_\alpha(s) x_{-\alpha}\left(-s^{-1}\right) x_\alpha(s) .
\end{aligned}
$$
Note that the $C_{ij}^{\alpha, \beta}$ are integers called structure constants that can be found in \cite{carter1989simple}.
\end{definition}
\begin{remark}
Let $K$ be a field and let $K[t]$ denote the polynomial ring in one variable over $K$. The simply-connected Chevalley group $\chev_{\Ao}(K[t])$ is, similar to the case of a Chevalley group over a field, generated by elements $x_\alpha(s)$ for $\alpha \in \mathring \Phi, s \in K[t]$ that satisfy the relations above, where for the third relation we have to add the extra assumption that $u,s$ are invertible in $K[t]$. This can be found e.g. in \cite{rehmann1975prasentationen}.
\end{remark}

\subsection{Kac--Moody--Steinberg groups}
We can now put the ingredients together to define what a Kac--Moody--Steinberg group is and present some of its properties. For more details, see \cite[Chapter 3]{hdxfromkms} and references therein. 
\begin{definition}
Let $K$ be a field. Let $A=(A_{i,j})_{i,j\in I}$ be a GCM over the index set $I$ which is 2-spherical. Let $Q_A=\{J \subseteq I \mid A_J:=(A_{ij})_{i,j \in J} \text{ is spherical}\}$. Let $\Phi$ be the root system associated to $A$ with simple roots $\Pi = \{\alpha_i \mid i \in I\}$. For $J\in Q_A$ set 
$$U_J:= \langle x_{\alpha_i}(s) \mid s \in K, i \in J \rangle \leq \chev_{A_J}(K) $$
the group generated by all simple roots in the Chevalley group of type $A_J$ over $K$. Note that if $L \subset J$ then we have a natural inclusion $U_L \hookrightarrow U_J$ by sending $x_{\alpha_i}(s)\in U_L$ to the same generator in $U_J$.  
The KMS group of type $A$ over $K$ is defined as the free product of the $U_J, J\in Q_A$ modulo the natural inclusions:
$$\uak = \freeprod_{J \in Q_A} U_J/(U_L \hookrightarrow U_J, L\subseteq J).$$ 
\end{definition}
Note that $U_J \hookrightarrow \uak$ and we will denote the image of $U_J$ in $\uak$ again by $U_J$. These subgroups are called the local groups of $\uak$. For $i,j\in I$ we will write $U_i:= U_{\{i\}}$ and $U_{i,j}=U_{\{i,j\}}$.

\begin{remark} \label{rmk: presentation of KMS}
For two roots $\alpha, \beta \in \Phi$ we write $]\alpha, \beta[_\N = \left\{n_1 \alpha + n_2 \beta \in \Phi \mid n_1,n_2 \in \N^* \right\}$ and $[\alpha,\beta]_\N = \left]\alpha, \beta\right[_\N \cup \left\{\alpha,\beta\right\}$.

We have the following abstract presentation for the KMS groups:
$$\mathcal{U}_A(K)= \Bigl\langle u_\beta(t) \text{ for } t \in K,  i,j \in I , \beta \in [\alpha_i,\alpha_j]_{\mathbb{N}} \ | \ \mathcal{R} \Bigr\rangle$$
where the set of relations $\mathcal{R}$ is defined as:

\indent for all $i,j \in I , \  \{\alpha,\beta \} \subseteq [\alpha_i, \alpha_j]_{\mathbb{N}},   \ t,s\in K$: 
\begin{align*}
u_\alpha(t) u_\alpha(s) &= u_\alpha(t+s) \\
[u_\alpha(t),u_\beta(s)] &= \prod_{\gamma = k\alpha + l\beta \in ]\alpha,\beta[_{\mathbb{N}}} u_\gamma \left(C_{k,l}^{\alpha \beta} t^k s^l\right).
\end{align*}
The constants $C_{k,l}^{\alpha,\beta}$ are the same structure constants as in the presentation of the Chevalley group.

Since $A$ is 2-spherical, the subgroups $U_\beta$, for $ \beta \in \left[ \alpha_i, \alpha_j\right] $, are contained in the group generated by the root groups $U_{\alpha_i}, U_{\alpha_j}$ (see \cite[Proposition 7]{Abr}).
\end{remark}

\subsection{KMS complexes}
In this section, we recall some facts about coset complexes and put everything together to define KMS complexes.

KMS complexes are coset complexes. More details on coset complexes can be found for example in \cite{KO2023high},  and we will only give a brief overview. 

\begin{definition}
Let $G$ be a group and $\mathcal{H}=(H_0,\dots , H_n)$ be a family of subgroups of $G$. Then the coset complex $\CC(G;(H_0, \dots, H_n))$ is defined to be the simplicial complex 
\begin{itemize}
\item with vertex set $\bigsqcup_{i=0}^n G/H_i$
\item such that a set of vertices $\left\{ g_1H_{i_1},\dots,g_kH_{i_k} \right\}$ forms a $(k-1)$-simplex in $\CC(G; \mathcal{H})$ if and only if $\bigcap_{j=1}^k g_j H_{i_j} \neq \emptyset$.
\end{itemize}
\end{definition}

Note that the set of vertices is partitioned into $n+1$ subsets $G/H_i$ for $i = 0,\dots,n$ such that for any simplex $\sigma \in \CC(G;\mathcal{H})$ we have $|\sigma \cap G/H_i| \leq 1$ for all $i$. Hence, we say that $\CC(G;\mathcal{H})$ is a $(n+1)$-partite simplicial complex.


\begin{definition} \label{def: type of a face}
    The \emph{type} of a simplex $\left\{ g_1H_{i_1},\dots,g_kH_{i_k} \right\}$ in $\CC(G;\mathcal{H})$ is the set of indices $\{i_1,\dots,i_k\} \subseteq \{0,\dots,n\}$. Given a type $\emptyset \neq T \subseteq \{0,\dots,n\}$ we write 
    $$H_T:= \bigcap_{i \in T} H_i \qquad \text{ and set } \qquad H_{\emptyset}:= \langle H_0,\dots,H_n \rangle \leq G.$$
\end{definition}

The following well-known fact will be very useful. 

\begin{proposition} \label{prop: links in CC} 
Let $\sigma$ be a face in $\CC (G , \mathcal{H})$ of type $T \neq \emptyset$. Then the link of $\sigma$ is isomorphic to the coset complex $\CC \left(H_T ,\left(H_{T \cup\{i\}}: i \notin T\right)\right)$.
\end{proposition}

KMS complexes are coset complexes over certain finite quotients of KMS groups. The precise definition is as follows. 

\begin{definition} \label{def: KMS complex}
Let $A$ be a GCM over the index set $I$ which is $(\vert I \vert-1)$-spherical but non-spherical.  Let $K$ be a finite field of size $q \geq 4$. Let $\phi: \uak \to G$ be a finite quotient of $\uak$ such that
\begin{enumerate}
\item $\phi|_{U_J}$ is injective for all $J \in Q_A$ (i.e. $\phi$ is injective on the local groups), 
\item $\phi(U_J \cap U_L) = \phi(U_J) \cap \phi(U_L)$ for all $J,L \in Q_A$.
\end{enumerate}
Then we set 
$$X= \CC(G; (\phi(U_{I \setminus \{j \}}))_{j \in I})$$ 
and call $X$ a KMS complex of type $A$ over $K$. 
\end{definition}

\begin{theorem}{\cite[Theorem 4.3]{hdxfromkms}}
\label{KMS spectral thm}
Let $A$ be a GCM over the index set $I$ which is non-spherical but $\lvert I \rvert -1$-spherical. Let $K$ be a finite field of size $q \geq 4$. Let $\phi_i: \uak \to G_i, i\in I$ be a family of finite quotients of $\uak$ such that $X_i = \CC(G_i, (\phi_i(U_{I \setminus \{j \}}))_{j \in I}), i \in \N$ is a family of KMS complexes and $\lvert G_i \rvert \overset{i \to \infty}{\longrightarrow} \infty$.
If $q$ is such that $\sqrt{\frac{3}{q}} \leq \frac{1}{\lvert I \rvert -1}$ then the family $(X_i)_{i\in \N}$ is a family of bounded degree local spectral expanders. 

Moreover,  if $q \gg (\vert I \vert -1)^2$,  then the family $(X_i)_{i\in \N}$ is a family of bounded degree $\frac{2}{\sqrt{q}}$-local spectral expanders. 
\end{theorem}

\begin{example} \label{expl: explicit construction KMS complex}
	The following is an example of an $n$-dimensional, $n$-classical KMS complex. Consider $\SL_{n+1} (\mathbb{F}_q [t])$ where $n \geq 2$ and $q$ is an odd prime power such that $q \geq 5$. For $1 \leq i,j \leq n+1,  i \neq j$ and $f \in  \mathbb{F}_q [t]$,  we denote $e_{i,j} (f) \in \SL_{n+1} (\mathbb{F}_q [t])$ to be the matrix with $1$'s along the main diagonal,  $f$ in the $(i,j)$-th entry and $0$'s in all other entries.  Let
$$B = \left\{e_{i,i+1} (1); i = 1,\dots n\right\} \cup \{ e_{n+1,1} (t)\}$$ and let $\mathcal{U} (\mathbb{F}_q) < \SL_{n+1} (\mathbb{F}_q [t])$ be the subgroup generated by $B$, i.e.
$$\mathcal{U} (\mathbb{F}_q) = \langle e_{1,2} (1),  e_{2,3} (1),\dots ,e_{n,n+1} (1),  e_{n+1,1} (t) \rangle.$$
We define subgroups $H_0,\dots, H_n < \mathcal{U} (\mathbb{F}_q)$ as 
\begin{align*}
H_{0} &=  \langle e_{1,2} (1),  e_{2,3} (1),\dots,e_{n,n+1} (1) \rangle = \langle B \setminus \{e_{n+1,1}(t)\} \rangle, \\
H_{1} &=  \langle  e_{2,3} (1),\dots,e_{n,n+1} (1),  e_{n+1,1} (t) \rangle = \langle B \setminus \{e_{1,2}(1)\} \rangle, \\
&\vdots \\
H_{i} &=  \langle  B \setminus \{e_{i,i+1}\} \rangle, \\
& \vdots \\
 H_{n} &= \langle  e_{1,2} (1),  e_{2,3} (1),\dots,e_{n-1,n} (1),  e_{n+1,1} (t) \rangle ,
\end{align*}
and denote $\mathcal{H} = (H_0,\dots , H_n)$.

Let $f \in \mathbb{F}_q [t]$ be an irreducible polynomial of degree $s >1$. Denote $\phi_{f} : \SL_{n+1} (\mathbb{F}_q [t]) \rightarrow \SL_{n+1}(\mathbb{F}_q[t]/(f)) \cong \SL_{n+1} (\mathbb{F}_{q^s})$ defined as follows: For every matrix $A \in  \SL_{n+1} (\mathbb{F}_q [t]),  A= (A_{i,j})$ where $A_{i,j} \in \mathbb{F}_q [t]$ define 
$$\phi_f (A) =  (A_{i,j} + ( f)).$$
Restricting $\phi_f$ to $\mathcal{U} (\mathbb{F}_q)$, we get $G_f = \phi_f (\mathcal{U} (\mathbb{F}_q))$ which is a finite quotient of $\mathcal{U} (\mathbb{F}_q)$. It is proven in \cite{hdxfromkms}, that $\phi_f$ is injective on $H_0,\dots,H_n$. The KMS complex of $G_f = \phi_f (\mathcal{U} (\mathbb{F}_q))$ is defined to be the coset complex $X_f = \CC (G_f , \mathcal{H}_f)$, where $\mathcal{H}_f=(\phi_f(H_0),\dots, \phi_f(H_n))$.   

\end{example}
\begin{remark}
	Similar examples of $n$-classical KMS complexes can be constructed for all Chevalley groups of classical type, see \cite[Chapter 5]{hdxfromkms} for more details.
\end{remark}

\section{Describing classical buildings and their opposite complexes}
\label{sec: describing buildings}
In the first part of this section, we briefly discuss Coxeter complexes. We then describe buildings and their opposite complexes from a group theoretic view point. In this setting, we show that the links of KMS complexes are opposite complexes of spherical buildings. In the third part, we describe a more geometric view point using flag complexes. This context will be used to show that opposite complexes are coboundary expanders.

\subsection{Coxeter complex} \label{subsec: types of buildings}

To each Dynkin diagram with set of vertices $I$ as in \Cref{def: GCM and DyDi} we can associate a group, called a Coxeter group, in the following way:
$$\langle s_i; i \in I \mid s_i^2 = 1, (s_is_j)^{m_{ij}} = 1 ; i,j \in I, i \neq j\rangle$$
where 
\begin{align*}
\begin{tikzpicture}
		\diagnode{(0,0)}
		\diagnode{(1,0)}
		\node[above] at (0,0.2) {$i$};
		\node[above] at (1,0.1) {$j$};
\end{tikzpicture} &\Rightarrow m_{ij} = 2, \qquad \begin{tikzpicture}
		\draw (0,0) -- (1,0);
		\diagnode{(0,0)}
		\diagnode{(1,0)}
		\node[above] at (0,0.2) {$i$};
		\node[above] at (1,0.1) {$j$};
\end{tikzpicture} \Rightarrow m_{ij} = 3\\
\begin{tikzpicture}
		\draw[double distance=2.3pt] (0,0) -- (1,0);
		\draw[line width=.7pt] (0.3,.3) -- (0.6,0) -- (0.3,-.3);
		\diagnode{(0,0)}
		\diagnode{(1,0)}
		\node[above] at (0,0.2) {$i$};
		\node[above] at (1,0.1) {$j$};		
\end{tikzpicture} &\Rightarrow m_{ij} = 4,  \qquad \begin{tikzpicture}
		\draw[double distance=3pt] (0,0) -- (1,0);
		\draw (0,0) -- (1,0);
		\draw[line width=.7pt] (0.3,.3) -- (0.6,0) -- (0.3,-.3);
		\diagnode{(0,0)}
		\diagnode{(1,0)}
		\node[above] at (0,0.2) {$i$};
		\node[above] at (1,0.1) {$j$};
\end{tikzpicture}\Rightarrow m_{ij} = 6
\end{align*}
and $m_{ij} = +\infty$ if $i,j$ are connected by a bold-faced labelled edge, meaning that in this case there is no relation between $s_i$ and $s_j$.
Note that $m_{ij}$ is defined to be symmetric in $i,j$.  This implies that the root systems of type $B_n$ and $C_n$ give rise to the same Coxeter group.
The Coxeter complex associated to the Coxeter group $W$ with set of generators $S$ can be viewed as the coset complex
$$\CC(W,(\langle S\setminus \{t\}\rangle)_{t \in S}).$$

More details can be found for example in \cite{AB2008}.

\subsection{Buildings and BN-pairs}
One way to define a building is as a simplicial complex with a family of subcomplexes called apartments such that each apartment is isomorphic to the same Coxeter complex. The type of the Coxeter complex, which can be described by its Coxeter matrix, is also what we call the type of the building. Any two simplices of the building have to be contained in a common apartment and any two apartments are isomorphic with an isomorphism that acts as the identity on the intersection.  

One way to construct buildings is via groups with BN-pair (and for thick, irreducible buildings of dimension 2 and larger, all buildings can be constructed that way, as was proven by Jacques Tits in \cite{TitsBuildings}). Note that if we construct a building from a BN-pair of a Chevalley group with underlying root system of some type $X$, then the building will have the same type. In particular, since type $B_n$ and $C_n$ give rise to the same Coxeter complex, they also determine the same type of building.

Given a group $G$, a BN-pair of $G$ is a pair of subgroups $B,N \leq G$ satisfying certain properties, see \cite[Definition 6.55]{AB2008}. In particular, $G = BN$ and setting $T=B \cap N$ we have that $W= N/T$ admits a set of generators $S$, such that $(W,S)$ is a Coxeter system. $W$ is also called the Weyl group of the BN-pair. 
To describe the building, we need to define the parabolic subgroups $P_J:= B \langle J \rangle B$ for $J \subset S$ (here $\langle J \rangle$ is the subgroup generated by the elements of $J$ inside $W$). The building can then be described as the following coset complex:
$$\Delta(G,B) := \CC(G; (P_{S\setminus \{j \}})_{j \in S}).$$

Since $P_J \cap P_L = P_{J \cap L}$ for any $J,L \subseteq S$ we can see $\Delta(G,B)$ also as the poset $\{gP_J \mid g \in G, J \subseteq S\}$ ordered by reversed inclusion. Then the simplex $\{g P_{S\setminus \{j\}} \mid j \in J\} \in \CC(G; (P_{S\setminus \{j \}})_{j \in S})$ corresponds to $gP_{S\setminus J}$ for $J \subseteq S, g \in G$.
In particular, the set of maximal simplices, called chambers, corresponds to $G/B$.

If the Weyl group is finite, we call the building spherical. This coincides with our definition of calling a GCM/root system spherical. A finite Weyl group has a unique longest element denoted by $w_0$, where longest means that is generated by the largest number of generators, counted with multiplicities. 

Two chambers $gB, hB, g,h\in G$ are called opposite if and only if $g^{-1}h \in Bw_0B$, denoted by $gB \text{ op } hB$.

Given a simplex $\sigma \in \Delta(G,B)$, the complex opposite $\sigma$ is defined as follows: 
$$\Delta^0(\sigma):= \{\tau \in \Delta(G,B) \mid \text{ there exist chambers } c,d \in \Delta(G,B) \text{ such that } \tau \subseteq c, \sigma \subseteq d, c \text{ op } d\}.$$

To see the connection between opposite complexes and the links of KMS complexes, note that every Chevalley group has a BN-pair in the following way. Let $A=(A_{ij})_{i,j \in I}$ be a spherical GCM of rank at least 2, $K$ a field and $G= \chev_A(K) = \langle x_\alpha(\lambda) \mid \alpha \in \Phi, \lambda \in K \rangle$, using the notation of \Cref{def: chevalley group}.
We set
\begin{alignat*}{3}
U_\alpha &:= \langle x_{\alpha}(\lambda) \mid \lambda \in K \rangle, \alpha \in \Phi,  \qquad &&U^+ &&:= \langle U_ \alpha \mid \alpha \in \Phi^+ \rangle, \qquad T := \langle h_\alpha(\lambda) \mid \alpha \in \Phi, \lambda \in K\rangle,\\
B &:= U^+ T, \qquad &&N &&:= \langle T, n_{-\alpha_i}(-\lambda^{-1}) \mid i \in I, \lambda \in K^* \rangle.
\end{alignat*}
For example \cite[Theorem 7.115]{AB2008} shows that this is indeed a BN-pair. 

Let $C_0= 1B$ denote the fundamental chamber of $\Delta(G,B)$. Recall that $w_0$ denotes the longest element on the Weyl group $W$. It satisfies $Bw_0B = U^+ w_0 B$. We can describe the opposite complex in this set-up as 
\begin{align*}
\Delta^{0}(C_{0}) &= \{ gP_{J} \mid g \in G, J \subseteq S: \exists h \in U^+w_{0}B: gP_{J}\supseteq hB  \}\\
&= \{ gP_{J} \mid \exists h \in U^+w_{0}B: h \in gP_{J} \} \\
&= \{ hP_{J} \mid h \in U^+w_{0}B \} \\
&= \{ uw_{0}P_{J} \mid u \in U^+ \}.
\end{align*}

This leads to the following proposition. 
\begin{proposition}{\cite[Chapter 3.2]{hdxfromkms}} \label{prop: opposite complex and U+}
Let $G$ be a Chevalley group with BN-pair as described above. Then 
$$\CC(U^+, (U_{I\setminus \{i\}})_{i \in I}) \cong \Delta^0(C_0).$$
\end{proposition}

Since we assume that $A$ is 2-spherical and $\lvert K \rvert \geq 4$ we have $U^+ = \langle U_\alpha \mid \alpha \in \Phi^+ \rangle = \langle U_{\alpha_i} \mid i \in I \rangle$ (see e.g. the comment before \cite[Lemma 3.7]{hdxfromkms}). Combining this with the above proposition and \Cref{prop: links in CC} gives the following corollary. 
\begin{corollary}
Let $X$ be a KMS complex, and $\sigma \in X$ a face of dimension less than or equal to $ \dim X-2$ of type $\emptyset \neq T \subset I$. Then 
$$\lk_X(\sigma) \cong \Delta^0(C_0)$$
where $\Delta$ is the building obtained from the BN-pair structure of the Chevalley group of type $A_{I \setminus T}$.
\end{corollary}

\subsection{Geometric constructions of buildings} \label{subsec: geometric construction of buildings}
In the following section, we describe how buildings of type $A_n$, $C_n$ and $D_n$ can be constructed as flag complexes of certain sets equipped with an incidence relation. We furthermore describe the opposite complexes in this set-up. The general framework is as follows.

\begin{definition}
Let $X$ be a set and let $I \subseteq X \times X$ be a reflexive and symmetric relation. For $a,b \in X$ we write $aIb$ if and only if $(a,b)\in I$ and call $a$ and $b$ incident. The structure gives rise to the following simplicial complex, called the flag complex of $X$: 
$$\flag(X) = \{\sigma \subseteq X \mid \forall a,b \in \sigma: aIb \}.$$
\end{definition} 

\begin{remark}
We note that by definition,  $\flag(X)$ is always a clique complex.
\end{remark}

\subsubsection{Buildings of type $A_n$} \label{subsubsec: geom descr A_n}
Let $n\in \N, n\geq 1$. Any building of type $A_n$ can be described in the following way, see e.g. \cite{TitsBuildings}. 
Let $K$ be a field or skew field, $V$ an $(n+2)$-dimensional vector space over $K$. Set $X:= X(V) := \{0 <U<V \mid U \text{ is a } K \text{-subspace of } V \}$. Two subspaces $U,W \in X$ are called incident if $U \subseteq W$ or $W \subseteq U$. Then $\Delta = \flag X$ is a building of type $A_{n+1}$. 

To describe the opposite complex in this context, we need the following notation, see \cite[Definition 9]{Abr}.
\begin{definition} \label{def: A_n case X_E}
\begin{enumerate}
\item Two subspaces $U,W$ of $V$ are called transversal (in $V$) if $U\cap W = 0$ or $U+W = V$. If this is the case, we write $U\pitchfork_V W$ or simply $U \pitchfork W$.
\item Let $\EE$ be a set of subspaces of $V$. Then we write $U \pitchfork \EE$ if $U\pitchfork E$ for all $E \in \EE$. We set
$$X_{\EE}(V) := \{U \in X \mid U \pitchfork \EE \} \qquad T_{\EE}(V) = \flag X_{\EE}(V).$$
\end{enumerate}
\end{definition}  
Note that if $0 <U \leq W <V$ then $U$ and $W$ are not transversal in $V$.
We get the following result. 
\begin{proposition}{\cite[Corollary 12]{Abr}} \label{prop: An oppo as flag}
For any simplex $\sigma = \{E_1 < \dots < E_r\} \in \Delta$ set $\EE(\sigma) = \{E_i \mid 1 \leq i \leq r \}$. Then 
$$\Delta^0(\sigma) = T_{\EE(\sigma)}(V).$$
\end{proposition}

\begin{definition} \label{def: class C_A}
We denote by $\mathbf{C}_A$ the class of all simplicial complexes $T_{\EE}(V)$ where $\dim(T_{\EE}(V)) = n\geq 0$, $K$ is a (skew) field, $V$ an $(n+2)$-dimensional vector space over $K$, $\mathcal{E}$ is a finite set of subspaces of $V$ such that if we set $e_{j}=\lvert \left\{ E \in \mathcal{E} \mid \dim E=j \right\} \rvert$ we have $\sum_{j=1}^{n+1} {{n}\choose{j-1}} e_{j} \leq \lvert K \rvert$. Furthermore, $T_{\EE}(V)$ is defined as described in \Cref{def: A_n case X_E}. 
\end{definition}

\subsubsection{Hermitian and pseudo-quadratic forms}
In this part, we recall facts about hermitian and pseudo-quadratic forms that are needed to describe the buildings of type $C_n$ and $D_n$. We follow \cite[Chapter 5]{Abr}, who in turn mostly follows \cite{TitsBuildings}.

Let $K$ be a skew field, $\sigma: K \to K, a \mapsto a^\sigma$ be an involution, i.e. an anti-automorphism of $K$ such that $\sigma^2 = \operatorname{id}_K$. Let $\epsilon \in \{1,-1\} \subset K$. If $\sigma \neq \operatorname{id}_K$ then we require $\epsilon = -1$. Let $V$ be a right $K$-vector space of dimension $m \in \N \cup \{\infty\}$.
Let $K_{\sigma, \epsilon}:=\left\{\alpha-\alpha^\sigma \varepsilon \mid \alpha \in K\right\}$ and $ K^{\sigma, \epsilon}:=\left\{\alpha \in K \mid \alpha+\alpha^\sigma \varepsilon=0\right\}$. Let $\Lambda$ be a form parameter relative to $(\sigma, \varepsilon)$, i.e. $\Lambda$ is a subgroup of $(K,+)$ satisfying $K_{\sigma, \varepsilon} \subseteq \Lambda \subseteq K^{\sigma, \varepsilon}$ and $\alpha^\sigma \Lambda \alpha \subseteq \Lambda$ for all $\alpha \in K$.

Let $f: V \times V \to K$ be a $(\sigma, \epsilon)$-hermitian form, i.e. $f$ is biadditive and for all $x,y\in V, a,b \in K$ we have $f(xa,yb) = a^\sigma f(x,y) b$ and $f(y,x) = f(x,y)^\sigma \epsilon$.
Let $Q: V \to K/\Lambda$ be a $(\sigma,\epsilon)$-quadratic form with associated $(\sigma, \epsilon)$-hermitian form $f$, i.e. $Q(xa) = a^\sigma Q(x) a +\Lambda$ and $Q(x+y) -Q(x) -Q(y) = f(x,y) + \Lambda$ for all $x,y \in V, a \in K$. 
If $\Lambda = K$ we require that $f$ is alternating. If $\Lambda \neq K$ then $f$ is uniquely defined by $Q$.

For $M\subseteq V$ we write $M^\perp := \{x \in V \mid f(x,M) = 0\}$.
A subspace $U<V$ is called 
\begin{itemize}
\item non-degenerate if $U \cap U^\perp = 0$,
\item totally isotropic if $U \subseteq U^\perp$ and $Q(U) = 0$.
\end{itemize}   
We require that there is at least one finite-dimensional maximal totally isotropic subspace. In that case, all maximal totally isotropic subspace have the same dimension. We denote this dimension by $n$ and call it the Witt index of $(V,Q,f)$. We require $0 < n< \infty$.
Furthermore, we require that $V^\perp = 0$. 
We say that the triple $(V,Q,f)$ is a pseudo-quadratic space if it is of the form described above. 

If it additionally satisfies $(m, \Lambda) \neq (2n, 0)$ then we say $(V,Q,f)$ is a thick pseudo-quadratic space.

Note that if $(V,Q,f)$ is a thick pseudo-quadratic space and $K$ is a finite field, then $\dim V \leq 2n+1$.
In general, if one only considers finite fields, especially of characteristic different from 2, the set-up can be substantially simplified, see \cite[Remark 9.3,9.4]{AB2008}. 

\subsubsection{Buildings of type $C_n$} \label{subsubsec: geom descr C_n}
Let $n \in \N, n\geq 1$ and $(V,Q,f)$ be a thick pseudo-quadratic space with Witt index $n$. Set $X := X(V) = \{0 < U < V \mid U \text{ is totally isotropic}\}$ with incidence relation given by containment as in the $A_n$ case and $\Delta = \flag X$. Then $\Delta$ is a thick building of type $C_n$ and every classical $C_n$ building can be obtained in this way (here classical means that the links of type $A_2$ correspond to Desarguesian planes, which is always the case for $n\geq 4$), see \cite[Theorem 8.22]{TitsBuildings}. 

\begin{definition}\label{def: C_n case N(E)}
Set $\mathcal{U} := \{0\leq U \leq V \mid \dim U < \infty \}, \mathcal{U}^\perp := \{U^\perp \mid u \in \mathcal{U}\}$ and $\mathcal{W}:= \mathcal{U} \cup \mathcal{U}^\perp$.  

Let $\EE$ be a finite subset of $\mathcal{W}$ such that $\EE^\perp = \EE$. Assume $K = \mathbb{F}_q$ is a finite field. 
Set $\mathcal{E}_j:=\{E \in \mathcal{E} \mid \operatorname{dim} E=j\}, e_j:= \lvert \mathcal{E}_j \rvert$ and $e_h^{(s)}:=\sum_{j=0}^{2 s}\binom{2 s}{j} e_{h+j}$ for $h \in \mathbb{N}, s \in \mathbb{N}_0$ and $h+2 s<m$. We define 
\begin{itemize}
\item $N(\mathcal{E}):=e_1^{(n-1)}$ if $f$ is alternating $(\Rightarrow m=2 n)$ and $Q=0$, 
\item $N(\mathcal{E}):=\left(e_1^{(n-1)}\right)^2$ if $m=2 n$ and $\sigma \neq \mathrm{id}$, 
\item $N(\mathcal{E}):=2 e_2^{(n-1)}$ if $m=2 n+1$,
\item $N(\mathcal{E}):=\max \left\{e_2^{(n-1)}+e_3^{(n-1)}+1,2 e_3^{(n-1)}\right\}$ if $m=2n+2$.
\end{itemize}
\end{definition}

\begin{definition}\label{def: X_E in C_n case}
For a set $\EE$ of subspaces of $V$ we define, similar to \Cref{def: A_n case X_E}, 
$$X_{\EE}(V) := \{U\in X \mid U \pitchfork \EE\}, \qquad T_{\EE}(V) = \flag X_{\EE}(V).$$ 
\end{definition}

\begin{proposition}{\cite[Corollary 15]{Abr}}
For any simplex $\tau = \{E_1 < \dots < E_r\} \in \Delta$ set $\EE(\tau) = \{E_i, E_i^\perp \mid 1 \leq i\leq r\}$. Then $\Delta^0(\tau) = T_{\EE(\tau)}(V)$.

\end{proposition}

\begin{definition}\label{def: class C_C}
With $\mathbf{C}_C$ we denote the class of simplicial complexes $T_{\EE}(V)$ where $n=\dim(T_{\EE}(V) +1 \geq 1$, $K$ a (skew) field, $(V,Q,f)$ is a thick pseudo-quadratic space with Witt index $n$, $\EE$ a finite subset of $\mathcal{W}$ (see \Cref{def: C_n case N(E)}) such that $\EE^\perp = \EE$. If $K$ is a finite field, we further require that $\lvert K \rvert \geq N(\EE)$. We set $T_{\EE}(V)$ as in \Cref{def: X_E in C_n case}.
\end{definition}

\subsubsection{Buildings of type $D_n$} \label{subsubsec: geom descr D_n}
For this section let $(V,Q,f)$ be a pseudo-quadratic space but in the specific case where $K$ is a field, $\sigma=\operatorname{id}_K, \epsilon = 1, \Lambda = 0$, $V$ a $K$-vector space of dimension $m=2n \geq 4$, $Q$ an ordinary quadratic form and $f$ the corresponding symmetric bilinear form. We further assume that $V$ is non-degenerate and of Witt index $n$. 

We set as before $X = X(V)=\{0<U<V \mid U \text{ is totally isotropic}\}$. Then $\Delta = \flag X$ is a weak $C_n$ building, in particular it is not thick, since each totally isotropic subspace of dimension $n-1$ is contained in exactly two totally isotropic subspaces of dimension $n$.  Furthermore, the space $Y= \{U\in X \mid \dim U = n\}$ is partitioned into two sets $Y = Y_1\sqcup Y_2$ such that for $U_1,U_2 \in Y_i$ we have $n-\dim(U_1 \cap U_2)\in 2\Z$ for $i=1,2$ and for $U_1 \in Y_1, U_2 \in Y_2$ we have $n - \dim(U_1 \cap U_2) \in 2\Z +1$. In particular, if for $U_1,U_2 \in Y$ we have $\dim(U_1 \cap U_2) = n-1$ then $U_1$ and $U_2$ are in different sets of the partition. 

\begin{example}
	We illustrate the situation with the following example for $n=2, m=4$ (i.e. $V$ is a 4-dimensional vector space over the field $K$) and bilinear form $f$ given by $f(x,y) = x_1y_2 + x_2y_1+ x_3y_4 +x_4y_3$. Then the line generated by $e_1 = (1,0,0,0)$ is totally isotropic (since $f(e_1,e_1)=0$). Let $y = (y_1,y_2,y_3,y_4)$ be another vector and assume that the plane spanned by $e_1$ and $y$ is totally isotropic. This gives the conditions $f(e_1,y)=f(y,y)=0$. The first one implies $y_2=0$ and the latter $y_3y_4=0$. Hence there are exactly two isotropic planes $\langle e_1,e_3\rangle$ and $\langle e_1,e_4\rangle$ containing $\langle e_1 \rangle$.
\end{example}

To obtain a thick building of type $D_n$ we set $\tilde{X}=\tilde{X}(V) = \{U \in X \mid \dim U \neq n-1\}$ and we say that $U,W \in \tilde{X	}$ are incident (write $UIW$) if and only if
$$ U \subseteq W \text{ or } W \subseteq U \text{ or } \dim (W \cap U) = n-1.$$
Then $\tilde{\Delta} = \flag \tilde{X}$ is desired building. We will write $\orifl(\tilde{X})$ for $\flag \tilde{X}$ to stress that we consider $\tilde{X}$ not with the usual incidence relation given by inclusion, but with this new one. Note that $\orifl \tilde{X}$ is called the oriflamme complex of $\tilde{X}$.  

For $n\geq 4$ every building of type $D_n$ is of the form $\tilde{\Delta}$ for some field $K$ \cite[Proposition 8.4.3]{TitsBuildings}. We also consider the construction for $n=2,3$ in which case we get certain (but not all) buildings of type $D_2 = A_1 \times A_1$ and $D_3 = A_3$, which were already covered earlier.

To describe the opposite complexes, we need to introduce some further notation. 
\begin{definition} 
For arbitrary subspaces $U, W$ of $V$, we define
$$U\mathbin{\tilde{\pitchfork}}_V W \iff U \pitchfork_V W \text{ or } (U,W \in \tilde{X}, U = U^\perp, W = W^\perp \text{ and } \dim(U \cap W) = 1).$$
For a set $\EE$ of subsets of $V$ we define (different from \Cref{def: A_n case X_E})
\begin{align*}
X_{\EE}(V) &= \{U \in X \mid U \mathbin{\tilde{\pitchfork}} \EE \}, \qquad T_{\EE}(V) = \flag X_{\EE}(V) \\
\tilde{X}_{\EE}(V) &= \{U \in \tilde{X} \mid U \mathbin{\tilde{\pitchfork}} \EE \}, \qquad \tilde{T}_{\EE}(V) = \orifl \tilde{X}_{\EE}(V) 
\end{align*}
For $T_{\EE}(V)$ we consider the usual inclusion of subspaces as incidence relation, while for $\tilde{T}_{\EE}(V)$ we consider the incidence relation specified above.
\end{definition}

\begin{proposition}{\cite[Corollary 17]{Abr}} \label{prop: description oppo Dn}
Let $\tau = \{E_1<\dots < E_r\} \in \tilde{\Delta}$ and set $\EE(\tau) = \{E_i, E_i^\perp \mid 1 \leq i\leq r\}$ then $\tilde{\Delta}^0(\tau) = \tilde{T}_{\EE(\tau)}(V)$.
\end{proposition}

\section{Coboundary expansion for sub-complexes of spherical buildings} 
\label{sec: cbe for subcomplexes of spherical buildings}
The goal of this section is to show that the opposite complex of an irreducible spherical building is a coboundary expander by constructing a cone function. We achieve this for the opposite complexes of type $A_n,C_n$ and $D_n$ using their geometric description from \Cref{subsec: geometric construction of buildings} (recall that for Coxeter complexes and hence buildings there is no difference between type $B_n$ and type $C_n$). Using that we can construct cone functions for joins given we have cone functions for the factors, we cover all opposite complexes of classical type, i.e. where each connected component of the Dynkin diagram is of type $A_n,B_n, C_n$ or $D_n$. We start by generalizing an induction argument from \cite[Lemma 22]{Abr} to our setting of cone functions. It gives a list of conditions on a class of simplicial complexes $\mathbf{C}$ that can be checked by induction on the dimension of the complexes which results in the existence of cone functions with bounded radius for all complexes in the class. We then show that this theorem can be applied to the classes $\mathbf{C}_A, \mathbf{C}_A \cup \mathbf{C}_C$ from \Cref{def: class C_A,def: class C_C} which shows the existence of cone functions for opposite complexes of type $A_n$ and $C_n$. The argument for type $D_n$ is more involved and we will discuss it in the corresponding section.

\subsection{General induction argument}

The general procedure will be similar for all three cases of buildings. We use a quantitative version of \cite[Lemma 22]{Abr}. We basically only change assumption 1. and the conclusion. 

\begin{theorem} \label{thm: generalization of Lemma 22}
Let $\mathbf{C}$ be a class of non-empty simplicial complexes. 
Assume that for any $n \in \mathbb{N}, n\geq 1$ there exists $\ell_{n} \in \mathbb{N}$ such that for all $\kappa \in \mathbf{C}, \dim \kappa = n$, there exists a filtration $\kappa_{0} \subset \kappa_{1} \subset \dots \subset \kappa_{\ell} = \kappa$ of $\kappa$ of length $\ell \leq \ell_n$ such that the following is satisfied (set $V_{i}:= \kappa_{i}(0) \setminus \kappa_{i-1}(0) $):
\begin{enumerate}
\item There is an $(n-1)$-cone function $\Cone_{\kappa_{0}}$ such that for every $0\leq j \leq n-1$, 
 $$\Rad_{j}(\Cone_{\kappa_{0}}) \leq f(n),$$
	where $f: \mathbb{R} \to \mathbb{R}$ is a function depending only on $\mathbf{C}$. \label{gen L22: condition kappa0}
\item $\kappa_{i}$ is a full subcomplex of $\kappa$ and for all vertices $w_1,w_2 \in V_{i}$ we have $\{w_1,w_2\}\notin \kappa$, for all $0 \leq i \leq \ell$.
\label{gen L22: condition full and dif type}
\item  $\lk_{\kappa_{i}}(w) \cap \kappa_{i-1}$ is either in $\mathbf{C}$ or can be written as the join of two elements from $\mathbf{C}$ for any $w \in V_{i}, i \geq 1$. \label{gen L22: condition links}
\end{enumerate}
Then for every $n \in \mathbb{N}_{0}$ there exists a constant $\mathcal{R}(n)$ such that for every $\kappa \in \mathbf{C}$ with $\dim \kappa = n$ we have:
$\kappa$ has an $(n-1)$-cone function $\Cone_{\kappa}$ such that for all $-1 \leq j \leq n-1$ we  have
$$\Rad_{j}(\Cone_{\kappa}) \leq \mathcal{R}(n).$$
We can describe $\mathcal{R}(n)$ recursively. We have $\mathcal{R}(0) = 1$. For $n\geq 1$, assume $\mathcal{R}(k)$ is known for $k <n$. Set 
$$S(n) = \max_{a,b \in \N: a+b+1=n} ((a+1)\mathcal{R}(a) +1)\mathcal{R}(b).$$
Then 
$$\mathcal{R}(n) = S(n)^{\ell_n} f(n) + \sum_{j=1}^{\ell_n} S(n)^j.$$
\end{theorem}

\begin{remark}
	Note that the theorem does not assume that the $\kappa_i$ are in $\mathbf{C}$. Furthermore, we have that 
	\[\lk_{\kappa_i}(w) \cap \kappa_{i-1} = \lk_{\kappa}(w)\cap \kappa_{i-1}.\]
	The $\subseteq$ inclusion is clear since $\kappa_i \subseteq \kappa$. The other inclusion $\supseteq$ holds since $\kappa_{i}$ is a full subcomplex of $\kappa$ thus for $\sigma \in \kappa_{i-1}\subseteq \kappa_i$ and $w \in \kappa_i$ if $\sigma \cup w \in \kappa$ we can deduce $\sigma \cup w \in \kappa_i$. 
\end{remark}

\begin{proof}
We will prove the result by induction on the dimension $n$.

For $n=0$, since $\kappa \neq \emptyset$ there exists $v \in \kappa(0)$ and thus we can define a $(-1)$-cone function via $\emptyset \mapsto \II_{[v]}$. Every cone function has $(-1)$-cone radius = 1, thus we take $\mathcal{R}(0)=1$. 

Fix $n\geq 1$ and $\kappa \in \mathbf{C}, \dim(\kappa)=n$ with filtration $\kappa_{0} \subset \dots \subset \kappa_{\ell}=\kappa$ as in the requirements.
We will prove by induction on $0 \leq i \leq \ell$, that there is an $(n-1)$-cone function $\text{Cone}_{\kappa_{i}}$ such that for every $-1 \leq k \leq n-1$
$$\text{Rad}_{k}(\text{Cone}_{\kappa_{i}}) \leq \mathcal{R}^{(i)}(n) $$
where $\mathcal{R}^{(i)}(n)$ is given inductively by
$$\mathcal{R}^{(0)}(n) = f(n), \qquad \mathcal{R}^{(i)}(n)= S(n)(\mathcal{R}^{(i-1)}(n)+1)=S(n)^if(n) + \sum_{j=1}^i S(n)^j.$$

For $i=0$ this follows from Assumption \ref{gen L22: condition kappa0}.
We proceed by induction on $i$. Fix $1 \leq i \leq \ell$ and assume there exists a constant $\mathcal{R}^{(i-1)}(n)$ and an $(n-1)$-cone function $\text{Cone}_{\kappa_{i-1}}$ such that $\text{Rad}_{j}(\text{Cone}_{\kappa_{i-1}}) \leq \mathcal{R}^{(i-1)}(n)$ for all $-1 \leq j \leq n-1$.
We want to apply \Cref{adding vertices to cone thm} to the following setting:

\begin{align*}
X=\kappa, X' = \kappa_{i-1},\  \operatorname{Cone}_{X'}= \operatorname{Cone}_{\kappa_{i-1}}, \ R_{j}' = \mathcal{R}^{(i-1)}(n), \ W = V_{i} = \kappa_{i}(0) \setminus \kappa_{i-1}(0).
\end{align*}
We check that the conditions of the theorem are satisfied.
\begin{enumerate}
\item $W\cap X' = \emptyset$ by definition.
\item For every $w_{1},w_{2} \in W: \left\{ w_{1},w_{2} \right\} \not\in \kappa_{i} \subset \kappa$ by assumption.
\item For every $w \in W, \lk_{X}(w) \cap X'$ is a non-empty simplicial complex, since by Assumption \ref{gen L22: condition links} we have $\lk_{X}(w) \cap X' \in \mathbf{C}$ or a join of two complexes in $\mathbf{C}$.  
\item By Assumption \ref{gen L22: condition links}, $\lk_{X}(w) \cap X'$ is either in $\textbf{C}$ or is a join of two elements of $\textbf{C}$. In the first case we use the induction hypothesis and get an $(n-2)$-cone function $\Cone_{\lk_{X}(w) \cap X'}$ such that 
$$\text{Rad}_{j}(\operatorname{Cone}_{\lk_{X}(w)\cap X'}) \leq \mathcal{R}(n-1).$$
In the second case, we can write $\lk_{X}(w) \cap X' = A *B$ for some $A,B \in \textbf{C}$. Set $n_A=\dim A, n_B= \dim B$ then $n-1 \geq \dim(\lk_{X}(w) \cap X') = n_A+n_B+1$, hence $n_A,n_B<n$. Thus by induction there exists an $(n_A - 1)$-cone function $\Cone_A$ with $\Rad_j(\Cone_A) \leq \mathcal{R}(n_A)$ for all $-1\leq j \leq n_A-1$ and similarly for $B$. Using \Cref{cone of a join - general prop} we get an $(n-1)$-cone function $\Cone_{\lk_{X}(w) \cap X'}$ such that 
$$\Rad_j(\Cone_{\lk_{X}(w) \cap X'}) \leq ((n_A+1)\mathcal{R}(n_A)+1)\mathcal{R}(n_B) \leq S(n) \text{ for all } -1 \leq j \leq n-2.$$
Note that we simplified the bound coming from \Cref{cone of a join - general prop} by only considering the largest possible of the cases. We chose $S(n)$ to be the maximum over all the bounds that can appear at this step. Note that in particular $\mathcal{R}(n-1) \leq S(n)$. We set $R_j'' = S(n)$.
\end{enumerate}

Applying \Cref{adding vertices to cone thm} we get an $(n-1)$-cone function $\operatorname{Cone}_{X' \cup W}$ such that 
$$\text{Rad}_{j}(\operatorname{Cone}_{X' \cup W}) \leq R_{j-1}''(R_{j}'+1)$$
which in the original formulation means
$$\text{Rad}_{j}(\operatorname{Cone}_{\kappa_{i}}) \leq S(n)(\mathcal{R}^{(i-1)}(n)+1).$$
Thus we have $\mathcal{R}^{(i)}(n) = S(n)(\mathcal{R}^{(i-1)}(n)+1)$ and solving the recursion gives the second expression of $\mathcal{R}^{(i)}(n)$. 

Now, since $\kappa_{\ell}=\kappa$, we get an $(n-1)$-cone function $\Cone_\kappa$ with $\Rad(\Cone_\kappa) \leq \mathcal{R}^{(\ell)}(n)$. We would like to set $\mathcal{R}(n)= \mathcal{R}^{(\ell)}(n)$. But the bound has to hold for all possible $\kappa \in \mathbf{C}$ of dimension $n$. Since $\mathcal{R}^{(i)}(n)$ is increasing in $i$, we take $\mathcal{R}(n) = \mathcal{R}^{(\ell_n)}(n) = S(n)^{\ell_{n}} f(n) + \sum_{j=1}^{\ell_{n}} S(n)^{j}$.
\end{proof}

\begin{remark}
Compared to \cite[Lemma 22]{Abr} we changed condition 1. from asking for $\kappa_0$ to be contractible to requiring an explicit cone function, and we removed the dependency on the underlying building, since the theorem also works in this wider generality. In particular, if we want to show that a certain class, for which Abramenko already showed that it satisfies the assumptions of \cite[Lemma 22]{Abr}, also satisfies the assumptions of this theorem, we need to check condition 1., while conditions 2. and 3. are usually already covered.
\end{remark}

We will apply this theorem to classes of simplicial complexes which are flag complexes over sets of vector subspaces. In this case, condition \ref{gen L22: condition full and dif type} can be achieved by only adding vector subspaces with the same dimension in one step of the filtration. By not altering the incidence relation considered for the flag complex, we can make sure that the $\kappa_i$ are full subcomplexes. 

To show that condition \ref{gen L22: condition links} is satisfied in our class of complexes, we use that for flag complexes over vector subspaces we can describe the links more explicitly using the following definition. 
\begin{definition}
Let $K$ be a field, $V$ a $K$-vector space, $X$ a collection of subspaces of $V$. Let $Y = \flag X$ and $U\in X$. Then we define
\begin{align*}
 X^{<U}&:=\{W\in X \mid W < U \}, \ Y^{<U} := \flag X^{<U}\\
 X^{>U}&:=\{W\in X \mid W > U \}, \ Y^{>U} := \flag X^{>U}. 
\end{align*}
\end{definition}
\begin{lemma} \label{lemma: link in flag complex}
Note that
$$\lk_Y(U) = Y^{<U} * Y^{>U}.$$
Moreover, if $Y_i = \flag X_i$ for $i=1,2$ such that $Y_1 \subseteq Y_2$ and $U \in X_2$, then 
\[\lk_{Y_2}(U) \cap Y_1 = Y_1^{<U} * Y_1^{>U}\]
\end{lemma}
\begin{proof}
	Straightforward application of the definitions. 
\end{proof}
Note that a similar statement holds for flag complexes over more general posets.

\subsection{Buildings of type $A_n$}
\label{An subcomplexes sec}

To show that the complex opposite the fundamental chamber in a spherical building of type $A_n$ has a cone function with bounded cone radius, we show that the class $\mathbf{C}_A$ from \Cref{def: class C_A} satisfies the assumptions of \Cref{thm: generalization of Lemma 22}.

We follow the machinery presented in \cite{Abr} closely, which itself builds on previous work \cite{Quill,Vogt,AApaper}.

We proceed as follows
  \begin{enumerate}
    \item Each complex $\kappa$ in the classes $\mathbf{C}_A$ is the flag complex of a certain set of vector subspaces. We recall the fact from \cite{Abr} that at least one of these subspaces has dimension 1. 
    \item We use this subspace of dimension 1 to define a filtration $\kappa_0 \subset \dots \subset \kappa_\ell = \kappa$ as in \cite{Abr}.
    \item We show that $\kappa_0$ admits a cone function and we bound its radius. \label{rmk proof structure part cone kappa0}
    \item We recall from the work in \cite{Abr} that $\lk_{\kappa_i}(U) \cap \kappa_{i-1}$ for all $U \in \kappa_i(0) \setminus \kappa_{i-1}(0)$ has the desired structure, i.e. is in the class we consider or can be written as a join of two elements in there.
    \item We conclude that the class satisfies the required assumptions and hence get a cone function with bounded radius for each element, in particular for the opposite complexes. 
  \end{enumerate}

We recall some facts that appear in the proof of \cite[Step 1, Proposition 12]{Abr}, see also \cite{AApaper}.  
\begin{proposition}
\label{Abr fact1}
If $T_{\EE}(V) \in \mathbf{C}_A$, then there exists $\ell \in X_{\mathcal{E}} (V)$ such that $\ell$ is a line,  i.e.,  $\dim (\ell) =1$. 
\end{proposition}
\begin{proof}
	Since similar statements will appear for each case, we will give a proof of this statement here. 

	A 1-dimensional subspace $\ell<V$ will be transversal to an arbitrary proper subspace $E<V$ if and only if $\ell \not \subseteq E$. Let $\EE$ be a collection of subspaces such that $T_{\EE}(V) \in \mathbf{C}_A$. In particular, $\EE$ contains finitely many subspaces. Thus, if the field $K$ is infinite, there will always exist a 1-dimensional subspace $\ell \not \subseteq \bigcup_{E \in \EE}E$. Now, if $K$ is a finite field of size $q$, then each $E$ contains $\frac{q^{\dim E}-1}{q-1}$ lines. Thus in total $\bigcup_{E \in \EE}E$ contains at most 
	\[\sum_{E \in \EE} \frac{q^{\dim E}-1}{q-1} \leq \frac{q^{n}-1}{q-1} \sum_{i=1}^n e_i \leq q \frac{q^n-1}{q-1} < \frac{q^{n+1}-1}{q-1} \]
	where $e_i = \vert \{E\in \EE \mid \dim E = i\} \vert$ as above.
	Thus there exists at least one line transversal to all of $\EE$. 
\end{proof}

The important property of $\ell$ being a line, that we will use, is that for any other subspace $W$ we have either $\ell \leq W$ or $\ell \cap W = 0$ in which case $W < W+\ell$.

Let $\ell \in X_{\mathcal{E}} (V)$ be a line.  Define 
$$Y_0 = \lbrace  U  \in X_{\mathcal{E}} (V) : \ell + U \in X_{\mathcal{E}} (V) \rbrace,$$
and further define $\kappa_0$ to be the subcomplex of $T_{\mathcal{E}} (V)$ spanned by $Y_0$.

Assume now that $\dim V = n+2$ and hence $\dim T_{\EE}(V) = n$ (and $\flag X$ is a building of type $A_{n+1}$). 
For $1 \leq i \leq n+1$,  define $Y_i \subseteq X_{\mathcal{E}} (V)$ as follows: 
$$Y_i = \lbrace U \in X_{\mathcal{E}} (V)\mid \dim (U) \geq n+2-i \text{ or } U \in Y_0 \rbrace.$$
For every such $i$,  let $\kappa_i$ be the subcomplex of $T_{\mathcal{E}} (V)$ spanned by $Y_i$. We have $\kappa_{n+1} = T_{\EE}(V)$.

\begin{lemma}
\label{X0 bound lemma}
There is an $(n-1)$-cone function $\Cone_{\kappa_0}$ such that for every $0 \leq k \leq n-1$,  
$$\Rad_{k} (\Cone_{\kappa_0}) \leq n+2 .$$
\end{lemma}

\begin{proof}
Let $\Gamma_0$ be the subcomplex of $\kappa_0$ spanned by $Z_0=\lbrace  U  \in X_{\mathcal{E}} (V)  : \ell \leq U \rbrace$.  For $1 \leq i \leq n+1$,  we denote 
$$Z_i = \lbrace U   \in Y_0 \mid   \dim (U) \leq i \text{ or } U \in Z_0 \rbrace.$$
Further denote $\Gamma_i = \flag Z_i$. Note that $\Gamma_{n+1} = \kappa_0$.  We will show by induction that for every $0 \leq i \leq n+1$, there is an  $(n-1)$-cone function $\Cone_{\Gamma_i}$ such that for every $-1 \leq k \leq n-1$,  $\Rad_k ( \Cone_{\Gamma_i}) \leq i+1$.

For $i=0$,  we note that  $\Gamma_0$ is the ball of radius $1$ around $\ell$ in $\kappa$ and can be written as $\Gamma_0 = \lbrace \ell \rbrace * \lk_{\kappa}(\ell)$.  By Proposition \ref{cone of a join - single vertex prop},  there is an $(n-1)$-cone function $\Cone_{\Gamma_0}$ such that for every $-1 \leq k \leq n-1$,  $\Rad_k (\Cone_{\Gamma_0}) \leq 1$.  

Let $1 \leq i \leq n+1$ and assume there is an $(n-1)$-cone function $\Cone_{\Gamma_{i-1}}$ such that for every $-1 \leq k \leq n-1$,  $\Rad_k ( \Cone_{\Gamma_{i-1}}) \leq i$.  We will apply Theorem \ref{adding vertices to cone thm} in order to bound the cone radii of $\Gamma_i$.  We note that any two $U_1,  U_2 \in Z_i \setminus Z_{i-1}, U_1 \neq U_2$ are of the same dimension and thus are not connected by an edge.  Let $U \in Z_i \setminus Z_{i-1}$.  We will show that $\lk_{\kappa_0}(U) \cap \Gamma_{i-1}$ is a join of the vertex $\lbrace U + \ell \rbrace$ with the complex spanned by all the other vertices in $\lk_{\kappa_0}(U) \cap \Gamma_{i-1}$. Note that $U+\ell \neq U$ since $U \notin Z_0$, and $U +\ell \in Z_0 \subseteq Z_{i-1}$. 

Since $\lk_{\kappa_0}(U) \cap \Gamma_{i-1}$ is a clique complex, it is enough to show that for every vertex $U'$ in $\lk_{\kappa_0}(U) \cap \Gamma_{i-1}$ that is not $\ell + U$,  there is an edge connecting $U'$ and $\ell + U$.  Let $U'$ be such a vertex.  

First, we will deal with the case where $\dim (U') < i$.  In that case (using the fact that $\dim (U) =i$),  $U'$ being in $\lk_{\kappa_0}(U)$ implies that $U' \leq U$ and thus $U' \leq U + \ell$,  i.e., $U'$ is connected by an edge to $U + \ell$ as needed.  

Second, assume that $\dim (U') >i$.  From the definition of $Z_{i-1}$,  it follows that $U' \in Z_0$,  i.e., that $\ell \leq U'$.  Also,  $U'$ is in $\lk_{\kappa_0}(U)$ and thus $U \leq U'$.  It follows that $U+\ell \leq U'$,  i.e.,  $U'$ is connected by an edge to $U + \ell$ as needed.  

By Proposition \ref{cone of a join - single vertex prop},  for every $U \in Z_i \setminus Z_{i-1}$,  there is an $(n-2)$-cone function $\Cone_{\lk_{\kappa_0}(U) \cap \Gamma_{i-1}}$ such that for every $-1 \leq k \leq n-2$,  $\Rad_k (\Cone_{\lk_{\kappa_0}(U) \cap \Gamma_{i-1}}) \leq 1$.  

By Theorem \ref{adding vertices to cone thm} (with $R_k '' =1$ for every $k$),  it follows that there is an $(n-1)$-cone function  $\Cone_{\Gamma_i}$ such that for every $-1 \leq k \leq n-1$,  
$\Rad_{k} (\Cone_{\Gamma_i}) \leq i+1$ as needed. 
\end{proof}

\begin{fact}[{\cite[Step 3, Proposition 12]{Abr}}]\label{fact: abr 3}
	Let $Y = X_{\mathcal{E}} (V), \kappa = T_{\EE}(V) \in \mathbf{C}_A$ and $\kappa_0 \subseteq \dots \subseteq \kappa_{n+1}=\kappa$ with $\kappa_i = \flag Y_i$ be the filtration as above. For every $1 \leq i \leq n+1$ and every $U \in Y_i \setminus Y_{i-1}$, we have
	\begin{enumerate}
		\item $\lk_{\kappa_i}(U) \cap \kappa_{i-1} = Y_{0}^{<U}* Y^{>U}$ and \label{fact: point 1}
		\item $\flag Y_{0}^{<U}, \flag Y^{>U} \in \mathbf{C}_A$. \label{fact: point 2}
	\end{enumerate}
\end{fact}
\begin{proof}
	We sketch the idea of the proof. For part \ref{fact: point 1}, we know from \Cref{lemma: link in flag complex} that $\lk_{\kappa_i}(U) \cap \kappa_{i-1} = Y_{i-1}^{<U}* Y_{i-1}^{>U}$. We have $Y_{i-1}^{<U}=Y_0^{<U}$ since $\dim U = n+2-i < n+2-i+1$. On the other hand, any $W>U$ has $\dim W > n+2-i$ and thus $Y_{i-1}^{>U} = Y^{>U}$.

	For part \ref{fact: point 2}, we observe
	\begin{itemize}
		\item if $W<U$ then $W\pitchfork_V \EE$ and $W+ \ell \pitchfork_V \EE$ if and only if $W \pitchfork_U \EE'$ with $\EE' = (\EE \cap U)\cup ((\EE + \ell) \cap U)$.
		\item If $U<W<V$ then $W\pitchfork_V \EE$ if and only if $W+U \pitchfork_{V/U}\bar{\EE}$ with $\bar{\EE} = \{E+U \mid E \in \EE\}$. 
	\end{itemize}
	Thus $Y_0^{<U} \cong X_{\EE'}(U)$ and $Y^{>U} \cong X_{\bar{\EE}}(V/U)$. In Step 3 \cite[Proposition 12]{Abr} it is shown that $X_{\EE'}(U), X_{\bar{\EE}}(V/U) \in \mathbf{C}_A$.
\end{proof}

\begin{corollary}
	The class $\mathbf{C}_A$ satisfies the assumptions of \Cref{thm: generalization of Lemma 22}.
\end{corollary}
\begin{proof}
The result follows from \Cref{X0 bound lemma} and Fact \ref{fact: abr 3}. In the notation of \Cref{thm: generalization of Lemma 22}, we have $f(n) = n+2$ and $\ell_n = n+1$.
\end{proof}

Using \Cref{prop: An oppo as flag} we can then get the desired conclusion for opposite complexes of type $A_n$.
\begin{corollary} \label{thm: cone for A_n case}
Let $\Delta$ be a building of type $A_n$ over a field $K$ with $\lvert K \rvert \geq 2^{n-1}$ and $a\in \Delta$ a simplex, then there exists an $(n-2)$-cone function $\Cone_{\Delta^0(a)}$ of $\Delta^0(a)$ such that  
$$\Rad_{j} (\Cone_{\Delta^0(a)}) \leq \mathcal{R} (n-1),  \forall -1 \leq j \leq n-2 .$$
\end{corollary}

\subsection{Buildings of type $C_n$}

The class $\textbf{C}$ that will cover the case of opposite complexes in buildings of type $C_n$ will be the union of class $\mathbf{C}_A$ (see \Cref{def: class C_A}) and $\mathbf{C}_C$ (see \Cref{def: class C_C}).

We start by considering $\kappa = T_{\EE}(V)=\flag X_{\EE}(V) \in \mathbf{C}_C$ of dimension $n \geq 0$ with $V,\EE$ as in \Cref{def: class C_C}. The corresponding building is then of type $C_{n+1}$.

We recall from \cite[Proposition 13, Step 1]{Abr} the following.
\begin{fact}
  There exists a 1-dimensional subspace $\El \in X_{\mathcal{E}}(V)$.
\end{fact}
We will also need the following fact, which can be found in \cite[Proposition 13, Step 2]{Abr}.
\begin{fact}\label{fact: pitchfork plus ell}
	Given $\EE$ and $\El$ as above, for $0<U<V$, if $\El \nleq U$ then we have $U\pitchfork \EE +\El$ if and only if $U + \El \pitchfork \EE$.
	Furthermore, if $W,E \nleq \El^\perp$ then $W \cap \El^\perp \pitchfork E \iff W \pitchfork E \cap \El^\perp$.
\end{fact}

We define the filtration $\kappa_0 \subset \dots \subset \kappa_{2n+2}$ of $\kappa$ as follows. We start by setting $X_{0}=X^1 \cup X^2 \cup X^3$ where 
\begin{enumerate}
\item $X^1 = \{U \in X \mid \El \leq U$ and $U \pitchfork \mathcal{E}\}$; 
\item $X^2 = \{U \in X \mid \El \nleq U \leq \El^\perp$ and $U \pitchfork \mathcal{E} \cup (\mathcal{E}+\El)\}$;
\item $X^3 = \{U \in X \mid U\nleq \El^\perp, \dim U >1$ and $U \pitchfork \mathcal{E}\cup(\mathcal{E} \cap \El^\perp) \cup (\mathcal{E} + \El) \cup ((\mathcal{E} \cap \El^\perp)+\El)\}$.
\end{enumerate}
Set $\kappa_0 = \flag X_0$.

Next, we define the rest of the desired filtration in two steps. Set $Z=\{U \in X_{\mathcal{E}}(V) \mid \El \leq U \text{ or } U \pitchfork \mathcal{E}+\El \}$ and $Y_0 = X_0$. Note that $Y_0 \subseteq Z$. For the first half of the filtration we define 
$$Y_i = \{U \in Z \mid U \in Y_0 \text{ or } \dim U \leq i \}, \kappa_i = \flag Y_i, \ 1 \leq i \leq n+1.$$
For the second half we set
$$Y_i = \{U \in X_{\mathcal{E}}(V) \mid U \in Z \text{ or } \dim U \geq 2n+3-i \}, \kappa_i = \flag Y_i, \ n\leq i \leq 2n+2.$$
Then $T_{\EE}(V) = \kappa_{2n+2}$. We start by showing the existence of a cone function for $\kappa_0$, where we use the filtration of $\kappa_0$ from \cite[Proposition 13]{Abr}.

\begin{lemma} \label{lemma: abelian kappa_0 C_n case}
With the definitions as above, we have that $\kappa_{0}$ has an $(n-1)$-cone function $\Cone_{\kappa_0}$ such that 
$$\Rad_j(\Cone_{\kappa_0}) \leq 2n+3 \qquad \forall -1\leq j \leq n-1.$$
\end{lemma}
\begin{proof}
Let $T_{\mathcal{E}}(V) \in \textbf{C}_C$ as above. Thus $\dim \kappa = n$. In particular, any totally isotropic subspace of $V$ has dimension at most $n+1$. 

We will apply \Cref{adding vertices to cone thm} inductively to two different filtrations of $\kappa_0$. The first one is defined as follows, using $X^1,X^2,X^3$ as defined above.
\begin{align*}
A_0 &= X^1 \\
A_i &= Y_{i-1} \cup \left\{U \in X^2 \mid \dim U = i \right\} = X^1 \cup \{U\in X^2 \mid \dim U \leq i \} \text{ for } 1 \leq i \leq n+1.
\end{align*} 
We set $\alpha_i = \flag A_i$. Note that $A_{n+1} = X^1 \cup X^2$. 

The second filtration is given by
\begin{align*}
B_0 &= A_{n+1} \\
B_i &= B_{i-1} \cup \left\{ U \in X^3 \mid \dim U = n+1-(i-1) \right\} \\
&= B_{0} \cup \left\{ U \in X^3 \mid \dim U \geq n+1-(i-1) \right\} \text{ for } 1 \leq i \leq n+1. 
\end{align*}
We set $\beta_i = \flag B_i$. Note that $B_{n+1} = X^1 \cup X^2 \cup X^3 = X_0$ and hence $\beta_{n+1}=\kappa_0$.

Next, we want to show that for each $\alpha_i$ there exists an $(n-1)$-cone function $\Cone_{\alpha_i}$ such that for every $-1\leq j \leq n-1$ we have $\Rad_j(\Cone_{\alpha_i}) \leq i+1$ for all $0 \leq i \leq n+1$. We do this by induction on $i$. 

For $i=0$ we have $\alpha_0 = \st_{\kappa}(\El) = \{\El\} * \operatorname{lk}_\kappa(\El)$. Hence by \Cref{cone of a join - single vertex prop} there exists an $(n-1)$-cone function $\Cone_{\alpha_0}$ with $\Rad_j(\Cone_{\alpha_0}) \leq 1$ for all $-1\leq j \leq n-1$.

For the following step, fix $1 \leq i \leq n+1$. We want to apply \Cref{adding vertices to cone thm} to the following set-up. 
Here $\kappa_{0}$ corresponds to the complex called $X$ in \Cref{adding vertices to cone thm}, $\alpha_{i-1}$ corresponds to $X'$ , and $\cW_i :=  A_{i} \setminus A_{i-1} $ to $W$ i.e. the set of vertices that we want to add. 
We check that the conditions of \Cref{adding vertices to cone thm} are satisfied.

For assumption \ref{add vert thm: as 1} note that $\cW_{i} \cap \alpha_{i-1}=\emptyset$ by definition.
Assumption \ref{add vert thm: as 2} holds since all subspaces in $\cW_{i}$ have the same dimension, hence they are not connected in $\kappa$.

For \ref{add vert thm: as 3} and \ref{add vert thm: as 4}, let $W\in \cW_{i}$. Then $\lk_{\kappa_{0}}(W) = \flag\left\{ U \in \kappa_{0}(0) \mid U < W \text{ or } W > U  \right\}$. Furthermore, for $U \in A_{i-1}$ we have $U+\El \in X_{0}$, since by Fact \ref{fact: pitchfork plus ell} we have $U \pitchfork \mathcal{E} + \El \iff U + \El \pitchfork \mathcal{E}$ if $\El \nleq U$. If $\El \leq U$ then we have $U \in X^1 \subseteq X_0$.  
Note that $\El \nleq W$ because $W \in X^2$. Furthermore, since $\El \leq W+\El$, we have $W + \El \in A_{0} \subset A_{i-1}$. Next, we want to show that $\lk_{\kappa_{0}}(W) \cap \alpha_{i-1}$ is a join of $\{W + \El\}$ with the flag complex of all other vertices from $\lk_{\kappa_{0}}(W) \cap \alpha_{i-1}$. Since $\lk_{\kappa_{0}}(W) \cap \alpha_{i-1}$ is a flag complex, it suffices to check that every vertex different from $W + \El$ is adjacent to $W + \El$. Let $U \in (\lk_{\kappa_{0}}(W) \cap \alpha_{i-1})(0) \setminus \{W + \El\}$. Then we differentiate two cases: 
\begin{itemize}
\item Case 1: $\dim U \leq i-1$: then $U \leq W \leq W + \El$;
\item Case 2: $\dim U >i$ then $U \in A_{0}$, thus $W \leq U$ and $\El \leq U$ and hence $W + \El \leq U$.
\end{itemize}
Hence, any other vertex is connected to $W+ \El$ by an edge and thus $$\lk_{\kappa_{0}}(W) \cap \alpha_{i-1} = \{W + \El \} *\flag  \left((\lk_{\kappa_{0}}(W) \cap \alpha_{i-1})(0) \setminus \{W + \El\}\right).$$ By \Cref{cone of a join - single vertex prop} we have
$\operatorname{Rad}_j(\operatorname{Cone}_{\lk_{\kappa_{0}}(W) \cap \alpha_{i-1}}) \leq 1$ for all $-1 \leq j \leq n-1$. 

Thus we can apply \Cref{adding vertices to cone thm} and get an $(n-1)$-cone function $\Cone_{\alpha_i}$ satisfying $\operatorname{Rad}_j(\operatorname{Cone}_{\alpha_{i}}) \leq i+1$.

We treat the second filtration similarly. We want to show that for every $0 \leq i \leq n+1$ there is an $(n-1)$-cone function $\Cone_{\beta_i}$ such that $\Rad_j(\Cone_{\beta_i}) \leq n+2 +i$. We again proceed by induction on $i$. 

For $i=0$ we have that $\beta_{0}=\alpha_{n+1}$ has a cone function with radius $\leq n+2$ by the above reasoning.

Now fix $1\leq i \leq n+1$. We want to use \Cref{adding vertices to cone thm} with 
$\kappa_{0}$ corresponding to $X$, $\beta_{i-1}$ to $X'$, and $B_{i}\setminus B_{i-1}$ to $W$.
Then assumptions \ref{add vert thm: as 1} and \ref{add vert thm: as 2} of \Cref{adding vertices to cone thm} are satisfied by the same argument as above. 

For assumptions \ref{add vert thm: as 3} and \ref{add vert thm: as 4}, fix $W \in B_i \setminus B_{i-1}$ and consider 
$$\lk_{\kappa_{0}}(W) \cap \beta_{i-1} = \flag \left\{ U \in B_{i-1} \mid U <W \text{ or } U>W \right\}.$$

First, note that $W \in X^3$, thus $\El^\perp \nleq W$ and hence $W \cap \El^\perp \neq W$. Next, we show that $W \cap \El^\perp \in B_0=X^1 \cup X^2$. Clearly $W\cap \El^\perp \leq \El^\perp$. Thus it suffices to show that $W\cap \El^\perp \pitchfork \EE \cup(\EE + \El)$. Since $W \in X^3$ we know $W \pitchfork \EE \cup (\EE +\El)$. Fix $E\in \EE$. Then $E \nleq \El^\perp$ since otherwise $E^\perp \geq \El$, but $E^\perp \in \EE$ and $\El \pitchfork \EE$. Similarly, $E + \El \nleq \El^\perp$ since otherwise $(E+\El)^\perp \geq \El$, while $(E+ \El)^\perp = E^\perp \cap \El^\perp  \in \EE \cap \El^\perp$ and $\El \pitchfork \EE \cap \El^\perp$. Thus, by Fact \ref{fact: pitchfork plus ell}, it suffices to show that $W \pitchfork (\EE \cap \El^\perp) \cup ((\EE + \El)\cap \El^\perp)$, which holds since $W \in X^3$.  

Next, we show that every vertex (different from $W \cap \El^\perp$) is connected to $W \cap \El^\perp$ hence 
$$\lk_{\kappa_{0}}(W) \cap \beta_{i-1}= \left\{ W\cap \El^\perp \right\}* \flag \left\{ U \in B_{i-1} \setminus \left\{ W \cap \El^\perp \right\} \mid U <W \text{ or } U>W \right\}.$$
Let $U \in \lk_{\kappa_{0}}(W) \cap \beta_{i-1}(0) = \left\{ U \in B_{i-1} \mid U<W \text{ or } U>W \right\}$. If $W \leq U$ then $W\cap \El^\perp \leq U$ and thus $\{U,W\cap \El^\perp \} \in \lk_{\kappa_{0}}(W) \cap \beta_{i-1}$. If $U \leq W$, then $\dim U < \dim W = n+1-i+1$ and hence $U \in B_0$. Therefore either $U\in X^1$ or $U \in X^2$. But if $U\in X^1$, then $\El \leq U \leq W$, a contradiction to $W \in B_i \setminus B_{i-1}$. Thus, $U\in X^2$, and in particular $U\leq \El ^\perp$. Therefore, $U \leq W \cap \El^\perp$ and hence $\{U,W\cap \El^\perp \} \in \lk_{\kappa_{0}}(W) \cap \beta_{i-1}$.

Therefore, by \Cref{cone of a join - single vertex prop}, $\lk_{\kappa_{0}}(W) \cap \beta_{i-1}$ has a cone radius bounded by 1. 
Applying \Cref{adding vertices to cone thm} we get a cone function $\Cone_{\beta_i}$ satisfying
$$\operatorname{Rad}(\operatorname{Cone}_{\beta_{i}}) \leq \operatorname{Rad}(\operatorname{Cone}_{\beta_{i-1}}) +1 \leq n+2 +i.$$

Since $\kappa_0 = \beta_{n+1}$ we get the desired result. 
\end{proof}

From \cite[Proposition 13]{Abr} we get the following.
\begin{fact}
  Let $\kappa_0 \subset \dots \subset \kappa_{2n+2} = \kappa$ be the filtration as above.
  Then for every $w \in \kappa_i(0) \setminus \kappa_{i-1}(0), 1\leq i\leq 2n+2$ we have $\lk_{\kappa_i}(w) \cap \kappa_{i-1}$ is in $\mathbf{C} = \mathbf{C}_A \cup \mathbf{C}_C$ or is a join of two elements from $\mathbf{C}$. 
  \end{fact}

Combining the above results, we get:
\begin{corollary}\label{cor: C satisfies assumptions}
The class $\textbf{C}$ satisfies the requirements of \Cref{thm: generalization of Lemma 22}.
\end{corollary}

\begin{corollary} \label{cor: cone for Cn case ab and nonab}
Let $\Delta=\flag X(V)$ be a classical $C_n$ building as described in \Cref{subsubsec: geom descr C_n} (note that $\dim \Delta = n-1$). Assume that $K$ is infinite or $K=\mathbb{F}_q$ and 
\begin{itemize}
\item $q \geq 2^{2n-2}$ if $f$ is alternating and $Q=0$,
\item $q \geq 2^{4n-4}$ if $m=2n$ and $\sigma \neq \operatorname{id}$,
\item $q \geq 2^{2n-1}$ if $m=2n+1$ or $2n+2$.
\end{itemize}
In particular, if $\Delta$ is the building coming from the BN-pair of a Chevalley group of type $B_n$ or $C_n$ over the field $\mathbb{F}_q$ we require $q \geq 2^{2n-1}$. 

Then for any $a \in \Delta$ there exists an $(n-2)$-cone function $\Cone_{\Delta^0(a)}$ of $\Delta^0(a)$ with 
$$\Rad_j(\Cone_{\Delta^0(a)}) \leq \mathcal{R}(n-1) \qquad \forall -1 \leq j \leq n-2$$
where $\mathcal{R}(n-1)$ is as in \Cref{thm: generalization of Lemma 22} with $f(n)=2n+1$ and $\ell_n = 2n+2$.
\end{corollary}

\subsection{Buildings of type $D_n$}
The proof that opposite complexes of buildings of type $D_n$ are coboundary expanders is the most involved case. We will use the notation introduced in \Cref{subsubsec: geom descr D_n}.
We will later apply an adapted version of the general induction argument to a class of complexes containing flag complexes of the form $T_{\EE}(V) = \flag X_{\EE}(V)$. But by \Cref{prop: description oppo Dn}, the opposite complexes in a building of type $D_n$ are the oriflamme complexes $\tilde{T}_{\EE}(V)$. Thus we start by showing that a cone function for $T_{\EE}(V)$ also gives rise to a cone function for $\tilde{T}_{\EE}(V)$. This can be seen as a quantitative version of \cite[Corollary 17 (ii)]{Abr}.

\begin{proposition} \label{prop: quantitative subdivision D_n}
Let $\tau = \{E_1<\dots < E_r\} \in \tilde{\Delta}$ and set $\EE(\tau) = \{E_i, E_i^\perp \mid 1 \leq i\leq r\}$. If $T_{\EE(\tau)}(V)$ has an $(n-1)$-cone function $\Cone_{T_{\EE(\tau)}(V)}$ with $\Rad(\Cone_{T_{\EE(\tau)}(V)}) \leq c$ for some $c \geq 0$ then there exists an $(n-1)$-cone function $\Cone_{\tilde{T}_{\EE(\tau)}(V)}$ with $\Rad(\Cone_{\tilde{T}_{\EE(\tau)}(V)}) \leq 2c$. 
\end{proposition}

\begin{proof}
For shorter notation, we write $T = T_{\EE(\tau)}(V), Y = X_{\EE(\tau)}(V)$ and $\tilde{T} = \tilde{T}_{\EE(\tau)}(V), \tilde{Y}= \tilde{X}_{\EE(\tau)}(V)$.

The idea of the proof is to construct an $(n-1)$-cone function for $\tilde{T}$ given a cone function 
$$\Cone_T: \bigoplus_{j=-1}^{n-1}C_j(T) \to \bigoplus_{j=0}^n C_j(T)$$ 
by defining maps $f_j: C_j(\tilde{T}) \to C_j(T)$ for $-1\leq j \leq n-1$ and $g_{j}: C_j(T) \to C_j(\tilde{T})$ for $0 \leq j \leq n$ such that $g_{j+1} \circ \Cone_T|_{C_j(T)} \circ f_j$ gives rise to a cone function of $\tilde{T}$. 

As already observed in \cite{Abr}, $T$ is a subdivision of $\tilde{T}$ in the following way. Note that $\tilde{Y} \subset Y$ and $Y \setminus \tilde{Y} = \{U \in Y \mid \dim U = n-1\}$. If two subspaces $U,W \in \tilde{Y}$ are connected by an edge in $T$ then the same is true in $\tilde{T}$, but we have additional edges in $\tilde{T}$, namely in the case when $\dim U = \dim W = n$ and $\dim (U \cap W) = n-1$. To go from $\tilde{T}$ to $T$, we subdivide the edge $\{U,W\}$ by adding the vertex $U \cap W$ and connecting $U,W$ and every vertex $U' \in \tilde{Y}$ for which $\{U',U,W\} \in \tilde{T}$ to $U\cap W$. 

Recall from the discussion in \Cref{subsubsec: geom descr D_n} that there are two different types of $n$-dimensional spaces in $\tilde{X}$ and that spaces from the same type cannot have an intersection of dimension $n-1$. Furthermore, each $U \in X$ of dimension $n-1$ is contained in exactly two totally isotropic spaces of dimension $n$. If $U \in Y$ of dimension $n-1$, i.e. $U \mathbin{\tilde{\pitchfork}} \EE$, with $W_1,W_2 \in X$ of dimension $n$ such that $U = W_1 \cap W_2$ then by \cite[Lemma 31 (iii)]{Abr} we have that $W_i \mathbin{\tilde{\pitchfork}}\EE$ and hence $W_1,W_2 \in Y$.  

In particular, any simplex in $\tilde{T}$ contains at most two vertices which are spaces of dimension $n$.  

For each $U \in Y \setminus \tilde{Y}$ pick one of the two $n$-dimensional subspaces containing it and denote it by $W_U$. Note that if $W \in Y$ is connected to $U$ by an edge in $T$, i.e. $W \subseteq U$ or $U\subseteq W$ then $\{W_U, W\} \in T(1)$ if and only if $W \neq W_U$.

We define the following two subsets of $T(j)$ for $0 \leq j \leq n$: 
$$I_j = \{\sigma \in T(j) \mid \exists U \in \sigma(0): \dim U = n-1\}, \ J_j = \{\sigma \in T(j) \mid \exists U \in \sigma(0): \dim U = n-1 \text{ and } W_U \notin \sigma(0) \}.$$
Note that $I_j = T(j)\setminus \tilde{T}(j)$.
Let $\sigma \in J_j$ then we define $\sigma \langle U, W_U \rangle \in \overrightarrow{T(j)}$ to be the oriented simplex with vertex $W_U$ at the position of $U$ instead of $U$. This is well defined by the above observation. 
Similarly, if $\sigma \in \overrightarrow{\tilde{T}(j)}$ with $W_1,W_2 \in \sigma(0), \dim W_1 = \dim W_2 = n$ and $U= W_1 \cap W_2$ has dimension $n-1$, we write $\sigma \langle W_i, U \rangle$ for the oriented simplex where we replaced $W_i$ with $U$. This is now a well-defined simplex in $\overrightarrow{T(j)}$.

Next, we define the function $f_j: C_j(\tilde{T}) \to C_j(T), 0\leq j \leq n-1$ by its action on chains of the form $\mathbbm{1}_\sigma$ and then extending $\Z$-linearly:
\begin{align*}
f_j(\II_\sigma) = \begin{cases} \mathbbm{1}_\sigma & \text{ if } \sigma \in \overrightarrow{T(j)} \cap \overrightarrow{\tilde{T}(j)}\\
\mathbbm{1}_{\sigma \langle W_1, U\rangle} + \mathbbm{1}_{\sigma \langle W_2, U \rangle} & \text{ if } \sigma \in \overrightarrow{\tilde{T}(j)} \setminus \overrightarrow{T(j)}, W_1,W_2 \in \sigma(0), \dim W_i = n, U = W_1 \cap W_2, \dim U = n-1.
\end{cases}
\end{align*} 
This is indeed well-defined since the only simplices in $\tilde{T}$ that are not simplices in $T$ are those that have two $n$-dimensional subspaces as vertices. 

On the other hand, we define 
\begin{align*}
g_j : C_j(T) & \to C_j(\tilde{T}) \\
\mathbbm{1}_\sigma & \mapsto \begin{cases} \mathbbm{1}_\sigma & \text{ if } \sigma \in \overrightarrow{T(j)} \cap \overrightarrow{\tilde{T}(j)}\\
\mathbbm{1}_{\sigma \langle U, W_U \rangle} & \text{ if } \sigma \in J_j \\
0 & \text{ if } \sigma \in I_j \setminus J_j.
\end{cases}
\end{align*}

We define $\Cone^j_{\tilde{T}} = g_{j+1} \circ \Cone_T^j \circ f_j$ for all $-1\leq j \leq n-1$, where $\Cone_T^j := \Cone_T|_{C_j(T)}$ and $\Cone_{\tilde{T}} = \oplus_{j=-1}^{n-1} \Cone_{\tilde{T}}^j$. The majority of the rest of the proof is dedicated to showing that this is indeed a cone function. 
As concatenation of $\Z$-linear functions it is $\Z$-linear. 
If $[U]$ is the apex of $\Cone_T$ then the apex of $\Cone_{\tilde{T}}$ is again $[U]$ if $\dim U \neq n-1$ and otherwise it is $[W_U]$.

Hence, we are left with showing the cone equation. 
Let $\sigma \in \overrightarrow{T(j)}$. We write $\Cone_T(\mathbbm{1}_\sigma) = \sum_{\tau \in \overrightarrow{T(j+1)}} \lambda_\tau \mathbbm{1}_\tau$ for $\lambda_\tau \in \Z$. In the following, if $\tau \in I_{j+1}$ is a simplex, then $U$ will always denote the subspace of dimension $n-1$ in $\tau (0)$. Additionally, $[\sigma:\tau]$ will denote the oriented incidence number, which is $0$ if $\tau \not \subseteq \sigma$ and 1 or $-1$ if $\tau \subset \sigma$ depending on the position of the vertex of $\sigma$ that was removed to obtain $\tau$. It is chosen in such a way that for $\sigma \in \overrightarrow{T(j+1)}$ we have $\partial_{j+1} \II_\sigma = \sum_{\tau \in \overrightarrow{T(j)}} [\sigma:\tau] \II_\tau$. In the calculations below we will slightly abuse notation as $\Cone_T(\mathbbm{1}_\sigma)$ is not a chain in $\tilde{T}$, but this is fine since the problematic terms sum to zero. We have
\begin{align*}
\partial_{j+1} (g_{j+1}(\Cone_T (\mathbbm{1}_\sigma))) &= \partial_{j+1}(\Cone_T (\II _\sigma)) - \sum_{\tau \in I_{j+1}} \lambda_\tau \partial_{j+1} \II_\tau + \sum_{\tau \in J_{j+1}} \lambda_\tau \partial_{j+1} \II_{\tau \langle U, W_U \rangle} \\
&= \II_\sigma - \Cone_T(\partial_j \II_\sigma)-\sum_{\tau \in I_{j+1}}\lambda_\tau \sum_{\gamma \in \overrightarrow{T(j)}} [\tau: \gamma] \II_\gamma\\
& + \sum_{\tau \in J_{j+1}} \lambda_\tau [\tau:\tau \setminus \{U\}]\II_{\tau \setminus \{U\}} + \sum_{\tau\in J_{j+1}} \lambda_\tau \sum_{\gamma \in \overrightarrow{T(j)}: \gamma \neq \tau \setminus \{U\}} [\tau: \gamma] \II_{\gamma \langle U, W_U \rangle}
\end{align*} 
where we used that 
$$\partial_{j+1} \II_{\tau \langle U, W_U \rangle} = \sum_{\gamma \in \overrightarrow{T(j)}: \gamma \neq \tau \setminus \{U\}} [\tau:\gamma] \II_{\gamma \langle U, W_U \rangle} + [\tau: \tau \setminus \{U\}] \II_{\tau \setminus \{U\}}.$$ 
Note that $(\tau \langle U, W_U \rangle)\setminus \{W_U\} = \tau \setminus \{U\}$ and $[\tau: \gamma \langle U, W_U \rangle] = [\tau: \gamma]$.

We take a closer look at the second half of that expression: 
\begin{align*}
&-\sum_{\tau \in I_{j+1}}\lambda_\tau \sum_{\gamma \in \overrightarrow{T(j)}} [\tau: \gamma] \II_\gamma + \sum_{\tau \in J_{j+1}} \lambda_\tau [\tau:\tau \setminus \{U\}]\II_{\tau \setminus \{U\}} + \sum_{\tau\in J_{j+1}} \lambda_\tau \sum_{\gamma \in \overrightarrow{T(j)}: \gamma \neq \tau \setminus \{U\}} [\tau: \gamma] \II_{\gamma \langle U, W_U \rangle} \\
= &-\sum_{\tau \in I_{j+1}}\lambda_\tau \sum_{\gamma \in \overrightarrow{T(j)}:\gamma \neq \tau \setminus \{U\}} [\tau: \gamma] \II_\gamma + \sum_{\tau\in J_{j+1}} \lambda_\tau \sum_{\gamma \in \overrightarrow{T(j)}: \gamma \neq \tau \setminus \{U\}} [\tau: \gamma] \II_{\gamma \langle U, W_U \rangle} \\
&-\sum_{\tau \in I_{j+1}}\lambda_\tau [\tau: \tau \setminus \{U\}] \II_{\tau \setminus \{U\}} + \sum_{\tau \in J_{j+1}} \lambda_\tau [\tau:\tau \setminus \{U\}]\II_{\tau \setminus \{U\}} \\
= & -\sum_{\tau \in I_{j+1}}\lambda_\tau \sum_{\gamma \in \overrightarrow{T(j)}:\gamma \neq \tau \setminus \{U\}} [\tau: \gamma] \II_\gamma + \sum_{\tau\in J_{j+1}} \lambda_\tau \sum_{\gamma \in \overrightarrow{T(j)}: \gamma \neq \tau \setminus \{U\}} [\tau: \gamma] \II_{\gamma \langle U, W_U \rangle} \\
&-\sum_{\tau \in I_{j+1}\setminus J_{j+1}}\lambda_\tau [\tau: \tau \setminus \{U\}] \II_{\tau \setminus \{U\}}
\end{align*}

On the other hand, let $\lambda'_\gamma \in \Z$ such that $\Cone_T(\partial \II_\sigma) = \sum_{\gamma \in \overrightarrow{T(j)}} \lambda'_\gamma \II_\gamma$. Then 
\begin{align*}
g_j(\Cone_T(\partial \II_\sigma))  = \Cone_T(\partial \II_\sigma) - \sum_{\gamma \in I_j} \lambda'_\gamma \II_\gamma + \sum_{\gamma \in J_j} \lambda'_\gamma \II_{\gamma \langle U, W_U \rangle}.
\end{align*}
What are the $\lambda'_\gamma$ in terms of the $\lambda_\tau$?
\begin{align*}
\Cone_T(\partial \II_\sigma) &= \II_\sigma - \partial_{j+1} \Cone_T(\II_\sigma) = \II_\sigma - \partial_{j+1} \sum_{\tau \in \overrightarrow{T(j+1)}} \lambda_\tau \II_\tau\\
&= \II_\sigma - \sum_{\tau \in \overrightarrow{T(j+1)}} \lambda_\tau \sum_{\gamma \in \overrightarrow{T(j)}} [\tau:\gamma] \II_\gamma \\
&= \II_\sigma - \sum_{\gamma \in \overrightarrow{T(j)}} \left( \sum_{\tau \in \overrightarrow{T(j+1)}} \lambda_\tau [\tau:\gamma] \right) \II_\gamma.
\end{align*}
Hence if $\gamma \neq \sigma$ then 
$$\lambda'_\gamma = - \sum_{\tau \in \overrightarrow{T(j+1)}}\lambda_\tau [\tau:\gamma]$$
and otherwise
$$\lambda'_\sigma = 1-\sum_{\tau \in \overrightarrow{T(j+1)}}\lambda_\tau [\tau:\sigma].$$
To compare the expression of $g_j(\Cone_T(\partial_j \II_\sigma))$ with $\partial_{j+1} (g_{j+1}(\Cone_T (\mathbbm{1}_\sigma)))$ we furthermore need the following observation comparing the sets $I_j$ with $I_{j+1}$ and $J_j$ with $J_{j+1}$.
\begin{align*}
\left\{(\tau,\gamma) : \gamma \in I_j, \tau \in \overrightarrow{T(j+1)}, \gamma \subset \tau \right\} = \left\{(\tau,\gamma) : \tau \in I_{j+1}, \gamma \in \overrightarrow{T(j)}, \gamma \subset \tau , \gamma \neq \tau \setminus \{U\} \right\} 
\end{align*}
and
\begin{align*}
&\left\{(\tau,\gamma) : \gamma \in J_j, \tau \in \overrightarrow{T(j+1)}, \gamma \subset \tau \right\}\\
= &\left\{(\tau,\gamma) : \tau \in J_{j+1}, \gamma \in \overrightarrow{T(j)}, \gamma \subset \tau , \gamma \neq \tau \setminus \{U\} \right\} \cup \left\{(\tau,\gamma) : \tau \in I_{j+1}\setminus J_{j+1}, \gamma = \tau \setminus \{W_U\} \right\}.
\end{align*}
For $\sigma \in \overrightarrow{T(j)}$ we define
\begin{align*}
	h_\sigma = \begin{cases}
		0 & \text{if } \sigma \notin I_j \\ \II_\sigma & \text{if } \sigma \in I_j \setminus J_j \\ \II_\sigma - \II_{\sigma\langle U,W_U \rangle} & \text{if } \sigma \in J_j
	\end{cases}.
\end{align*}
We use the above result and the new notation to further investigate the second half of the expression for $g_j(\Cone_T(\partial \II_\sigma))$.
\begin{align*}
& -\sum_{\gamma \in I_j} \lambda'_\gamma \II_\gamma + \sum_{\gamma \in J_j} \lambda'_\gamma \II_{\gamma \langle U, W_U \rangle} = -\sum_{\gamma \in I_j \setminus \{\sigma\}} \left(- \sum_{\tau \in \overrightarrow{T(j+1)}}\lambda_\tau [\tau:\gamma]\right) \II_\gamma + \sum_{\gamma \in J_j\setminus \{\sigma\}} \left(- \sum_{\tau \in \overrightarrow{T(j+1)}} \lambda_\tau [\tau:\gamma] \right) \II_{\gamma \langle U, W_U \rangle} \\[1em] 
&\qquad - \begin{cases}
	0 & \text{ if } \sigma \notin I_j \\ 1- \sum_{\tau \in \overrightarrow{T(j+1)}}\lambda_\tau [\tau:\sigma]  \II_\sigma & \text{ if } \sigma \in I_j
\end{cases} + \begin{cases}
	0 & \text{ if } \sigma \notin J_j \\  1- \sum_{\tau \in \overrightarrow{T(j+1)}}\lambda_\tau [\tau:\sigma]  \II_{\sigma \langle U, W_U \rangle} & \text{ if } \sigma \in J_j 
\end{cases}\\[1em]
&=- \sum_{\gamma \in I_j }\left( - \sum_{\tau \in \overrightarrow{T(j+1)}}\lambda_\tau [\tau:\gamma] \right)\II_\gamma + \sum_{\gamma \in J_j}\left(- \sum_{\tau \in \overrightarrow{T(j+1)}} \lambda_\tau [\tau:\gamma] \right) \II_{\gamma \langle U, W_U \rangle} - h_\sigma \\[1em]
 &=  \sum_{\tau \in I_{j+1}} \sum_{\gamma \in \overrightarrow{T(j)}: \gamma \neq \tau \setminus \{U\}} [\tau:\gamma]\lambda_\tau \II_\gamma - \sum_{\tau \in J_{j+1}} \sum_{\gamma \in \overrightarrow{T(j)}: \gamma \neq \tau \setminus \{U\}} [\tau:\gamma] \lambda_\tau \II_{\gamma \langle U, W_U \rangle} \\
 & \qquad- \sum_{\tau \in I_{j+1} \setminus J_{j+1}} [\tau:\tau \setminus \{W_U\}] \lambda_\tau \II_{(\tau \setminus \{W_U\})\langle U, W_U \rangle}-h_\sigma.
\end{align*}

We now look at almost the complete cone equation for $\Cone_{\tilde{T}}$, still not considering $f_j$. Let $\sigma \in \overrightarrow{T(j)}$, then we get
\begin{align*}
&\partial_{j+1} g_{j+1}(\Cone_T(\II_\sigma)) + g_{j}(\Cone_T(\partial_j \II_\sigma))  = \II_\sigma \underbrace{- \Cone_T(\partial_j \II_\sigma)}_{(1)}\\
&\underbrace{ -\sum_{\tau \in I_{j+1}}\lambda_\tau \sum_{\gamma \in \overrightarrow{T(j)}:\gamma \neq \tau \setminus \{U\}} [\tau: \gamma] \II_\gamma}_{(2)} \underbrace{+ \sum_{\tau\in J_{j+1}} \lambda_\tau \sum_{\gamma \in \overrightarrow{T(j)}: \gamma \neq \tau \setminus \{U\}} [\tau: \gamma] \II_{\gamma \langle U, W_U \rangle}}_{(3)} \\
&\underbrace{-\sum_{\tau \in I_{j+1}\setminus J_{j+1}}\lambda_\tau [\tau: \tau \setminus \{U\}] \II_{\tau \setminus \{U\}}}_{(4)}\\
& \underbrace{+ \Cone_T(\partial \II_\sigma)}_{(1)} \underbrace{ + \sum_{\tau \in I_{j+1}} \sum_{\gamma \in \overrightarrow{T(j)}: \gamma \neq \tau \setminus \{U\}} [\tau:\gamma]\lambda_\tau \II_\gamma}_{(2)} \underbrace{ - \sum_{\tau \in J_{j+1}} \sum_{\gamma \in \overrightarrow{T(j)}: \gamma \neq \tau \setminus \{U\}} [\tau:\gamma] \lambda_\tau \II_{\gamma \langle U, W_U \rangle} }_{(3)}\\
 &\underbrace{- \sum_{\tau \in I_{j+1} \setminus J_{j+1}} [\tau:\tau \setminus \{W_U\}] \lambda_\tau \II_{(\tau \setminus \{W_U\})\langle U, W_U \rangle}}_{(4)} -h_\sigma.
\end{align*}
The parts marked (1),(2), and (3) cancel each other out. It remains to show that the same is true for the parts marked (4). 
Note that the sums in the parts marked (4) run over the same elements. Hence we need to show that for any $\tau \in I_{j+1}\setminus J_{j+1}$ we have
$$[\tau: \tau \setminus \{U\}] \II_{\tau \setminus \{U\}}=-[\tau:\tau \setminus \{W_U\}]  \II_{(\tau \setminus \{W_U\})\langle U, W_U \rangle}$$
Note that as unoriented simplices $\tau \setminus \{U\}$ is the same as $(\tau \setminus \{W_U\})\langle U, W_U \rangle$ (in both cases we have all the vertices of $\tau$ except $U$).
Let $\tau = [U_0,\dots,U_{j+1}]$ and assume that $U=U_\ell$ and $W_U = U_k$ for some $0 \leq \ell, k \leq j+1$. Then $[\tau:\tau\setminus \{U\}] = (-1)^\ell$ and $[\tau: \tau \setminus \{W_U\}] = (-1)^k$. Further note that we can obtain $(\tau\setminus \{W_U\})\langle U, W_U \rangle$ from $\tau \setminus \{U\}$ by cyclically permuting the vertices of index between $\min\{\ell, k\}$ and $\max\{\ell,k\}-1$. The cycle has length $\lvert k-\ell \rvert -1$ and thus signature $(-1)^{k-\ell-1}$. Hence $(-1)^k \II_{(\tau \setminus \{W_U\})\langle U, W_U \rangle} = (-1)^k \cdot(-1)^{k - \ell -1} \II_{\tau \setminus \{U\}} = -(-1)^\ell \II_{\tau \setminus \{U\}}$ which is what we wanted to show.

To summarize, we get
$$\partial_{j+1} g_{j+1}(\Cone_T(\II_\sigma)) + g_{j}(\Cone_T(\partial_j \II_\sigma)) = \II_\sigma-h_\sigma = \begin{cases}
		\II_\sigma & \text{if } \sigma \notin I_j \\ 0 & \text{if } \sigma \in I_j \setminus J_j \\ \II_{\sigma\langle U,W_U \rangle} & \text{if } \sigma \in J_j
\end{cases}.$$

Note that if $\sigma \in \overrightarrow{T(j)}\cap \overrightarrow{\tilde{T}(j)}$, then $f_j(\II_\sigma) = \II_\sigma$ which shows that the cone equation holds in this case. 
To conclude the proof, we need to check the cone equation also for the case where $\sigma \in \overrightarrow{\tilde{T}(j)} \setminus \overrightarrow{T(j)}$. In this case $\sigma(0)$ contains two subspaces $W_1,W_2$ with $\dim W_i = n$ and $U:=W_1\cap W_2$ has dimension $n-1$. Without loss of generality, assume that $W_U = W_1$. We get $f_j(\II_\sigma) = \mathbbm{1}_{\sigma \langle W_1, U \rangle} + \mathbbm{1}_{\sigma \langle  W_2 , U \rangle}$. We have that $\sigma \langle W_1, U \rangle \in J_j$ and $\sigma \langle W_2, U \rangle \in I_j \setminus J_j$ thus we can apply the above reasoning to both of them separately and get, using linearity: 
\begin{align*}
&\partial_{j+1} g_{j+1}(\Cone_T(f_j(\II_\sigma))) + g_{j}(\Cone_T(f_{j-1}(\partial_j \II_\sigma))) \\
&= \partial_{j+1} g_{j+1}(\Cone_T(\II_{\sigma \langle W_1, U \rangle})) + g_{j}(\Cone_T(\partial_j \II_{\sigma \langle W_1, U \rangle}))\\
& \qquad + \partial_{j+1} g_{j+1}(\Cone_T(\II_{\sigma \langle W_2, U \rangle})) + g_{j}(\Cone_T(\partial_j \II_{\sigma \langle W_2, U \rangle})) \\
&= \II_{(\sigma \langle W_1, U \rangle)\langle U, W_U \rangle} + 0 = \II_\sigma.
\end{align*} 
For the last equality, we used that $(\sigma \langle W_1, U \rangle)\langle U, W_U \rangle = \sigma$ since we first replace $W_1$ in $\sigma$ by $U$ and by applying $\langle U, W_U \rangle$ we replace $U$ by $W_U = W_1$ and end up where we started. We furthermore use that $f_{j-1}\partial_j \II_\sigma = \partial_j f_j \II_\sigma$. This is indeed true as the following argument shows. We have
$$
\partial_j \II_\sigma=\sum_{\left\{W_1, W_2\right\} \subset \gamma \subset \sigma}[\sigma: \gamma] \II_\gamma+\sum_{i \in\{1,2\}}\left[\sigma: \sigma \backslash\left\{W_i\right\}\right] \II_{\sigma \backslash\left\{W_i\right\}} 
$$
and thus
$$
f_{j-1} \partial_j \II_\sigma=\sum_{i \in\{1,2\}}\left(\left[\sigma: \sigma \backslash\left\{W_i\right\}\right] \II_{\sigma \backslash\left\{W_i\right\}}+\sum_{\left\{W_1, W_2\right\} \subset \gamma \subset \sigma}[\sigma: \gamma] \II_{\gamma\langle W_i, U\rangle}\right) .
$$
On the other hand, 
$$
\partial_j f_j \II_\sigma=\sum_{i \in\{1,2\}} \partial_j \II_{\sigma\langle W_i, U\rangle}=\sum_{i \in\{1,2\}} \sum_{\gamma \subset \sigma\langle W_i, U\rangle}\left[\sigma\langle W_i, U\rangle: \gamma\right] \II_\gamma.
$$
The set $\left\{\gamma:\left\{W_1, W_2\right\} \subset \gamma \subset \sigma\right\}$ is in bijection with the set $\left\{\gamma:\left\{U, W_{3-i}\right\}  \subset \gamma \subset \sigma\langle W_i, U\rangle\right\}$ via the map $\gamma \mapsto \gamma\langle W_i, U\rangle$, which preserves the oriented incidence number. Thus,
$$
\begin{aligned}
& \sum_{\left\{W_1, W_2\right\} \subset \gamma \subset \sigma}[\sigma: \gamma] \II_{\gamma\langle W_i, U\rangle}=\sum_{\left\{U, W_{3-i}\right\} \subset \gamma \subset \sigma\langle W_i, U\rangle}\left[\sigma\langle W_i, U\rangle: \gamma\right] \II_\gamma \\
& \quad=\partial_j \II_{\sigma\langle W_i, U\rangle}-\left[\sigma\langle W_i, U\rangle: \sigma\langle W_i, U\rangle \backslash\left\{W_{3-i}\right\}\right] \II_{\sigma\langle W_i, U\rangle \backslash\left\{W_{3-i}\right\}}-\left[\sigma\langle W_i, U\rangle: \sigma \backslash\{U\}\right] \II_{\sigma\langle W_i, U\rangle \backslash\{U\}}.
\end{aligned}
$$

The last term is equal to $\left[\sigma: \sigma \backslash\left\{W_i\right\}\right] \II_{\sigma \backslash\left\{W_i\right\}}$ (since deleting $W_i$ or swapping $W_i$ to $U$ and then deleting $U$ does not change orientation). Adding everything up, it remains to check that
$$
\left[\sigma\langle W_1, U\rangle: \sigma\langle W_1, U\rangle \backslash\left\{W_2\right\}\right] \II_{\sigma\langle W_1, U\rangle \backslash\left\{W_2\right\}}+\left[\sigma\langle W_2, U\rangle: \sigma\langle W_2, U\rangle \backslash\left\{W_1\right\}\right] \II_{\sigma\langle W_2, U\rangle \backslash\left\{W_1\right\}}=0.
$$

To see that this holds, we assume without loss of generality (by applying cyclic shifts if necessary) that $\sigma=\left[U_0, \ldots, U_{j}\right]$ with $U_0=W_1$ and $U_{k}=W_2$. Then we have 
\begin{align*}
	\left[\sigma\langle W_1, U\rangle: \sigma\langle W_1, U\rangle \backslash\left\{W_2\right\}\right]&=(-1)^{k}, \\
\left[\sigma\langle W_2, U\rangle: \sigma\langle W_2, U\rangle \backslash\left\{W_1\right\}\right]&=1. \\
\end{align*}

Since $\sigma\langle W_2, U\rangle \backslash\left\{W_1\right\}=[U_2,\dots,U_{k-1},U,U_{k+1},\dots,U_j]$ may be obtained from $\sigma\langle W_1, U\rangle \backslash\left\{W_2\right\} = [U,U_2,\dots,U_{k-1},U_{k+1},\dots,U_j]$ via cyclically shifting the coordinates between 0 and $k-1$, which is a cycle of length $k$ and has signature $(-1)^{k-1}$, we get
 $$\II_{\sigma\langle W_2, U\rangle \backslash\left\{W_1\right\}}= (-1)^{k-1} \II_{\sigma\langle W_2, U\rangle \backslash\left\{W_2\right\}}$$
which concludes the argument.

For bounding the cone radius of $\Cone_{\tilde{T}}$ note that $f_j$ at most doubles the size of support of a chain, while applying $g_j$ will not increase the size of the support of a chain. Hence in the worst case, i.e. if $\sigma \in I_j$ then 
$$\lvert \supp (\Cone_{\tilde{T}}(\II_\sigma)) \rvert \leq \lvert \supp (\Cone_{T}(\mathbbm{1}_{\sigma \langle W_1, U \rangle} )) \rvert + \lvert \supp (\Cone_{T}(\mathbbm{1}_{\sigma \langle W_2,  U \rangle})) \rvert \leq 2 \Rad_j \Cone_T \leq 2c.$$
Hence we get the desired bound on the cone radius. 
\end{proof}

To show that the flag complexes $T_{\EE}(V)$ as above admit cone functions, we first need a weaker and more specialized version of \Cref{thm: generalization of Lemma 22} which allows for arguments via a filtration also for the links.

\begin{theorem}  \label{thm: the big induction for the D_n case}
	Let $\mathbf{C}$ be a class of pairs $(\kappa, Y)$ where $Y\neq \emptyset$ is a poset and $\kappa = \flag Y$ (with inclusion as incidence relation).  
	Let $n = \dim \kappa $, assume there exists $\ell_n \in \N$ and a filtration $\kappa_0 \subset \kappa_1 \subset \dots \subset \kappa_\ell = \kappa$ with $\ell \leq \ell_n$ satisfying:
\begin{enumerate}
\item There exists an $(n-1)$-cone function $\Cone_{\kappa_0}$ such that 
$$\Rad_j(\Cone_{\kappa_0}) \leq f(n) , \forall -1 \leq j \leq n-1$$
where $f: \R \to \R$ is a function depending only on $\mathbf{C}$.
\item Every $\kappa_i$ is a full subcomplex of $\kappa$, i.e. $\kappa_i = \flag Y_i$ for some $Y_i \subseteq Y$ and all elements in $Y_i \setminus Y_{i-1}$ are not comparable (i.e. not incident). 
\item For every $U \in Y_i \setminus Y_{i-1}$ we have that $\lk_{\kappa_i}(U) \cap \kappa_{i-1} = \kappa_{i-1}^{<U} * \kappa_{i-1}^{>U}$ is such that $\kappa_{i-1}^{<U}, \kappa_{i-1}^{>U} \in \mathbf{C}$ or one (or both) have again a filtration $R_0 \subset \dots R_h, h \leq \ell_{\dim R_0}$ which starts with an element of $\mathbf{C}$ (i.e. $R_0\in \mathbf{C}$), such that condition 2. is satisfied and such that for every $W \in R_j(0) \setminus R_{j-1}(0)$ we have $R_{j-1}^{<W},R_{j-1}^{>W} \in \mathbf{C}$. 
\end{enumerate} 
Then there exist constants $\mathcal{R}(n)$ depending only on $\mathbf{C}$ and $n$, and an $(n-1)$-cone function $\Cone_\kappa$ such that 
$$\Rad_j(\Cone_\kappa) \leq \mathcal{R}(n), \forall -1 \leq j \leq n-1.$$ 
\end{theorem}
\begin{proof}
The proof is analogous to the one of \Cref{thm: generalization of Lemma 22}.
The main difference is that to get a cone function for $\lk_{\kappa_i}(U) \cap \kappa_{i-1}$ it might happen that we cannot just use the induction hypothesis, but that we need to do another induction on the filtration of $\kappa_{i-1}^{<U}$ or $\kappa_{i-1}^{>U}$, using \Cref{adding vertices to cone thm} in the step of the induction. This gives an even weaker and more complicated bound on the cone radius of the final cone function, thus we decided to not keep track of it. But it will still only depend on $n$ and the class $\mathbf{C}$, not on the field $K$. 
\end{proof}

We now describe the set up for the rest of this section. 
  \begin{enumerate}
    \item We define a class $\mathbf{C}_1$.
    \item We define a filtration $\kappa_0 \subset \dots \subset \kappa_\ell = \kappa$ for each $\kappa \in \mathbf{C}_1$.
    \item We show that $\kappa_0$ admits a cone function with bounded radius.
    \item We conclude that $\mathbf{C}_1$ satisfies the assumptions of the adapted general induction argument. 
    \item We define a class $\mathbf{C}_2$.
    \item We proceed as for $\mathbf{C}_1$ (for this class the original induction argument suffices).
    \item We define a class $\mathbf{C_3}$.
    \item We define a filtration $\kappa_0 \subset \dots \subset \kappa_\ell = \kappa$ for each $\kappa \in \mathbf{C}_3$.
    \item We show that $\kappa_0$ admits a cone function with bounded radius.
    \item We show that $\mathbf{C_1}\cup \mathbf{C_2} \cup \mathbf{C}_3$ satisfies the assumptions of the adapted general induction argument. 
    \item We conclude that all opposite complexes of type $D_n$ admit a cone function.
  \end{enumerate}

We fix $n$ to be the Witt index of $(V,Q,f)$ hence $\dim V = 2n$ and we are dealing with a $D_n$ building, which is $n-1$-dimensional.
Recall $X= \{0<U<V \mid U \text{ totally isotropic}\}$.
\begin{definition}[{\cite[Definition 12]{Abr}}]
Let $U,E <V$, $\dim  U <n, \dim E = n$. Assume that $E$ is not totally isotropic. Then 
$$U\mathbin{@}E :\iff U^\perp \cap E \text{ is not totally isotropic.}$$
Set 
\begin{align*}
\mathcal{U}_n &= \mathcal{U}_n(V) = \{A<V \mid \dim A = n, A\notin X\} \\
\mathcal{M}(A) &= \{M<A \mid \dim M = n-1, M \in X\} \text{ for any } A \in \mathcal{U}_n.
\end{align*}
\end{definition}

\begin{definition} \label{def: class 1 D_n case}
We define the class $\mathbf{C}_1$ to consist of complexes $\flag X_{\EE; \mathcal{F}}(U;V)$ where 
\begin{itemize}
\item $(V,Q,f)$ is a thick pseudo-quadratic space of Witt index $n$, in particular $\dim V = 2n$;
\item $U\in X$ with $\dim U = k \geq 2$;
\item $\EE$ a finite set of subspaces of $U$, set $e_j = \lvert \EE_j \rvert, 1 \leq j \leq k-1$ like before;
\item Let $\mathcal{F} = \mathcal{F}_1 \supset \mathcal{F}_2 \supset \dots \supset \mathcal{F}_{k-1}$ be finite subsets of $\mathcal{U}_n(V)$ satisfying
\begin{enumerate}
\item $U \cap F \in \EE \cup \{0\}$ for all $F\in \mathcal{F}$;
\item $\mathcal{F}_j^\perp = \mathcal{F}_j$;
\item $\dim(U\cap F) \leq k-1-j$ for all $F \in \mathcal{F}_j, 1 \leq j \leq k-1$.
\end{enumerate}
\item Assume that $\lvert K \rvert \geq \sum_{j=1}^{k-1} \binom{k-2}{j-1} e_j + 2^{k-1}s$, where $s = \lvert \mathcal{F} \rvert$.
\end{itemize} 
Set 
$$X_{\EE; \mathcal{F}}(U;V):= \{0< W < U \mid W \pitchfork_U \EE \text{ and } W \mathbin{@}_V \mathcal{F}_j \text{ for } \dim W = j \}.$$
\end{definition}
\begin{fact}[{\cite[Lemma 33, Step 1]{Abr}}]\label{D1 the Y0 fact}
	Let $\kappa = \flag X_{\EE; \mathcal{F}}(U;V) \in \mathbf{C}_1$ with $ k = \dim U \geq 2$, in particular $\kappa$ has dimension $k-2$. Then there exists a 1-dimensional subspace $\ell \in X_{\EE; \mathcal{F}}(U;V)$.
\end{fact}

For $\kappa = \flag X_{\EE; \mathcal{F}}(U;V) \in \mathbf{C}_1$ with $ k = \dim U \geq 3$, we set $Y= X_{\EE; \mathcal{F}}(U;V)$ and define a filtration as follows: We fix a 1-dimensional subspace $\ell \in X_{\EE; \mathcal{F}}(U;V)$ and set
$$Y_0 := \{A \in Y \mid A + \ell \in Y\}.$$
The filtration will then consist of two dimension decreasing filtrations. First we define
\begin{align*}
Z &:= \{W \in Y \mid \ell \leq W \text{ or } W \pitchfork_U(\mathcal{F} \cap U) +\ell \}\\
Y_i &:= \{W \in Z \mid W \in Y_0 \text{ or } \dim W \geq k-i \},\\
\kappa_i &:= \flag Y_i, \ 1 \leq i \leq k-1.
\end{align*}
The second part of the filtration of $\kappa$ is defined as follows:
\begin{align*}
Z_i&:= \{W \in Y \mid W \in Z \text{ or } \dim W \geq k-1\}, \\
\kappa_{k+i} &:= \flag Z_i, 0 \leq i \leq k-1.
\end{align*}

\begin{fact}[{\cite[Step 2, Lemma 33]{Abr}}] \label{fact: D1 filtration well defined}
	We have
$$Y_0 = \left\{W \in Y \mid \El \leq W \lor (\dim W < k-1 \land W \pitchfork_U \EE +\El \land W @_V(\cF_{\dim W} \cup ((\cF_{\dim W +1} \cap \El^\perp)+\El) )) \right\} \subseteq Z$$
using that $W \pitchfork_U \EE +\El$ implies $W \pitchfork_U (\cF\cap U)+\El$ since by assumption $\cF \cap U \subseteq \EE$. 
\end{fact}

\begin{lemma}\label{lemma: Y_0 for class 1 D_n case}
With the definitions as above, we have that $\kappa_0$ has a $(k-3)$-cone function $\Cone_{\kappa_0}$ with 
$$\Rad_j(\Cone_{\kappa_0}) \leq k, \ -1 \leq j \leq k-3.$$
\end{lemma}

\begin{proof}
We proceed similar to the proof of \Cref{X0 bound lemma}. Set $B_0 :=\{A \in Y \mid \ell \leq A\}, \beta_0 = \flag B_0$. Then $\beta_0 =\{\ell\} * \lk_Y(\{\ell\})$ and hence it has a $(k-3)$-cone function with cone radius bounded by 1.
Next, set $B_i= \{A \in Y_0 \mid \ell \leq A \text{ or } \dim A \leq i \}, \beta_i = \flag B_i$ for $1 \leq i \leq k-1$. Then all subspaces in $B_i \setminus B_{i-1}$ have dimension $i$. Let $A \in B_i \setminus B_{i-1}$. 
Then $\lk_{\beta_i}(A)\cap \beta_{i-1}= \{A+\ell\} * \flag \{W \in B_{i-1}\mid W \neq A + \ell, A <W \text{ or } W <A \}$ by the same argument as in \Cref{X0 bound lemma}. 
Applying \Cref{adding vertices to cone thm} inductively, we get a $(k-3)$-cone function $\Cone_{\beta_{k-1}} = \Cone_{ \kappa_0}$ with 
$$\Rad_j(\Cone_{\kappa_0}) \leq k, \ -1 \leq j \leq k-3.$$
\end{proof}

We collect some facts about the filtration of $\kappa$ from \cite[Lemma 33]{Abr}.
\begin{fact}[{\cite[Step 3-5, Lemma 33]{Abr}}] \label{fact: D1 first part of filtration}
	Let $W \in Y_i \setminus Y_{i-1}$. Then $Y_{i-1}^{<W} = Y_{0}^{<W}$ (this can be seen directly by looking at the dimension). Furthermore, $Y_0^{<W} = X_{\EE'; \mathcal{F}'}(W;V)$ for some $\EE', \mathcal{F}'$. Thus $\kappa_{i-1}^{<W} \in \mathbf{C}_1$. On the other hand, $Y_{i-1}^{>W} = Z^{>W}$. While $\flag Z^{>W}$ is not necessarily in $\mathbf{C}_1$ it admits itself a filtration of length $i$ satisfying the necessary properties. 
\end{fact}

\begin{fact}[{\cite[Step 4-5, Lemma 33]{Abr}}] \label{fact: D1 second part of filtration}
	Let $W \in Z_i \setminus Z_{i-1}$. Then $Z_{i-1}^{>W} = Y^{>W}$ and $\flag Y^{>W} \in \mathbf{C}_1$. On the other hand, $Z_{i-1}^{<W}=Z^{<W}$ for which we  have $\flag Z^{<W} \in \mathbf{C}_1$.  
\end{fact}

\begin{lemma}\label{lemma: class 1 D_n case}
The class $\mathbf{C}_1$ satisfies the assumptions of \Cref{thm: the big induction for the D_n case}. In particular, any complex in $\mathbf{C}_1$ of dimension $k-2$ has a $(k-3)$-cone function with cone radius bounded by a constant depending only on $k$.  
\end{lemma}
\begin{proof}
Let $\kappa = \flag X_{\EE; \mathcal{F}}(U;V) \in \mathbf{C}_1$ and $k = \dim U = \dim \kappa +2$. By Fact \ref{D1 the Y0 fact} we can fix a 1-dimensional subspace $\ell \in X_{\EE; \mathcal{F}}(U;V) =:Y$. In particular, if $\kappa$ is 0-dimensional (i.e. $k=2$), it is not empty. Thus we can define a $(-1)$-cone function for $\kappa$ with cone radius 1. Hence, we get a filtration of $\kappa$ with $\ell_0 = 0$.
For $k \geq 3$, we use the filtration defined above. By \Cref{lemma: Y_0 for class 1 D_n case}, there exists a $(k-3)$-cone function $\Cone_{\kappa_0}$ such that 
$$\Rad_j(\Cone_{\kappa_0}) \leq k =:f(k),  \forall -1\leq j \leq k-3.$$

Combining Facts \ref{fact: D1 first part of filtration}, \ref{fact: D1 second part of filtration} we see that the filtration defined above indeed satisfies the requirements of \Cref{thm: the big induction for the D_n case}.
\end{proof}

\begin{definition}\label{def: class 2 D_n case}
We define the class $\mathbf{C}_2$ to consist of complexes $\flag Z_{\EE; \mathcal{F}}(U;V)$ where
\begin{itemize}
\item $(V,Q,f)$ is a thick pseudo-quadratic space of Witt index $n$, in particular $\dim V = 2n$;
\item $U\in X$ with $\dim U = k \geq 2$;
\item $\EE$ a finite set of subspaces of $U$, set $e_j = \lvert \EE_j \rvert, 1 \leq j \leq k-1$ like before;
\item $\mathcal{F} \subset \mathcal{U}_n(V)$ with $\lvert \mathcal{F}\rvert = s < \infty$ such that 
\begin{enumerate}
\item $U\cap \mathcal{F} \subseteq \EE \cup \{0\}$,
\item $U \cap F^\perp = 0$ and $\dim(U \cap F) \leq 1$ for all $F \in \mathcal{F}$.
\end{enumerate}
\item Assume $\lvert K \rvert \geq \sum_{j=1}^{k-1} \binom{k-2}{j-1} e_j + 2s$.
\item Set 
$$Z_{\EE; \mathcal{F}}(U;V) := \{0<W<U \mid W \pitchfork_U \EE \text{ and } W \mathbin{@}_V \mathcal{F}\}.$$
\end{itemize}
\end{definition}
Let $\kappa=\flag Z_{\EE; \mathcal{F}}(U;V) \in \mathbf{C}_2$, $\dim U = k \geq 2$, in particular $\dim \kappa = k-2$. Set $Y = Z_{\EE; \mathcal{F}}(U;V)$. 
\begin{fact}[{\cite[Lemma 34, Step 3]{Abr}}]\label{fact: D2 exists hyperplane}
	There exists $H <U$ with $\dim H = k-1$ and $H \in Y$.
\end{fact}
Using such an $H$, we define a filtration for $\kappa$ as follows.
We set
$$Y_0:= \{B\in Y \mid B \cap H \in Y\}$$  
and
\begin{align*}
Y_i &:= \{W \in Y \mid W \in Y_0 \text{ or } \dim W \leq i\},\\
\kappa_i&:= \flag Y_i, \ 1 \leq i \leq k-1. 
\end{align*}

\begin{lemma} \label{lemma: Y_0 for class 2 D_n case}
The complex $\kappa_0 = \flag Y_0$ has a $(k-3)$-cone function $\Cone_{\kappa_0}$ with 
$$\Rad_j(\Cone_{\kappa_0}) \leq k , \ -1 \leq j \leq k-3.$$
\end{lemma}

\begin{proof}
The proof is similar to the proofs of \Cref{X0 bound lemma} and \Cref{lemma: Y_0 for class 1 D_n case}, the difference being that we now look at the dual version, where we fix a hyperplane instead of a line.

Set $B_0 = \{A \in Y \mid A \leq H \}, \beta_0 = \flag B_0$. Then $\beta_0 = \{H \} * \lk_{Y_0}(\{H\})$ hence it has a $(k-3)$-cone function with cone radius bounded by $1$.

Next, set $B_i = \{A \in Y_0 \mid A \leq H \text{ or } \dim A \geq k-i\}, \beta_i = \flag B_i$ for $1 \leq i \leq k-1$. Then the subspaces in $B_i\setminus B_{i-1}$ are all of dimension $k-i$. Let $A \in B_i \setminus B_{i-1}$, then $\lk_{\beta_i}(A) \cap \beta_{i-1} = \{A \cap H \} * \flag \{W \in B_{i-1} \mid W \neq A \cap H, W < A \text{ or } A < W \} $. Hence it has a $(k-4)$-cone function with cone radius bounded by 1. Thus we can apply \Cref{adding vertices to cone thm} inductively and obtain a $(k-3)$-cone function $\Cone_{B_{k-1}} = \Cone_{\kappa_0}$ with 
$$\Rad_j(\Cone_{\kappa_0}) \leq k , \ -1 \leq j \leq k-3.$$

\end{proof}

We recall the following fact about the filtration of $\kappa$.
\begin{fact}[{\cite[Step 4, Step 5, Lemma 34]{Abr}}] \label{fact: D2 filtration good}
	Let $W \in Y_i \setminus Y_{i-1}$. Then $Y_{i-1}^{<W}= Y^{<W} = \{0<A<W \mid A \pitchfork_W \EE \cap W\} = Z_{\EE \cap W; \emptyset}(W;V)$ and hence $\flag Y^{<W} \in \mathbf{C}_2$.
On the other hand, $Y_{i-1}^{>W} = Y_0^{>W}$ can be described in such a way that $\flag Y_0^{>W} \in \mathbf{C}_2$.
\end{fact}

\begin{lemma}\label{lemma: class 2 D_n case}
The class $\mathbf{C}_2$ satisfies the assumptions of \Cref{thm: generalization of Lemma 22}. In particular, any complex in $\mathbf{C}_2$ of dimension $k-2$ has a $(k-3)$-cone function with cone radius bounded depending only on $k$.  
\end{lemma}
\begin{proof}
Let $\kappa= \flag Z_{\EE; \mathcal{F}}(U;V) \in \mathbf{C}_2$. Then using Fact \ref{fact: D2 exists hyperplane} we fix $H<U$ such that $\dim H = \dim U- 1= k-1$ and $H \in Z_{\EE; \mathcal{F}}(U;V)=:Y$. In particular, $\kappa \neq \emptyset$.
We prove the existence of a desired filtration by induction on $k$. 
For $k =2$ we have that $\dim \kappa = 0$ hence the desired filtration and cone function exist trivially. 
Fix $k\geq 3$. We consider the filtration for $\kappa$ as defined above.
This clearly satisfies assumption 2. of \Cref{thm: generalization of Lemma 22}.
Assumption 3. is satisfied by Fact \ref{fact: D2 filtration good}.
Since $\kappa_{k-1}=\kappa$, the filtration is indeed as desired. 
\end{proof}

Next, we will define a third class, which will be the one containing the relevant opposite complexes, see also \Cref{rmk: connection to op compl D_n}.

\begin{definition}\label{def: class 3 D_n case}
We define the class $\mathbf{C}_3$ to consist of complexes of the form $\flag Y_{\EE}(V)$ where
\begin{itemize}
\item $(V,Q,f)$ is a thick pseudo-quadratic space of Witt index $n$, in particular $\dim V = 2n$;
\item $\EE = \EE^\perp$ is a finite set of subspaces of $V$, set $\EE_j= \{E \in \EE : \dim E=j\}, e_j = \lvert \EE_j \rvert$;
\item $\lvert K \rvert \geq 2 \sum_{j=1}^{2n-1}\binom{2n-2}{j-1}e_j$;
\item $\widehat{\EE}_n:= \EE_n \cap \mathcal{U}_n(V)$;
\item $X = \{0<U<V \mid U \text{ is totally isotropic} \}$.
\end{itemize}
Set $$Y_{\EE}(V) :=\{U \in X \mid U \mathbin{\tilde{\pitchfork}}\EE \text{ and if } \dim U < n \text{ then } U\mathbin{@} \hat{\EE}_n \}.$$
\end{definition}

Let $\kappa= \flag Y_{\EE}(V) \in \mathbf{C}_3$ with $\dim V = 2n$.
\begin{fact}[{\cite[Proposition 14, Step 1]{Abr}}] \label{fact: D3 line}
	There exists a 1-dimensional subspace $\ell \in Y=Y_{\EE}(V)$. 
\end{fact}

Assume $n\geq 3$ and fix an $\ell$ as above. Set
$$
\begin{aligned}
& \mathcal{F}:=\mathcal{E} \cup\left(\mathcal{E} \cap \ell^{\perp}\right) \cup(\mathcal{E}+\ell) \cup\left(\left(\mathcal{E} \cap \ell^{\perp}\right)+\ell\right) \text { and } \\
& \widehat{\mathcal{F}}_n:=\mathcal{F}_n \cap \mathcal{U}_n(V). 
\end{aligned}
$$
We further define
\begin{align*}
	X^1 &= \{U \in X \mid \ell \leq U, U \mathbin{\tilde{\pitchfork}} \mathcal{E} \text{ and if } \dim U <n \text{ then } U \mathbin{@} \widehat{\mathcal{E}}_n\}\\
	X^2&= \left\{U\in X \mid \ell \not \leq U \leq \ell^{\perp} \text{ and } \begin{cases} U \pitchfork \mathcal{E}, U \mathbin{@} \widehat{\mathcal{E}}_n & \text{ if } \dim U = n-1 \\ U \pitchfork \mathcal{E} \cup(\mathcal{E}+\ell), U \mathbin{@} \widehat{\mathcal{F}}_n & \text{ if } \dim U \leq n-2 \end{cases} \right\}\\
	X^3 &= \left\{U \in X \mid U \not \leq \ell^{\perp}, \operatorname{dim} U>1, U \widetilde{\pitchfork} \mathcal{F} \text{ and if } \dim U <n \text{ then } U \mathbin{@} \widehat{\mathcal{F}}_n \right\}\\
	X^4 &= \left\{ U \in X \mid U \not \leq \ell ^\perp  \text{ and } \begin{cases}U \mathbin{\tilde{\pitchfork}} \mathcal{F} &\text{ if } \dim U = n \\
U \pitchfork \mathcal{F} \text{ and } U\mathbin{@}\widehat{\mathcal{F}}_n &\text{ if } \dim U = n-1 \\
U \pitchfork \EE \cup (\EE + \ell) \text{ and } U \mathbin{@} \widehat{\mathcal{F}}_n &\text{ if } \dim U = n-2.
 \end{cases} \right\}
\end{align*}

We can now define the filtration for $\kappa$. It will contain two parts, one dimension increasing and one dimension decreasing filtration. We start with
\begin{align*}
	Y_0 = \{ U \in X \mid U \in X^1 \cup X^2 \cup X^3\},
\end{align*}
and set 
\begin{align*}
Z&:=  X^1 \cup X^2 \cup X^4 \\
Y_i&:= \{U \in Z \mid U \in Y_0 \text{ or } \dim U \leq i \}, 1 \leq i \leq n-2\\
Y_{n-1} &:= Z \cup \{U \in Y \mid \dim U = n-1 \text{ and } Y_0^{>U} \neq \emptyset\}.
\end{align*}
Note that $Y_0 \subseteq Z$ since if $\dim U < n$ then $U \mathbin{\tilde{\pitchfork}} \EE \iff U \pitchfork \EE$.
The dimension increasing filtration is defined by setting
$$Y_i:= \{U \in Y \mid U \in Y_{n-1} \text{ or } \dim U \geq 2n-i\}, n \leq i \leq 2n-1.$$
As usual, we define $\kappa_i = \flag Y_i$ and have $\kappa_{2n-1} = \kappa$.

Before constructing a cone function for $\kappa_0$ we recall the following fact.
\begin{fact}[{\cite[Step 3, Proposition 14]{Abr}}] \label{fact: D3 technical}
	For every $U \in Y_0$, we have that $U\cap \ell^\perp \in Y_0$ and $(U\cap \ell^\perp)+\ell \in Y_0$.
\end{fact}

\begin{lemma}\label{lemma: Y_0 for class 3 D_n case}
The complex $\kappa_0 = \flag Y_0$ has an $(n-2)$-cone function $\Cone_{\kappa_0}$ with 
$$\Rad_j(\Cone_{\kappa_0}) \leq 2n+1 , \ -1 \leq j \leq n-2.$$
\end{lemma}
\begin{proof}
Let $\kappa= \flag Y_{\EE}(V) \in \mathbf{C}_3$ with $\dim V = 2n$ and let $Y_0$ be defined as above. Recall that $X = \{0<U<V \mid U \text{is totally isotropic subspace}\}$. We define a filtration for $Y_0$ and we want to apply \Cref{adding vertices to cone thm} inductively, using Fact \ref{fact: D3 technical}.

Set
\begin{align*}
A_0 &:=X^1 \\
A_i &:= X^1 \cup \{U \in X^2 \mid \dim U \leq i  \}, \ 1 \leq i \leq n; \\
\alpha_i &:= \flag A_i, 0 \leq i \leq n.
\end{align*}
Then $A_0 = \{\ell\}* \flag \{U \in A_0 \mid U \neq \ell\}$ and hence has a cone function with cone radius bounded by 1. 
Clearly elements in $A_i \setminus A_{i-1}$ have the same dimension and are hence not adjacent. Let $U \in A_i \setminus A_{i-1}$. Then  
$$\lk_{\alpha_i}(U)\cap \alpha_{i-1} = \flag \{W \in X \mid (W \in A_0 \text{ and } U <W) \text{ or } (W \in A_{i-1} \text{ and } W <U)\}. $$ 
We want to show that each subspace in $\{W \in X \mid (W \in A_0 \text{ and } U <W) \text{ or } (W \in A_{i-1} \text{ and } W <U)\}$ is either contained in or contains $U + \ell$. Note that we have $U\cap \ell^\perp + \ell = U +\ell \in A_0 \subseteq A_{i-1}$, since we know that $U \leq \El^\perp$ because $U\in X^2$, $U + \El \in Y_0$ and $\El \leq \El^\perp$ (because $\El$ is totally isotropic) and thus $U+ \El$ cannot lie in $X^1 \cup X^3$ and hence has to lie in $X^2$.
If $W \in A_0$ and $U<W$ then $\ell \leq W$ and hence $U + \ell \leq W$. If $W \in A_{i-1}$ and $W<U$ then $W< U+ \ell$ as needed. 
The second part of the filtration is defined as follows. We set
\begin{align*}
B_i&:= A_n \cup \{U \in X^3 \mid \dim U \geq n-i+1) \},\\
\beta_i &:= \flag B_i, 0 \leq i \leq n.
\end{align*}
Then $A_n =B_0$. Let $U \in B_i \setminus B_{i-1}$. We want to show that $\lk_{\beta_i}(U) \cap \beta_{i-1}$ can be written as the join of $\{U \cap \ell ^\perp \}$ and $\flag \{W \in B_{i-1} \mid W\neq U \cap \ell^\perp, W < U \text{ or } U <W\}$.
First, note that $U \cap \ell^\perp \in Y_0$ and hence $U \cap \ell^\perp \in A_n \subseteq B_{i-1}$ (since $U\cap \El^\perp \leq \El^\perp$ and thus $U \cap \El^\perp \notin X^3$). Let $W \in B_{i-1}\setminus \{U \cap \ell^\perp\}$. If $W<U$ then $\dim W < n-i$, hence $W \in A_n\setminus A_0$, in particular $W \leq \ell^\perp$. Hence $W = W \cap \ell^\perp \leq U \cap \ell^\perp$.
If $U<W$ then $U \cap \ell^\perp \leq U \leq W$.   
Thus we get a filtration of $\kappa_0$ of length $2n$. By applying \Cref{adding vertices to cone thm} at each step, we get an $(n-2)$-cone function $\Cone_{\kappa_0}$ with 
$$\Rad_j(\Cone_{\kappa_0}) \leq 2n+1, -1 \leq j \leq n-2.$$
  
\end{proof}

We collect some of the relevant facts from \cite{Abr} to show that the filtration above has the desired properties. 

\begin{fact}[{\cite[Step 4, Step 5, Proposition 14]{Abr}}] \label{fact: D3 first filtration}
	Let $U \in Y_i \setminus Y_{i-1}$. We have that $U \not \leq \ell^\perp$ and $\dim U = i \leq n-2$. 
We get $Y_{i-1}^{<U}=Z^{<U}$. Since $\ell \not \leq U$ (otherwise $U \leq U^\perp \leq \ell^\perp$), we get that $\flag Z^{<U} \in \mathbf{C}_2$. 

On the other hand, $Y_{i-1}^{>U} = Y_0^{>U}$ and we have $\flag Y_0^{>U} \in \mathbf{C}_3$.
\end{fact}

\begin{fact}[{\cite[Step 20 (b), Proposition 14]{Abr}}] \label{fact: D3 last}
	Let $n \leq i \leq 2n-1$, $U \in Y_i \setminus Y_{i-1}$ with $\dim U < n$. Then $Y^{>U} \in \mathbf{C}_3$. 
\end{fact}

\begin{fact}[{\cite[Step 2, Proposition 14]{Abr}}] \label{fact: D3 graph case}
	Let $\kappa=\flag Y_{\EE}(V) \in \mathbf{C}_3$ in the case where $n=2$, i.e $\dim \kappa =1$. Fix a 1-dimensional subspace $\ell \in Y_{\EE}(V) =:Y$ (see Fact \ref{fact: D3 line}). Then each vertex $U \in Y$ can be connected to $\ell$ by a path of length at most 5.
\end{fact}

We now have all the necessary ingredients to show that the relevant class of simplicial complexes satisfies the conditions of \Cref{thm: the big induction for the D_n case}. The proof will be similar to the previous proofs of similar results. 

\begin{theorem} \label{thm: class for D_n case}
The class $\mathbf{C}_D= \mathbf{C}_1 \cup \mathbf{C}_2 \cup \mathbf{C}_3$ satisfies the conditions of \Cref{thm: the big induction for the D_n case}.
\end{theorem}
\begin{proof}
Let $\kappa = \flag Y \in \mathbf{C}_D$. If $\kappa \in \mathbf{C}_1 \cup \mathbf{C}_2$ then the result follows from \Cref{lemma: class 1 D_n case} and \Cref{lemma: class 2 D_n case}.

Now assume $\kappa=\flag Y_{\EE}(V) \in \mathbf{C}_3$. In this case the proof follows closely the proof of \cite[Proposition 14]{Abr}.
First of all, we fix a 1-dimensional subspace $\ell \in Y_{\EE}(V) =:Y$ (see Fact \ref{fact: D3 line}).
Let $n=\frac{1}{2}\dim V = \dim \kappa +1$. We will again prove the existence of a filtration by induction on $n$, but this time starting with $n=2$. 

Let $n=2$, hence $\dim \kappa = 1$. By Fact \ref{fact: D3 graph case} each vertex $U \in Y$ can be connected to $\ell$ by a path of length at most 5. Thus we define a $0$-cone function of $\kappa$ as follows:
\begin{align*}
\Cone_\kappa(\emptyset) &= \II_\ell \\
\Cone_\kappa(\II_{[U]})&= \sum_{i=0}^{k-1} \II_{[U_i,U_{i+1}]},
\end{align*}
where $U=U_0,\dots, U_k=\ell$ is the path from $U$ to $\ell$ for $U \in Y \setminus \{\ell\}$. In particular, $k\leq 5$. 
By \Cref{exmpl: cone for graph}, this is indeed a cone function. 
We can read off the cone radius from the explicit description and get
$$\Rad_{-1}(\Cone_\kappa)=1, \Rad_0(\Cone_\kappa) \leq 5.$$

Now assume $n\geq 3$. Consider the filtration of $\kappa$ described above. By \Cref{lemma: Y_0 for class 3 D_n case} we know that there exists an $(n-2)$-cone function $\Cone_{\kappa_0}$ with $\Rad_j(\Cone_{\kappa_0}) \leq 2n+1 = f(n)$ for all $-1 \leq j\leq n-2$. 

We now check that the rest of the filtration has the required form. Start with the first half, i.e. $1\leq i \leq n-2$. Let $U \in Y_i \setminus Y_{i-1}$. Fact \ref{fact: D3 first filtration} implies that the link of $U$ has the desired form.

Next we consider $i=n-1$. Let $U \in Y_{n-1} \setminus Y_{n-2}$. Then $Y_{n-2}^{<U} = Z^{<U}$ and $\flag Z^{<U} \in \mathbf{C}_2$. 
Furthermore, we have $Y_{n-2}^{>U}=Y_0^{>U} \neq \emptyset$. Note that $\dim \flag Y_0^{>U} = 0$ since $Y_0^{>U}$ only contains subspaces of dimension $n$. Hence it trivially has a $(-1)$-cone function with cone radius 1.

We now consider the dimension increasing filtration, i.e. $n \leq i \leq 2n-1$.
Let $U \in Y_i \setminus Y_{i-1}$ and set $k = 2n-i = \dim U$. We distinguish the cases $k=n$ and $k<n$. 

If $k=n$, then $Y_{n-1}^{>U } = \emptyset$ since there is no totally isotropic subspace of dimension $>n$. This is fine for our sake, since $\emptyset * A=A$ for an arbitrary simplicial complex $A$.
On the other hand, $Y_{n-1}^{<U}$ has again a filtration starting from an element in $\mathbf{C}_1$. To describe the filtration, we set $\EE':= (\EE \cap U) \cup ((\EE + \ell)\cap U), \mathcal{F}':= \FF \cap U, \HH :=\HH_1 := \dots := \HH_{n-2} = \widehat{\FF_n}$ and $\HH_{n-1} := \{F \in \widehat{\FF}_n \mid U \cap F =0\}$. Furthermore, we set 
$$S:= X_{\EE';\HH}(U;V) \cap \{0 < W < U \mid F' \not \leq W \text{ for all } 0 \neq F' \in \FF' \}.$$
By Step (20) in \cite[Proposition 14]{Abr}, we know that $S \leq Y_{n-1}^{<U}$ and given a cone function of $\flag S$ we get a cone function of $\flag Y_{n-1}^{<U}$ by one application of \Cref{adding vertices to cone thm}.
To see that $\flag S$ has a cone function, we note that if $W \pitchfork_U F'$ then this implies that $F' \not \leq W$. Thus we set $\EE'' = \EE' \cup \FF'$ and get 
$X_{\EE'',\HH}(U;V) \subseteq S $.
We set 
$$A_0 := X_{\EE'',\HH}(U;V), A_i = \{W \in S \mid W \in A_0 \text{ or } \dim W \leq i\}, 1 \leq i \leq n-1.$$
We follow the by now standard procedure to check that this is a valid filtration. 
Let $W \in A_i \setminus A_{i-1}$. Then $A_{i-1}^{>W} = A_0^{>W} =  X_{\EE'',\HH}(U;V)^{>W}$ and by Step 4. in \cite[Lemma 33]{Abr} (see also Fact \ref{fact: D1 first part of filtration}) we have $\flag X_{\EE'',\HH}(U;V)^{>W} \in \mathbf{C}_1$ .
On the other hand, $A_{i-1}^{<W} = S^{<W}$. Note that if $W \in S$ and $W' \in X_{\EE';\HH}(U;V)$ such that $W' \leq W$ then $W' \in S$ since if $F' \not \leq W$ then $F'\not \leq W' \leq W$. Thus
$$S^{<W} = X_{\EE';\HH}(U;V)^{<W} = \{0<W'<W \mid W' \cap \EE', W' \mathbin{@} \HH_j \text{ for } \dim W' = j \} = X_{\EE' \cap W, \HH}(W;V)$$
and $\flag X_{\EE' \cap W, \HH}(W;V) \in \mathbf{C}_1$.
Hence we got a filtration 
$$\flag A_0 \subseteq \dots \flag A_{n-1} = \flag S \subseteq \flag Y_{n-1}^{<U}$$
of length $n=\dim U$ which satisfies the assumptions.

If $\dim U <n$, then $Y_{i-1}^{<U} = Z^{<U}$ and we again have $\flag Z^{<U} \in \mathbf{C}_3$. Additionally, we have $Y_{i-1}^{>U} = Y^{>U}$ for which by Fact \ref{fact: D3 last} we have $\flag Y^{>U} \in \mathbf{C}_3$.

Hence all the conditions of \Cref{thm: the big induction for the D_n case} are satisfied. 
\end{proof}
The following remark from \cite[page 105]{Abr} gives the connection between $X_\EE(V)$ and $Y_\EE(V)$. 
\begin{remark} \label{rmk: connection to op compl D_n}
Let $\tau = \{E_1, \dots E_r \}$ be a simplex of $\tilde{\Delta}= \orifl \tilde{X}$ and set $\EE = \EE(\tau)=\{E_i,E_i^\perp \mid 1 \leq i \leq r\}$. Then $\EE \cap \mathcal{U}_n = \emptyset$, $Y_{\EE}(V) = X_{\EE}(V)$, $e_n \leq 2, e_{n-1} = e_{n+1}=0$ and $e_i \leq 1$ for all other $1 \leq i \leq 2n-1$.
\end{remark}

\begin{corollary}
Let $\tilde{\Delta}$ be a building of type $D_n$ over the field $K$, and $a \in \tilde{\Delta}$ be a simplex. 
If $\lvert K \rvert \geq 2^{2n-1}$, then $\tilde{\Delta}^0(a)$ has an $(n-2)$-cone function with
$\Rad_j(\Cone_{\tilde{\Delta}^0(a)}) \leq 2 \mathcal{R}(n-1), -1 \leq j \leq n-2$, where $\mathcal{R}(n)$ does not depend on $K$.
\end{corollary}
\begin{proof}
By \Cref{thm: class for D_n case}, we have that $Y_{\EE(a)}(V)$ has a cone function with radius $\leq \mathcal{R}(n)$, where $n=\dim Y_{\EE(a)}(V)$. By \Cref{rmk: connection to op compl D_n}, we have $X_{\EE(a)}(V) = Y_{\EE(a)}(V)$. Thus \Cref{prop: quantitative subdivision D_n} yields a cone function for $\tilde{\Delta}^0(a) = \tilde{T}_{\EE(a)}(V)$ with cone radius bounded by $2\mathcal{R}(n)$.
\end{proof}

\bibliographystyle{alpha}
\bibliography{bibl}
\end{document}